\documentclass{amsart}
\usepackage[foot]{amsaddr}
\usepackage{amssymb,latexsym,amsmath,amsfonts,graphicx,xcolor,pb-diagram}
\usepackage[mathscr]{eucal}
\usepackage{tikz, tikz-cd}
\usepackage{relsize}
\usepackage{slashed}
\usepackage{pgfplots}
\pgfplotsset{compat=1.18}
\usepgfplotslibrary{fillbetween}
\newcommand{\mpar}{\par \medskip \par }

\numberwithin{equation}{section}
\swapnumbers
\theoremstyle{definition}

\newtheorem{theorem}{Theorem}[section]

\newtheorem{remark}[theorem]{Remark}

\makeindex

\begin{document}

\title{Build Your Own Neural Network From Scratch}

\author{Samantha Donner$^1$, Ethan Geary$^2$, Peter March$^3$, Mathew Umano$^1$}

\begin{abstract}
These notes are based on our reading of Michael Nielsen's online book \textit{neuralnetworksanddeeplearning.com}. They represent our effort to elaborate his exposition of neural networks and to provide supporting mathematical detail. Our goal is to provide self-contained instructions for building a basic neural network in the form of pseudocode, so that an interested reader can build a network from scratch and know why it works. No prior knowledge of neural networks or of Python is assumed. Some familiarity with  linear algebra,  calculus, and probability is helpful but we explain the relevant mathematics in a set of appendices designed to be accessible to newcomers.
\end{abstract} 

\address{$^1$Department of Environmental Sciences, Rutgers University - New Brunswick} 
\address{$^2$Metuchan High School, Metuchan NJ} 
\address{$^3$Department of Mathematics, Rutgers University - New Brunswick}

\maketitle
\section{Introduction}
These notes were written as a byproduct of a high school summer research experience reading Michael Nielsen's online book \textit{neuralnetworksanddeeplearning.com}. \cite{N} Both of the initial team members, Geary and March, were newcomers to the field of deep learning and neither had any functional knowledge of Python, the language of choice in Nielsen's book. To fill the gap they recruited undergraduate researchers Donner and Umano to join the team.

\mpar
The original impetus for these notes was to elaborate various ideas we encountered while reading the book and to backfill or supplement the mathematical content. To give the notes narrative shape, we decided to address them to an audience of high school seniors who, we assumed, were not necessarily familiar with any of the basic ingredients: neural networks; linear algebra, calculus and probability; and Python. The goal was to begin at the beginning and end with understandable pseudocode that implements a neural network on a computer, while elaborating key mathematical ideas along the way.

\mpar
It was not part of the goal to provide a detailed account of neural networks and deep learning or to provide full proofs of the necessary mathematics. Rather, the goal was to create an opening for interested beginners to engage with the subject by providing enough technical background so they can see why the pseudocode actually works, rather than just accepting it on faith as a black box.

\mpar
Nielsen's work was published in 2015 and there has been an enormous amount of progress since then. In particular, the concept of \textit{transformers} \cite{AIAYN} introduced in 2017, revolutionized the field. That does not mean that the network architecture we will discuss is irrelevant. To the contrary, all the core ideas we will discuss in Nielsen's setting are present in current frontier models and continue to play a critical role in neural networks and artificial intelligence.

\mpar
The notes are divided into four sections and three appendices. Section 2 describes in qualitative terms what a neural network is and how it works. Here qualitative means without using mathematics in a significant way. The idea is to provide a mental picture of a neural network as a mechanism for processing information. This picture we paint is not accurate in all respects but we think is a useful caricature of real neural networks.

\mpar
Section 3 shows how mathematics converts qualitative ideas about neural networks into quantitive ones suitable for computing. It is organized into subsections that address important aspects of neural networks. Each subsection begins with a statement about what kind of mathematics will be needed and directs the reader to appropriate parts of the appendices, if the mathematics is unfamiliar. 

\mpar
 Section 4  presents pseudocode that serves as a guide to writing actual code, say in Python, that implements a functioning neural network. It is quite detailed and even repetitive in places, unlike typical pseudocode. Our purpose isn't efficiency but rather pedagogy. Rendering correct mathematics into correct code isn't as straightforward as one might think so we've erred on the side of repetition and attention to detail at the expense of concision.

\mpar
Appendices A, B, and C concern linear algebra, calculus, and probability, respectively. The style is informal with a strong emphasis on understanding basic ideas. Necessarily, the mathematics rapidly gets sophisticated as we move, for example, from vectors in the plane to vectors in very high dimensional spaces; from the meaning of the derivative as the slope of the tangent line to the meaning of the gradient as the direction of steepest ascent; and from flipping coins to the law of large numbers and the central limit theorem. Thus, the appendices start slowly so as to engage newcomers but build momentum quickly in order to reach the mathematical results needed to describe how neural networks learn.

\mpar
These notes do not cover all of Nielsen's book. We stop short of discussing how to train the network in practice, once the network's been built, or using a trained network to accomplish a task, like recognizing handwritten digits. Knowing how to train a network effectively is essential to the success of deep learning. But training is as much of an art as it is a science, so we content ourselves with showing how to build a network from scratch and hope that will inspire interested readers to pursue the subject further and perfect their art.

\mpar
\section*{Acknowledgements}
Our thanks are due to Moulik Kallupalam Balasubram, Midge Cozzens, Matthew March, David Pennock, and Larry Susanka for valuable feedback that helped improve the text. Any errors, misconceptions, and infelicities that remain are solely the authors' responsibility.

\mpar\mpar
\section {What Are Neural Networks?} 

We're going to take this in three steps. First, we'll describe how neural networks mimic what we know about the how the brain processes information via cascades of stimulation and response of neurons. We describe the so-called \textit{architecture} of the neural network as \textit{input} and \textit{output} layers connected by intermediate or \textit{hidden} layers of artificial neurons. This leads to a qualitative picture of neural networks that is correct in spirit but naive in practice. Second, we start to quantify these ideas by thinking of neurons as containing numbers, called \textit{activations}, and connections between neurons as defining transformations of activations. Third, we discuss how a neural network learns from data, which is called \textit{training} the network. The main idea is to choose a measure of how well the network performs and use that measure to create a positive feedback loop that improves performance.
 
\mpar
\subsection{Simplified Picture of Neural Networks} 

\textit{The main things to know and remember from this section are (1) neural networks are collections of \textit{artifical neurons} arranged in a sequence layers, (2) the \textit{input layer} is on the left, the \textit{output layer} is on the right and they are connected by one or more intermediate \textit{hidden layers}, (2) some neurons in a given layer are connected to some neurons in the next layer to the right, (3) the number of layers, the number of neurons in each layer and the pattern of connections between layers is called the \textit{architecture} of the neural network, (4) a \textit{feedforward} network operates from left to right on data entered into the input layer and produces data in the output layer, and (5) the role of hidden layers is to detect patterns in the input data.}

\mpar\mpar
To gain a basic understanding of neural networks we must first understand something about the human brain. Imagine you are looking at a picture of a dog. How do we perceive it as a dog? What does this correspond to in the brain? 

\mpar
The visual area of the human contains five cortices each consisting of millions of neurons that process information within the picture. When the brain perceives a certain feature in the image, the corresponding neurons fire and stimulate other neurons, confirming this feature is present. Roughly speaking, the cascade of stimulation from neurons in the retina to neurons in the visual cortex is what processes visual information in our brains. Neural networks work in the same manner. They are networks of artificial neurons programmed into computers to process data in a way similar to the brain.

\mpar
To explain further what neural networks are, let's create an imaginary one that can receive an image of an animal as input and infer the image to be either a dog, bear, or a human. How would that work? Let's imagine that each artificial neuron represents a specific feature of animals - ranging, say,  from gills to whiskers to antlers. If the neurons representing paws, fur, four legs, and a tail all fire, then the  network may infer it's looking at a dog. This isn't an accurate picture, to be honest, because both cortical units in the brain and artificial neurons in a computer respond to much finer features in the visual field than anatomical parts. But it's a useful fiction that helps us get the main ideas across, so we'll stick with it for the time being.

\mpar If this imaginary network going to work at all, it must be able to discern combinations of features in a picture of an animal that distinguish dogs, bears, and humans from one another. For example, an animal with four legs is not likely a human so perhaps the network can recognize that having four legs is a key feature that can be used to distinguish humans from dogs and bears. But it needs additional features other than the number of legs to distinguish a dog from a bear. If networks could also detect that each animal has something different attached to their limbs, for example dogs have paws, bears have claws, and humans have fingers, then hands would be a very helpful distinguishing feature. If networks could also detect hair and distinguish fur on the body from hair only on the head, that would also be helpful.

\mpar
On the other hand, there can be ambiguities in the images that cause problems. At first sight it seems like detecting the kind of hand will always lead to the right answer. But what if the human is clasping her hands behind her back or the bear is hibernating, all curled up in a ball? We want the network to recognize that a crawling baby is human and a freshly shaved poodle with a fluffy head of hair is not human. When you think about these kinds of ambiguities you realize you can't expect the network to classify pictures with perfect accuracy. 

\mpar
However, we can hope to do better than pure chance by having lots and lots of examples to learn from. A neural network is provided with information through \textit{training data}, from which the network develops the ability to recognize certain features of animals and draw  inferences from them. Training data are digital images or other forms of data that are already labeled with  answers. In our example, training data are images of animals plus a label that the image is of a dog, or a bear, or a human. With enough training data the hope is the neural network can learn to correctly classify animals in images. 

\mpar
The neural network receives data from images in an input layer of neurons, and this data is then sent to hidden layers of neurons. The hidden layers contain neurons that can detect certain features like hairiness, leg number, and hand differences. For each detected feature, its corresponding neurons fire and send information to an output layer which provides the final result: dog, bear, or human. This is how a \textit{feedforward neural network} works.

\mpar
How would we judge whether the network is any good at sorting images of dogs, bears and, humans? An obvious thing to do is feed the network new images it hasn't seen before and see how well it does. These new images are called the \textit{evaluation data}. If the network correctly sorts the evaluation data about a third of the time then its performance is no better than blind guessing. On the other hand, if the network gets it right, say  9 times out of 10, then that's pretty persuasive evidence that the network has learned to sort images of animals based on the training data.

\begin{figure}[h]
\begin{tikzpicture}
\draw (0,0) circle (0.25);
\draw (0,1) circle (0.25);
\draw (0,2) circle (0.25);

\draw[thick] (0.25, 0) -- (1.75,0);
\draw[thick] (0.25, 0) -- (1.75,-1);
\draw[thick] (0.25, 0) -- (1.75,-2);
\draw[thick] (0.25, 0) -- (1.75,1);
\draw[thick] (0.25, 0) -- (1.75,2);
\draw[thick] (0.25, 0) -- (1.75,3);
\draw[thick] (0.25, 0) -- (1.75,4);

\draw[thick] (0.25, 1) -- (1.75,0);
\draw[thick] (0.25, 1) -- (1.75,-1);
\draw[thick] (0.25, 1) -- (1.75,-2);
\draw[thick] (0.25, 1) -- (1.75,1);
\draw[thick] (0.25, 1) -- (1.75,2);
\draw[thick] (0.25, 1) -- (1.75,3);
\draw[thick] (0.25, 1) -- (1.75,4);

\draw[thick] (0.25, 2) -- (1.75,0);
\draw[thick] (0.25, 2) -- (1.75,-1);
\draw[thick] (0.25, 2) -- (1.75,-2);
\draw[thick] (0.25, 2) -- (1.75,1);
\draw[thick] (0.25, 2) -- (1.75,2);
\draw[thick] (0.25, 2) -- (1.75,3);
\draw[thick] (0.25, 2) -- (1.75,4);

\draw (2,-2) circle (0.25);
\draw (2,-1) circle (0.25);
\draw (2,0) circle (0.25);
\draw (2,1) circle (0.25);
\draw (2,2) circle (0.25);
\draw (2,3) circle (0.25);
\draw (2,4) circle (0.25);

\draw (4,0) circle (0.25);
\draw (4,1) circle (0.25);
\draw (4,2) circle (0.25);

\draw[thick] (2.25, -2) -- (3.75,0);
\draw[thick] (2.25, -1) -- (3.75, 0);
\draw[thick] (2.25, 0) -- (3.75,0);
\draw[thick] (2.25, 1) -- (3.75,0);
\draw[thick] (2.25, 2) -- (3.75,0);
\draw[thick] (2.25, 3) -- (3.75,0);
\draw[thick] (2.25, 4) -- (3.75,0);

\draw[thick] (2.25, -2) -- (3.75,1);
\draw[thick] (2.25, -1) -- (3.75, 1);
\draw[thick] (2.25, 0) -- (3.75,1);
\draw[thick] (2.25, 1) -- (3.75,1);
\draw[thick] (2.25, 2) -- (3.75,1);
\draw[thick] (2.25, 3) -- (3.75,1);
\draw[thick] (2.25, 4) -- (3.75,1);

\draw[thick] (2.25, -2) -- (3.75,2);
\draw[thick] (2.25, -1) -- (3.75, 2);
\draw[thick] (2.25, 0) -- (3.75,2);
\draw[thick] (2.25, 1) -- (3.75,2);
\draw[thick] (2.25, 2) -- (3.75,2);
\draw[thick] (2.25, 3) -- (3.75,2);
\draw[thick] (2.25, 4) -- (3.75,2);

\node at (-1.5,1) {\text{input layer}};
\node at (2,4.6) {\text{hidden layer}};
\node at (5.5,1) {\text{output layer}};

\end{tikzpicture}
\caption{An example of a neural network. It has three nodes in the input layer and three nodes in the output layer. It has just one hidden layer consisting of seven nodes. Our example is a \textit{shallow} network because it only has one hidden layer. In contrast, \textit{deep networks} have multiple hidden layers.}
\end{figure}
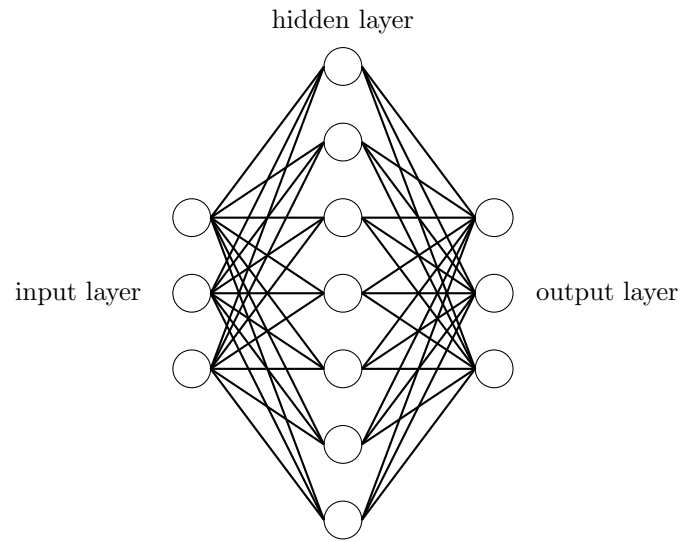

\mpar\mpar
You might have just read that and been completely lost. What is an input layer, hidden layer, and an output layer? What even is a layer? And what does feedforward mean? Well, to better explain it, let's draw a picture - a visual representation of a neural network. See Figure 1. We imagine connections between nodes go from left to right - from input to hidden to output. A given node may be connected to some nodes to the left and connected to some other nodes on the right, but in this example every node in the input and output layers is connected every node in the hidden layer. 

\mpar
Real life industrial scale neural networks can have millions of neurons arranged in multiple hidden layers. For the sake of illustration our network has just 13 neurons and one hidden layer.

\mpar
\textbf{Caveat Lector.} \textit{Before going further with our imaginary network, we reiterate that we're making some compromises to get the main ideas across. So reader beware: real neural networks don't work exactly like our imaginary network! It's really a caricature - or maybe even a cartoon of a caricature - but it helps us get across some ideas about the flow of information in a neural network at the risk, perhaps, of creating a false impression. To sharpen the point, networks do not have the a priori ability to detect designated features, as we assume below. Rather they have the ability to find patterns in data that humans subsequently perceive as features, without being told in advance what to look for. We temporarily set this distinction aside for narrative purposes.}
\mpar
To get started, let's assume the network can detect variants of \textit{Hands, Legs, Hair} and label the input neurons accordingly as in Figure 2.

\mpar
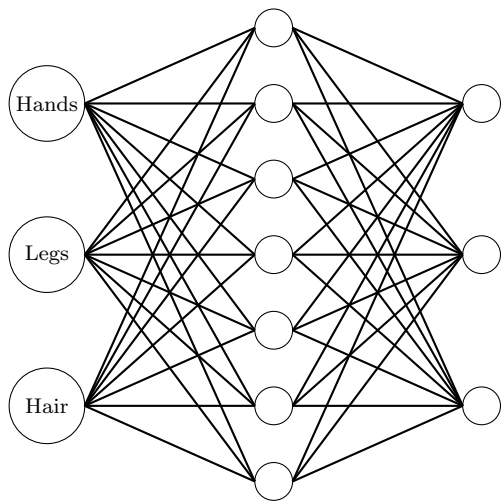
\begin{figure} [h]
\begin{tikzpicture}
\draw (0,0) circle (0.5) node{\footnotesize{Hair}};
\draw (0,2) circle (0.5) node{\footnotesize{Legs}};
\draw (0,4) circle (0.5) node{\footnotesize{Hands}};

\draw (3,-1) circle (0.25);
\draw (3,0) circle (0.25);
\draw (3,1) circle (0.25);
\draw (3,2) circle (0.25);
\draw (3,3) circle (0.25);
\draw (3,4) circle (0.25);
\draw (3,5) circle (0.25);

\draw (5.75,0) circle (0.25);
\draw (5.75,2) circle (0.25); 
\draw (5.75,4) circle (0.25);

\draw[thick] (3.25, -1) -- (5.5,0);
\draw[thick] (3.25, 0) -- (5.5, 0);
\draw[thick] (3.25, 1) -- (5.5,0);
\draw[thick] (3.25, 2) -- (5.5,0);
\draw[thick] (3.25, 3) -- (5.5,0);
\draw[thick] (3.25, 4) -- (5.5,0);
\draw[thick] (3.25, 5) -- (5.5,0);

\draw[thick] (3.25, -1) -- (5.5,2);
\draw[thick] (3.25, 0) -- (5.5, 2);
\draw[thick] (3.25, 1) -- (5.5,2);
\draw[thick] (3.25, 2) -- (5.5,2);
\draw[thick] (3.25, 3) -- (5.5,2);
\draw[thick] (3.25, 4) -- (5.5,2);
\draw[thick] (3.25, 5) -- (5.5,2);

\draw[thick] (3.25, -1) -- (5.5,4);
\draw[thick] (3.25, 0) -- (5.5, 4);
\draw[thick] (3.25, 1) -- (5.5,4);
\draw[thick] (3.25, 2) -- (5.5,4);
\draw[thick] (3.25, 3) -- (5.5,4);
\draw[thick] (3.25, 4) -- (5.5,4);
\draw[thick] (3.25, 5) -- (5.5,4);

\draw[thick] (2.75, -1) -- (0.5,0);
\draw[thick] (2.75, 0) -- (0.5, 0);
\draw[thick] (2.75, 1) -- (0.5,0);
\draw[thick] (2.75, 2) -- (0.5,0);
\draw[thick] (2.75, 3) -- (0.5,0);
\draw[thick] (2.75, 4) -- (0.5,0);
\draw[thick] (2.75, 5) -- (0.5,0);

\draw[thick] (2.75, -1) -- (0.5,2);
\draw[thick] (2.75, 0) -- (0.5, 2);
\draw[thick] (2.75, 1) -- (0.5,2);
\draw[thick] (2.75, 2) -- (0.5,2);
\draw[thick] (2.75, 3) -- (0.5,2);
\draw[thick] (2.75, 4) -- (0.5,2);
\draw[thick] (2.75, 5) -- (0.5,2);

\draw[thick] (2.75, -1) -- (0.5,4);
\draw[thick] (2.75, 0) -- (0.5, 4);
\draw[thick] (2.75, 1) -- (0.5,4);
\draw[thick] (2.75, 2) -- (0.5,4);
\draw[thick] (2.75, 3) -- (0.5,4);
\draw[thick] (2.75, 4) -- (0.5,4);
\draw[thick] (2.75, 5) -- (0.5,4);
\end{tikzpicture}
\caption{Input layer neurons enlarged to accommodate labels.}
\end{figure}

Next, we assume variants of Hair are \textit{Finger, Claw, Paw}; the variants of Legs are \textit{2 Legs, 4 Legs}; and the variants of Hair are \textit{Fur, Head}. 
We imagine these label the neurons of the hidden layer. See Figure 3.

\begin{figure}[h]
\begin{tikzpicture}
\draw (0.25,0) circle (0.25);  
\draw (0.25,2.5) circle (0.25);  
\draw (0.25,5) circle (0.25);  

\draw (3,-1.25) circle (0.5) node{\footnotesize{4Legs}};
\draw (3,0) circle (0.5) node{\footnotesize{2Legs}};
\draw (3,1.25) circle (0.5) node{\footnotesize{Head}};
\draw (3,2.5) circle (0.5) node{\footnotesize{Fur}};
\draw (3,3.75) circle (0.5) node{\footnotesize{Paw}};
\draw (3,5) circle (0.5) node{\footnotesize{Claw}};
\draw (3,6.25) circle (0.5) node{\footnotesize{Finger}};

\draw (5.75,0) circle (0.25);
\draw (5.75,2.5) circle (0.25); 
\draw (5.75,5) circle (0.25);

\draw[thick] (3.5, -1.25) -- (5.5,0);
\draw[thick] (3.5, 0) -- (5.5, 0);
\draw[thick] (3.5, 1.25) -- (5.5,0);
\draw[thick] (3.5, 2.5) -- (5.5,0);
\draw[thick] (3.5, 3.75) -- (5.5,0);
\draw[thick] (3.5, 5) -- (5.5,0);
\draw[thick] (3.5, 6.25) -- (5.5,0);

\draw[thick] (3.5, -1.25) -- (5.5,2.5);
\draw[thick] (3.5, 0) -- (5.5, 2.5);
\draw[thick] (3.5, 1.25) -- (5.5,2.5);
\draw[thick] (3.5, 2.5) -- (5.5,2.5);
\draw[thick] (3.5, 3.75) -- (5.5,2.5);
\draw[thick] (3.5, 5) -- (5.5,2.5);
\draw[thick] (3.5, 6.25) -- (5.5,2.5);

\draw[thick] (3.5, -1.25) -- (5.5,5);
\draw[thick] (3.5, 0) -- (5.5, 5);
\draw[thick] (3.5, 1.25) -- (5.5,5);
\draw[thick] (3.5, 2.5) -- (5.5,5);
\draw[thick] (3.5, 3.75) -- (5.5,5);
\draw[thick] (3.5, 5) -- (5.5,5);
\draw[thick] (3.5, 6.25) -- (5.5,5);

\draw[thick] (2.5, -1.25) -- (0.5,0);
\draw[thick] (2.5, 0) -- (0.5, 0);
\draw[thick] (2.5, 1.25) -- (0.5,0);
\draw[thick] (2.5, 2.5) -- (0.5,0);
\draw[thick] (2.5, 3.75) -- (0.5,0);
\draw[thick] (2.5, 5) -- (0.5,0);
\draw[thick] (2.5, 6.25) -- (0.5,0);

\draw[thick] (2.5, -1.25) -- (0.5,2.5);
\draw[thick] (2.5, 0) -- (0.5, 2.5);
\draw[thick] (2.5, 1.25) -- (0.5,2.5);
\draw[thick] (2.5, 2.5) -- (0.5,2.5);
\draw[thick] (2.5, 3.75) -- (0.5,2.5);
\draw[thick] (2.5, 5) -- (0.5,2.5);
\draw[thick] (2.5, 6.25) -- (0.5,2.5);

\draw[thick] (2.5, -1.25) -- (0.5,5);
\draw[thick] (2.5, 0) -- (0.5, 5);
\draw[thick] (2.5, 1.25) -- (0.5,5);
\draw[thick] (2.5, 2.5) -- (0.5,5);
\draw[thick] (2.5, 3.75) -- (0.5,5);
\draw[thick] (2.5, 5) -- (0.5,5);
\draw[thick] (2.5, 6.25) -- (0.5,5);
\end{tikzpicture}
\caption{Hidden layer neurons enlarged to accommodate labels.}
\end{figure}
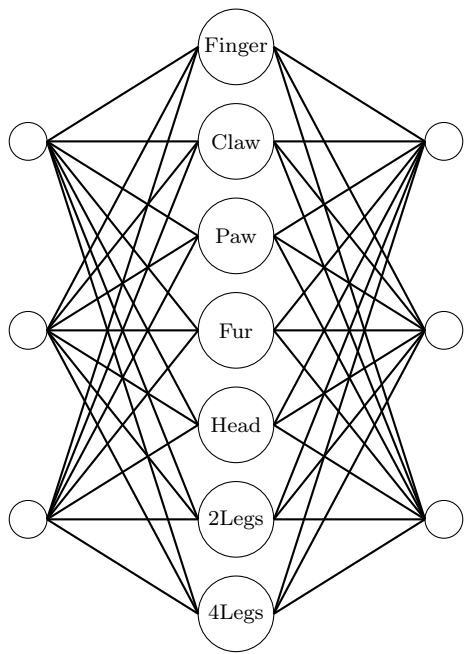

\mpar
Finally, we have our output layer, the simplest to understand. There are three possible outcomes in our imaginary network, \textit{Dog, Bear} or \textit{Human}. See Figure 4.

\begin{figure}[h]
\begin{tikzpicture}
\draw (0.25,0) circle (0.25);  
\draw (0.25,2) circle (0.25);  
\draw (0.25,4) circle (0.25);  

\draw (3,-1) circle (0.25);
\draw (3,0) circle (0.25);
\draw (3,1) circle (0.25);
\draw (3,2) circle (0.25);
\draw (3,3) circle (0.25);
\draw (3,4) circle (0.25);
\draw (3,5) circle (0.25);

\draw (6,0) circle (0.5) node{\footnotesize{Human}};
\draw (6,2) circle (0.5) node{\footnotesize{Bear}};
\draw (6,4) circle (0.5) node{\footnotesize{Dog}};

\draw[thick] (3.25, -1) -- (5.5,0);
\draw[thick] (3.25, 0) -- (5.5, 0);
\draw[thick] (3.25, 1) -- (5.5,0);
\draw[thick] (3.25, 2) -- (5.5,0);
\draw[thick] (3.25, 3) -- (5.5,0);
\draw[thick] (3.25, 4) -- (5.5,0);
\draw[thick] (3.25, 5) -- (5.5,0);

\draw[thick] (3.25, -1) -- (5.5,2);
\draw[thick] (3.25, 0) -- (5.5, 2);
\draw[thick] (3.25, 1) -- (5.5,2);
\draw[thick] (3.25, 2) -- (5.5,2);
\draw[thick] (3.25, 3) -- (5.5,2);
\draw[thick] (3.25, 4) -- (5.5,2);
\draw[thick] (3.25, 5) -- (5.5,2);

\draw[thick] (3.25, -1) -- (5.5,4);
\draw[thick] (3.25, 0) -- (5.5, 4);
\draw[thick] (3.25, 1) -- (5.5,4);
\draw[thick] (3.25, 2) -- (5.5,4);
\draw[thick] (3.25, 3) -- (5.5,4);
\draw[thick] (3.25, 4) -- (5.5,4);
\draw[thick] (3.25, 5) -- (5.5,4);

\draw[thick] (2.75, -1) -- (0.5,0);
\draw[thick] (2.75, 0) -- (0.5, 0);
\draw[thick] (2.75, 1) -- (0.5,0);
\draw[thick] (2.75, 2) -- (0.5,0);
\draw[thick] (2.75, 3) -- (0.5,0);
\draw[thick] (2.75, 4) -- (0.5,0);
\draw[thick] (2.75, 5) -- (0.5,0);

\draw[thick] (2.75, -1) -- (0.5,2);
\draw[thick] (2.75, 0) -- (0.5, 2);
\draw[thick] (2.75, 1) -- (0.5,2);
\draw[thick] (2.75, 2) -- (0.5,2);
\draw[thick] (2.75, 3) -- (0.5,2);
\draw[thick] (2.75, 4) -- (0.5,2);
\draw[thick] (2.75, 5) -- (0.5,2);

\draw[thick] (2.75, -1) -- (0.5,4);
\draw[thick] (2.75, 0) -- (0.5, 4);
\draw[thick] (2.75, 1) -- (0.5,4);
\draw[thick] (2.75, 2) -- (0.5,4);
\draw[thick] (2.75, 3) -- (0.5,4);
\draw[thick] (2.75, 4) -- (0.5,4);
\draw[thick] (2.75, 5) -- (0.5,4);
\end{tikzpicture}
\caption{Output layer neurons enlarged to accommodate labels.}
\end{figure}
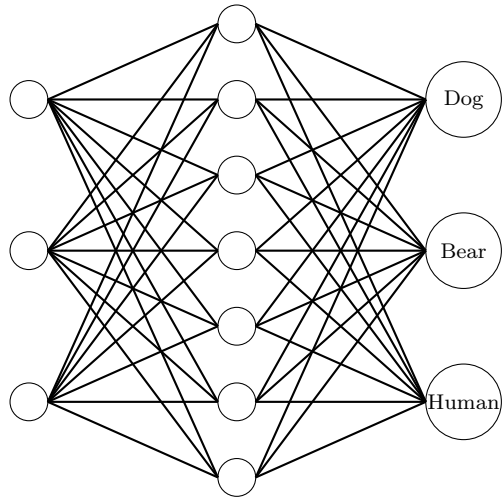

\mpar
What to the connections between neurons mean? We will use something called a \textit{perceptron} - essentially just an extremely basic neuron. Unlike general neurons, which we will delve into later, perceptrons can only send two answers to the next layer of neurons: \textit{Yes or No; Fire or Not Fire; 0 or 1.} 

\mpar
The main thing to remember is that in real neural networks, the purpose of hidden layers is to detect patterns in the input data, draw inferences or conclusions about the input based on those patterns, and register the conclusion in the output layer. We're using our cartoon network to exemplify this process and gain a mental picture of how networks process information.

\mpar
OK, our cartoon network is up and running. This corresponds to having trained the network so it is able to detect patterns and make inferences - how this is actually done will be discussed in subsequent sections.

\mpar
To see how a neural network works, let's give it a picture of an animal with \textit{Hair, Legs} and \textit{Hands}. This information appears in the input layer. Next, the network detects that \textit{Hair} is \textit{Fur} not \textit{Head}; it detects that \textit{Hand} is  \textit{Claw}, not \textit{Fingers} or \textit{Paws}; and it detects \textit{Legs} to be \textit{4Legs} not \textit{2Legs}. This causes the \text{Hair} neuron to fire in the direction of the \textit{Fur} neuron, \textit{Hand} to fire in the direction of \textit{Claw}, and \textit{Legs} to fire in the direction of \textit{4Legs}. So, the \textit{Fur, Claw} and \textit{4Legs} neurons all fire in the direction of the \textit{Bear} neuron because this is the only animal of the three that possess these three features. A picture of the path of fired neurons is in Figure 5.

\begin{figure}[h]
\begin{tikzpicture}
\draw (0,0) circle (0.5) node{\footnotesize{Hair}};  
\draw (0,2.5) circle (0.5) node{\footnotesize{Legs}}; 
\draw (0,5) circle (0.5) node{\footnotesize{Hands}}; 

\draw (3,-1.25) circle (0.5) node{\footnotesize{4Legs}};
\draw (3,0) circle (0.5) node{\footnotesize{2Legs}};
\draw (3,1.25) circle (0.5) node{\footnotesize{Head}};
\draw (3,2.5) circle (0.5) node{\footnotesize{Fur}};
\draw (3,3.75) circle (0.5) node{\footnotesize{Paw}};
\draw (3,5) circle (0.5) node{\footnotesize{Claw}};
\draw (3,6.25) circle (0.5) node{\footnotesize{Finger}};

\draw (6,0) circle (0.5) node{\footnotesize{Human}};
\draw (6,2.5) circle (0.5) node{\footnotesize{Bear}};
\draw (6,5) circle (0.5) node{\footnotesize{Dog}};

\draw[very thin] (3.5, -1.25) -- (5.5,0);
\draw[very thin] (3.5, 0) -- (5.5, 0);
\draw[very thin] (3.5, 1.25) -- (5.5,0);
\draw[very thin] (3.5, 2.5) -- (5.5,0);
\draw[very thin] (3.5, 3.75) -- (5.5,0);
\draw[very thin] (3.5, 5) -- (5.5,0);
\draw[very thin] (3.5, 6.25) -- (5.5,0);

\draw[ultra thick] (3.5, -1.25) -- (5.5,2.5);
\draw[very thin] (3.5, 0) -- (5.5, 2.5);
\draw[ultra thick] (3.5, 2.5) -- (5.5,2.5);
\draw[ultra thick] (3.5, 5) -- (5.5,2.5);

\draw[very thin] (3.5, -1.25) -- (5.5,5);
\draw[very thin] (3.5, 0) -- (5.5, 5);
\draw[very thin] (3.5, 1.25) -- (5.5,5);
\draw[very thin] (3.5, 2.5) -- (5.5,5);
\draw[very thin] (3.5, 3.75) -- (5.5,5);
\draw[very thin] (3.5, 5) -- (5.5,5);
\draw[very thin] (3.5, 6.25) -- (5.5,5);

\draw[very thin] (2.5, -1.25) -- (0.5,0);
\draw[very thin] (2.5, 0) -- (0.5, 0);
\draw[very thin] (2.5, 1.25) -- (0.5,0);
\draw[ultra thick] (2.5, 2.5) -- (0.5,0);
\draw[very thin] (2.5, 3.75) -- (0.5,0);
\draw[very thin] (2.5, 5) -- (0.5,0);
\draw[very thin] (2.5, 6.25) -- (0.5,0);

\draw[ultra thick] (2.5, -1.25) -- (0.5,2.5);
\draw[very thin] (2.5, 0) -- (0.5, 2.5);
\draw[very thin] (2.5, 1.25) -- (0.5,2.5);
\draw[very thin] (2.5, 2.5) -- (0.5,2.5);
\draw[very thin] (2.5, 3.75) -- (0.5,2.5);
\draw[very thin] (2.5, 5) -- (0.5,2.5);
\draw[very thin] (2.5, 6.25) -- (0.5,2.5);

\draw[very thin] (2.5, -1.25) -- (0.5,5);
\draw[very thin] (2.5, 0) -- (0.5, 5);
\draw[very thin] (2.5, 1.25) -- (0.5,5);
\draw[very thin] (2.5, 2.5) -- (0.5,5);
\draw[very thin] (2.5, 3.75) -- (0.5,5);
\draw[ultra thick] (2.5, 5) -- (0.5,5);
\draw[very thin] (2.5, 6.25) -- (0.5,5);
\end{tikzpicture}
\caption{Firing between neurons depicted in black; non-firing between neurons us in light grey.}
\end{figure}
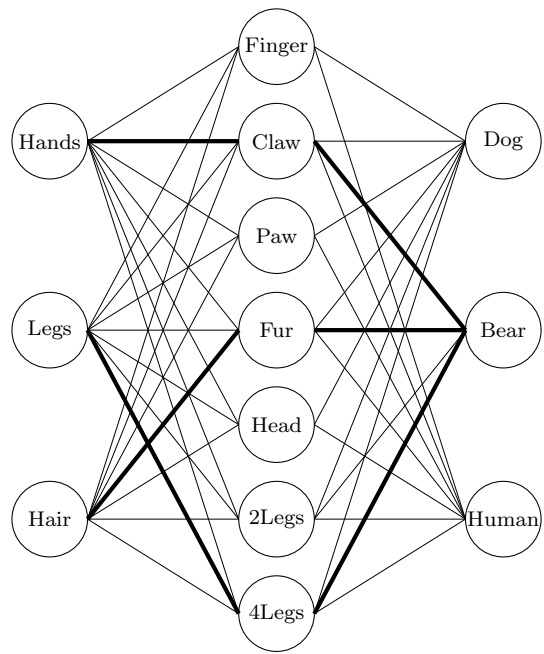

\mpar
But there is a potential problem, even in this cartoon. Imagine a distorted image is given to the system where it is hard to observe each feature. And let's say the network is having difficulty detecting if the hands are claws, paws, or fingers. For example, imagine there is a picture of a human with extremely sharp costume nails and the system cannot determine if they are fingers or claws. Having to make a decision, a network using perceptrons must make a choice to detect them as claws, say, not fingers. This could completely alter the outcome and  give us the incorrect output of \textit{Bear} rather than \textit{Human}. This instability is the flaw of perceptrons: having only the options of \textit{Yes} or \textit{No} means even one mistake can change everything drastically. Luckily, there are other kinds of neurons that don't have this flaw.

\mpar\mpar
\subsection{Sigmoidal Neurons}
\textit{The main things to know and remember from this section are (1) sigmoidal neurons contain activations, (2) arrows between neurons provide the means to modify activations of neurons in a given layer based on activations of neurons in the previous layer, (3) modifications are made using weights, biases and a sigmoidal function, and (4) networks of so-called sigmoidal neurons work better than networks of perceptrons do.}

\mpar\mpar
Let's continue with our ``sharp nails'' scenario. The network can see both fingers and claws, so it's ambiguous, but perceptrons must choose one or the other. So, despite the ambiguity, a perceptron can only say its 100\% claws and there is 0\% resemblance to fingers, or perhaps vice versa. But what sigmoidal neurons can do is give a value between 0 and 1 that can take ambiguity into account. In this scenario, we could have the value of fingers be say, 0.31, the value of claws be 0.64, and the value of paws be 0.05.  While the network still detects claws with higher confidence than fingers or paws, all possibilities can be taken into account. We will see the importance of these numbers, called \textit{activations}, more clearly as we go along.

\mpar
You may have noticed from the figures that each neuron is connected to every neuron of the next layer to the right so you might ask: if the \textit{Hands} neuron from the input layer was fired then why couldn't \textit{Hands} fire at least partially towards \textit{Paw} and \textit{Finger}, not just \textit{Claw}? This could handle potential ambiguity about the kind of hands in the picture. Clearly perceptrons can't split things up like this but sigmoidal neurons can. This is where \textit{weights} come in. Weights are what express the relative importance of inputs to the final output.

\mpar
For example, let's take the \textit{Legs} neuron from our imaginary network and suppose it has a sigmoidal value of 0.9 (basically a 90\% confidence level in detecting legs in the image). The presence or absence of legs in the image does not impact the network's decision between \textit{Fur} and \textit{Head}, so the arrows connecting the \textit{Legs} with \textit{Fur} and \textit{Head} could have an extremely low but nonzero weight. But the presence of \textit{Legs} is necessary in determining if there are \textit{2 legs} or \textit{4 legs}  For this reason, there could extremely large weight between the \textit{Legs} and  \textit{2 Legs} and \textit{4 Legs} in the hidden layer. 

\mpar
Here is another example to demonstrate the importance of weights. To make the point clear, we only have to look at part of the network, as in Figure 6. 
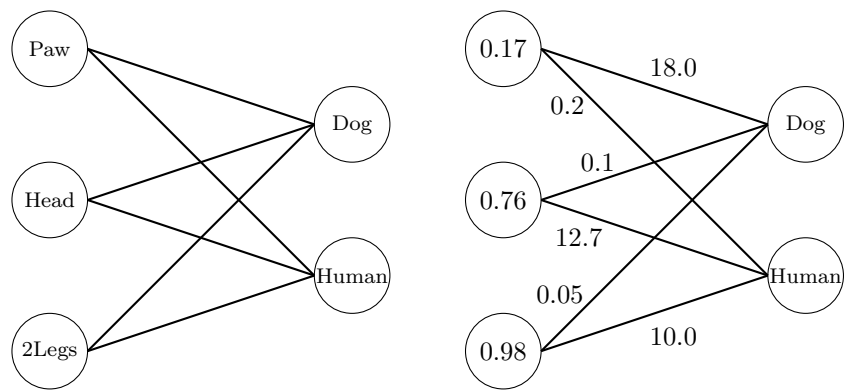
\begin{figure}[h]
\begin{tikzpicture}

\draw (0,0) circle (0.5) node{\footnotesize{2Legs}};
\draw (0,2) circle (0.5) node{\footnotesize{Head}};
\draw (0,4) circle (0.5) node{\footnotesize{Paw}};
\draw (4,1) circle (0.5) node{\footnotesize{Human}};
\draw (4,3) circle (0.5) node{\footnotesize{Dog}};
\draw[thick] (0.5,0) -- (3.5,1);
\draw[thick] (0.5,0) -- (3.5,3);
\draw[thick] (0.5,2) -- (3.5,1);
\draw[thick] (0.5,2) -- (3.5,3);
\draw[thick] (0.5,4) -- (3.5,1);
\draw[thick] (0.5,4) -- (3.5,3);
\draw (6,0) circle (0.5) node{0.98};
\draw (6,2) circle (0.5) node{0.76};
\draw (6,4) circle (0.5) node{0.17};
\draw (10,1) circle (0.5) node{\footnotesize{Human}};
\draw (10,3) circle (0.5) node{\footnotesize{Dog}};
\draw[thick] (6.5,0) -- (9.5,1);
\draw[thick] (6.5,0) -- (9.5,3);
\draw[thick] (6.5,2) -- (9.5,1);
\draw[thick] (6.5,2) -- (9.5,3);
\draw[thick] (6.5,4) -- (9.5,1);
\draw[thick] (6.5,4) -- (9.5,3);
\node at (8.25, 0.2) {10.0};
\node at (8.25, 3.75) {18.0};
\node at (6.75, 0.75) {0.05};
\node at (7.0, 1.5) {12.7};
\node at (7.25, 2.5) {0.1};
\node at (6.85, 3.25) {0.2};
\end{tikzpicture}
\caption{\textit{Left:} Part of the network. \textit{Right:} The same part of the network with neurons labeled with activations and arrows labeled with weights. }
\end{figure}

\mpar
In this scenario, we focus on three sigmoidal neurons in the hidden layer: \textit{Paw} with a value of 0.17 (not much resemblance to paws in the image), \textit{2Legs} with a value of 0.98 (almost 100\% sure there is a pair of legs), and \textit{Head} with a value of 0.76 (solid indication of a head of hair in the image). Paws are an indicator of dogs, which is why the weight of that arrow is high at 18. Paws do not indicate a human, so the weight there is only 0.2. On the opposite side, a head of hair is a strong indicator of humans (unless the person is bald or unless it's that freshly shaved poodle with a fluffy head of hair again) which is why the weight of the arrow to \textit{Human} is 10.0, while the weight of the arrow to \textit{Dog} is only 0.05. Finally, we have the \textit{2Legs} neuron. Two legs can may indicate either dogs and humans (for example, if the dog is doing a trick!), but it seems much more likely to be a human. So the weight of the arrow towards \textit{Dog} is quite low at 0.05 and the arrow toward \textit{Human} is quite high at 12.7.

\mpar
Here's one of the key steps. Let's aggregate all this information to get a score for \textit{Dog} and a score for \textit{Human}. If the scores are close, the network is having a hard time distinguishing between them and if the scores are very different then the network should chose the one with the higher score. We aggregate the information by multiplying the activations by the weights and adding them up over all the arrows. Here are the calculations,
\begin{align*}
Dog & = 0.17\times 18.0\, +\, 0.76\times 0.1\, +\, 9.8\times 0.05 = 3.596\\
Human & = 0.17\times 0.2\,+\, 0.76\times 12.7\,+\, .98\times 10.0 = 19.486
\end{align*}

On one hand, the score for \textit{Human} is significantly bigger than the score for \textit{Dog} so it seems natural to conclude the picture is of a human. On the other hand, we may want more evidence for concluding there's a person in the picture rather than a dog. This is where \textit{bias} comes in. Biases are thresholds and being above threshold is good. If the threshold was 5.0 for a dog, let's say, and 9.0 for a human then \textit{Dog} is below threshold and \textit{Human} is above.

\mpar
But how would you tell a neuron it needs a value of more than 5.0 or more than 9.0? The easiest way to demonstrate how bias works is by subtraction. Bias is a number which subtracts from the sum of weights and activations multiplied together. For this situation, if the final score after bias is subtracted is negative, the new activation is set to zero. But if the score is positive, the new activation is exactly that positive number. This process is captured by introducing a special function called RELU for Rectified Linear Unit,
$$
RELU(x) = \begin{cases}
x, & \text{if $x\geq 0$,}\\
0, &\text{if $x<0.$}
\end{cases}
$$
Then, the final activation for \textit{Dog} is $a_{Dog}=RELU(3.596 - 5)=0$ and the final activation for \textit{Human} its $a_{Human}=RELU(19.486-9)=10.496.$ 

\mpar
That's the arithmetic of the situation but it's important to look at the corresponding algebra. What does that mean? Well, it just means finding an abstract way of writing these calculations using symbols, rather than numbers. That way we have a hope of extending our insights from this cartoon model to more realistic neural networks. 

\mpar
Notice that each term in the score calculation has the form $w_{ij}\times a_j$ where $w_{ij}$ is a weight attached to the arrow connecting neuron $j$ in the hidden layer with neuron $i$ in the output layer and $a_j$ is the activation of neuron $j$ in the hidden layer. For example, $$
a_{Paw}= 0.17 \quad\text{and}\quad w_{Dog, \,Paw}=18.0.
$$
Then we subtract the bias $b_i$ and apply RELU. Rather than subtracting the bias we think of the bias as negative so we add it. Thus,
$$
b_{Dog}=-5\quad\textit{and}\quad b_{Human}=-9
$$ 
This neatens up the calculations so there aren't so many negative signs to keep track of. Altogether we can write this as,
$$
a_i=RELU\left(\sum_j w_{ij}a_j +b_i\right)
$$

To be completely specific, for each index $i=Dog,\, Human$ the sum is over indices $j= Paw,\, Head,\, 2Legs.$  So, if you stare at this formula long enough you'll see that it's just telling you how to do the calculations we did to get the activations for the neurons  \textit{Dog} and \textit{Human} without using the actual numbers.\,

\mpar
In calculating the activation of a neuron it's important to keep the functional form,
$$
\textit{activation = function of (sum of $weights\times activations$ plus bias)}
$$
It turns out that we can use functions other than $RELU.$ Observe that  $RELU$ can have very large outputs, if it has large inputs, so the value of an activation could get out of hand in a very large network. One way to deal with this issue is to cut off RELU at a fixed value, which we take to be the value 1 for convenience. This leads to the function,

$$
RELU^+(x) = \begin{cases}
1, & \text{if $x\geq 1,$}\\
x, & \text{if $0\leq x\leq 1$,}\\
0, &\text{if $x\leq 0.$}
\end{cases}
$$
Using $RELU^+$ ensures that the activation of every neuron is a number between 0 and 1, which is helpful.

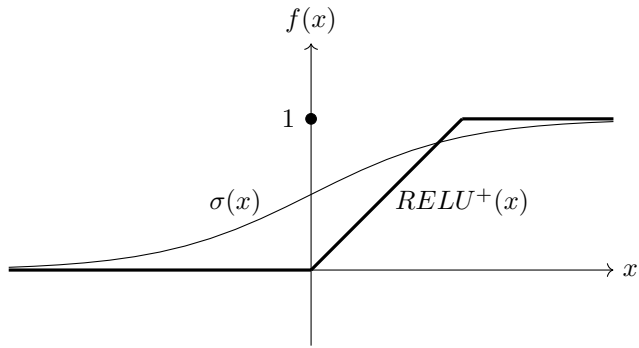
\begin{figure}[h]
\begin{tikzpicture} [domain =-4:4]
\draw[->](-3,0) -- (4,0) node[right] {$x$};
\draw[->](0,-1) -- (0,3) node[above] {$f(x)$};

\draw[thin] plot(\x, {2/(1+ exp(-\x))}); 
\draw[very thick] (-4,0)--(0,0);
\draw[very thick] (0,0)--(2,2);
\draw[very thick] (2,2)--(4,2);
\node at (2,0.9) {$RELU^+(x)$};
\node at (-1,0.9) {$\sigma(x)$};
\draw[fill=black] (0, 2) circle (2pt);
\node at (-0.3,2) {1};
\end{tikzpicture}
\caption{The graph of $\sigma(x)=1/(1+e^{-x})$ is in grey and the graph of $RELU^+(x)$ is in black. Both functions increase from 0 to 1 but $RELU^+$does it abruptly while $\sigma$ does it smoothly.}
\end{figure}

\mpar
When we discuss how networks are trained on data we'll see that a function whose graph is a smoothed out version of the graph of $RELU^+$ has certain advantages. One example of such a smoothed out version of $RELU^+$ is function is the sigmoid function given by the formula,
$$
\sigma(x)=\frac{1}{1+e^{-x}}.
$$

Observe that when $x$ is large and positive, $e^{-x}$ is very close to 0 hence $\sigma(x)$ is close to 1. On the other hand, when $x$ is large and negative, $e^{-x}$ is large and positive hence $\sigma(x)$ is close to 0. See Figure 7 for the graphs of $RELU^+$ and $\sigma.$

\mpar
This particular function $\sigma$ is the one we will use but other choices are possible. In general, a \textit{sigmoidal} neuron is a neuron whose activations are calculated using a \textit{sigmoidal} function which is any function, like $\sigma,$ that increases from 0 to 1 and whose graph is a smooth $S$-shaped curve. Summing it all up, 

\begin{itemize}
\item{if $i$ indexes a neuron in a given layer and $a_i$ is its activation}
\item{if $j$ indexes neurons in the previous layer and their activations are $a_j$}
\item{if $w_{ij}$ is the weight of the arrow from neuron $j$ to neuron $i$}
\item{if $b_i$ is the bias of neuron $i$}
\item{and if $\sigma$ is a sigmoidal function then},
\end{itemize}
$$
a_i=\sigma\left(\sum_j w_{ij}a_j +b_i\right).
$$

\mpar\mpar
\subsection{Training a Neural Network} \textit{The main things to know and remember from this subsection are (1) a neural network is a function of all the weights and biases that are used to calculate neural activations, (2) training data consist of pairs $(x,y)$ where $x$ is the input and $y$ is the output which we think of as providing the ground truth, (3) for a given input $x$ the network calculates an output $\hat{y}$(x),  called the predicted output, which is a function of the current weights and biases, (4) network performance is measured by comparing the true output $y$ with the predicted output $\hat{y}$, (5) learning occurs by changing the weights and bias so as to improve network performance, and (6) the algorithm for changing the weights and biases uses calculus and linear algebra.}

\mpar\mpar
The purpose of the imaginary \textit{Dog-Bear-Human} network was to create  the mental image of a neural network as a collection of sigmoidal neurons arranged in multiple hidden layers that connect the input layer to the output layer. By some means the hidden layers are able to detect patterns in the input data, relate the patterns, and draw conclusions from them all by modifying neural activations using the recipe above. How does that actually happen? That's where mathematics comes in. If you don't recognize the math, don't worry, we'll review it in the next sections as needed and lay it out coherently in appendices. For the time being, just read on as if you do!

\mpar
In practical situations we have a collection of input/output pairs $(x,y)$ that will be used to train the network. Input $x$ is a string of $d$ numerical values $x=(x_1, x_2, \cdots, x_d)$ and output $y$ is a string of $D$ numerical values $y=(y_1, y_2, \cdots, y_D).$ When we want to indicate that output $y$ belongs to input $x$ we write $y=y(x).$ Numerical strings like this are called \textit{vectors} and the individual numbers are called \textit{coordinates}. The number of coordinates in the string is called the \textit{dimension} of the vector so $x$ is a vector of dimension $d$ and both $y$ and $\hat{y}$ are vectors of dimension $D.$ 

\mpar
For example, $x$ might be the greyscale numbers of all the pixels in a digital image of an animal so it could potentially be a vector with millions of dimensions. On the other hand, in our imaginary network $y$ is a three dimensional vector with one dimension each for dog, bear and human. We imagine that human experts have already classified the data by labelling a picture of a dog with the vector $(1,0,0)$, a picture of a bear with $(0,1,0),$ and a picture of a human with $(0,0,1).$ We assume humans have done a similar classification in the general case where there are $D$ possible outcomes. So, this establishes \textit{ground truth} by providing many examples of correct image classification.

\mpar
The first step is to think of the network as a function that accepts  a vector $x$ as input and produces another vector $\hat{y}=(\hat{y}_1, \hat{y}_2, \cdots, \hat{y}_D)$ as output. The vector $\hat{y}$ is called the \textit{predicted output} or sometimes the \textit{predicted value.} Obviously, the predicted output is a function of $x$ and we show this dependence by writing $\hat{y}=\hat{y}(x).$ It's less obvious but extremely important to realize $\hat{y}$ is also a function of all the weights and biases of the network. Remember there is a bias for each neuron (except for neurons in the input layer) and a weight for each arrow between neurons.  Even for our cartoon network there are 10 biases and $6\times 7 =42$ weights, so $\hat{y}$ is a function of a lot of variables even in the simplest cases.

\mpar
Let the symbol $W$ stand for the set of the weights of all arrows in the network and the symbol $B$ stand for the set of the biases of every neuron. Then the predicted output is actually a function of the variables $x, W,$ and $B$  and we sometimes write,
$$
\hat{y}=\hat{y}(x) = \hat{y}(x, W,B)
$$
to make it clear that the coordinates of the vector $\hat{y}$ depend on the coordinates of the vector $x$ and also on the coordinates of $W$ and $B$.

\mpar
The next step is to assume the numerical values for the coordinates in $W$ and $B$ have been given. With those numbers kept fixed, the predicted value is now just a function of $x$ and we can compare predicted and actual outputs. There are many ways to make the comparison but the simplest is the sum of the squares of the difference in coordinates,

$$
\delta^2(x)=|\hat{y}_1(x)_1-y(x)_1|^2+ |\hat{y}(x)_2-y(x)_2|^2 + \cdots + |\hat{y}(x)_D-y(x)_D|^2.
$$

\mpar
Now, suppose there are $N$ input/output pairs in the training data set. Then we could average $\delta^2(x)$ over all inputs in the training set can get a quantity,
$$
\Delta^2 =\frac{1}{N}\sum_{n=1}^N\delta^2(x_n)
$$
which is a number that measures network performance. If $\Delta^2$ is small then predicted outputs are close to actual outputs on average, so the classification is pretty good. But if $\Delta^2$ is large then predicted outputs are generally not close to actual outputs so the performance leaves a lot to be desired.

\mpar
So, here's the key insight that makes network training possible. At first glance, $\Delta^2$ is just a number but in reality it is a function of the weights and biases,
 $$
 \Delta^2=\Delta^2(W,B).
 $$
What this means is $\Delta^2(W,B)$ changes as the weights and biases change. What if we could find a way to update the current weights and biases by adding or subtracting a little bit from each number so that $\Delta^2$ shrinks? That would mean the network performs better with the updated weights than with the previous weights, which is good. 

\mpar
But why stop there? Maybe there's a way to repeat the updating procedure on the new weights and biases to get newer weights and biases that perform even  better than before. The hope is that by repeating the update procedure many times we'll get to the point where we're satisfied with the network's performance. 

\mpar
This is where some serious mathematics enters the picture. In abstract terms we have a function of many variables (the weights and biases) and we're trying to find values of the variables that minimize the function. For the function $\Delta^2$ the minimum value is 0 which corresponds to the network getting it exactly right for every data point in the training set: $\hat{y}(x_n) = y(x_n),$ for all $1\leq n\leq N.$  In fact, weights and biases that make $\Delta^2=0$ is actually a bad thing, which may seem counterintuitive at first. It turns out that being super-accurate on the training set, what is called \textit{overfitting}, makes the network less able to handle data it hasn't seen yet, which defeats the whole purpose. We're not looking for perfection here, we're just looking for a way to find values $W_*$ and $B_*$ of the weights and biases such that $\Delta^2(W_*,B_*)$ is reasonably small and further updating doesn't really improve things.

\mpar
There is a calculus-based procedure that applies to situations like this called the \textit{method of steepest descent} and it's the key to successful training of neural networks. There are some serious challenges in getting the method to work for industrial scale neural networks and it took a lot of research effort to find a way. The breakthrough idea that overcomes the challenges of implementing gradient descent for neural networks is something called \textit{backpropagation}. We'll examine gradient descent and backpropagation thoroughly  in the next section.

\mpar\mpar
\section{How to Build a Neural Network in Five Easy Steps}

In flowchart terms a neural network has the following structure: \textit{initialization, feedforward, backpropagation, mini-batch,} and \textit{evaluation}. Let's discuss them briefly to see how they fit together then dig into each one in a separate subsection. This is where the mathematics comes in: linear algebra, calculus, and probability. 

\mpar
If you're not familiar with some, or any, of the mathematical ideas and notations, just accept them for now. We'll provide more detailed mathematical background a little later and you can always dig into the appendices first if you're curious right now. Either way works.

\mpar
First, the network must be created by choosing the size of the input and  output layers as well as the number of hidden layers and their sizes. These choices depend significantly on the form of the data to be used. Then corresponding weights and biases are randomly chosen and the data set of numerical input/output pairs is split into a training set and an evaluation set. This step is called \textit{initialization}. 

\mpar
Second, we describe how the weights and biases transform the input data from the training set into activations in the first layer, then from activations in the first layer to activations in the second layer, etc., eventually creating activations in the output layer which we call the \textit{predicted output} or sometimes the \textit{predicted value}. This process of taking data in the input layer and creating a predicted value in the output layer is called the \textit{feedforward} step. 

\mpar
It's good to think of feedforward as the assignment of a predicted value to an input value as a function of all the weights and biases collectively. It's natural to ask: how close is the predicted value to the actual value?  There are many ways to answer this question but for simplicity's sake we choose something called the \textit{mean squared error}, 
$$
\Delta^2=\frac{1}{N}\sum_{n=1}^N \delta^2(x_n)
$$  
Since $\Delta^2$ is a function of the weights and biases we write $\Delta^2=\Delta^2(W,B)$ to make this dependence explicit. 

\mpar
Third, we need to find a systematic way to update the weights and biases so that the new mean squared error $\Delta^2(W_{new}, B_{new})$ is smaller than the old mean squared error $\Delta^2(W_{old}, B_{old}).$ An important insight from multivariable calculus is that the \textit{gradient} of $\Delta^2$ evaluated at the current values of the weights and biases points in the direction of greatest increase of $\Delta^2,$ meaning the negative of the gradient points in the direction of greatest decrease of $\Delta^2.$ This suggests an update rule,
$$
(W_{new}, B_{new}) - (W_{old}, B_{old}) = -\eta\times\nabla \Delta^2(W_{old}, B_{old}),
$$

\mpar
where $\eta$ is a small, adjustable parameter called the \textit{learning rate} and $\nabla$ is the gradient operator. This approach, called \textit{gradient descent}, is an applied mathematics workhorse with an excellent track record of solving optimization problems in science and engineering. The gradient $\nabla \Delta^2$ is a collection of partial derivatives, one for each weight $w_{ij}$ and bias $b_i.$ Calculating these partial derivatives efficiently exploits the left-to-right nature of network architecture and results in the famous \textit{backpropagation} algorithm.

\mpar
Fourth, there is \textit{mini-batch}, a stochastic variant of the classic \textit{method of steepest descent} that takes advantage of randomness and the special nature to the mean squared error function. Classic gradient descent applied to $\Delta^2$ would require calculating partial derivatives with respect to all the weights and biases of,
$$
\delta^2(x_n)=||\hat{y}(x_n)-y(x_n)||^2
$$ 
for all $1\leq n\leq N.$ For a network to learn anything useful there has to be a large number of training data, meaning $N$ is a very large number, so that a lot of partial derivatives! To reduce to cost of computation, we choose a small percentage of the data points at random - called a \textit{mini-batch} -  and calculate partial derivatives of $\delta^2$ just for them. This is basic idea behind the method of mini-batch gradient descent: 

\mpar\mpar
\textit{Apply one step of classic gradient descent to a random approximation of $\Delta^2$ obtained by averaging $\delta^2(x_{r_j}),\, 1\leq j\leq R,$ where $x_{r_1}, x_{r_2},\, \dots,\, x_{r_R}$ is a randomly chosen subset of $R$ training data points, and repeat until we've used all the data.}

\mpar\mpar
Fifth, is the \textit{evaluation} step where we use the evaluation data set to measure network performance. This means measuring the extent to which predicted outputs $\hat{y}(x)$ match actual outputs $y(x)$ on data the network hasn't seen yet. The accuracy of the elavuation step demonstrates the extent to which the network has learned from the training data.

\mpar
Let's look more carefully at each step in the subsections below. Backpropagation and mini-batch are the most mathematically difficult steps which is why those subsections are much longer than the others. 

\mpar
Each subsection begins with an italicized statement indicating the kinds of mathematics that will be used. Readers are encouraged to consult the appendices to supplement their current knowledge, as needed. Appendices have multiple sections so, for example, $A(i)$ refers to section $(i)$ of appendix $A.$ The appendices are written in a cumulative fashion so you may need to refer back to previous sections, and perhaps previous appendices, when reading through a given section.
\mpar\mpar
\subsection{Initialization}

\textit{This subsection uses some ideas from linear algebra; especially vectors A(i), higher dimensions A(iv), and matrices and linear transformations A(iv). It also uses the derivative B(ii) and normal random variables C(v).}
\mpar

\mpar
Network architecture refers to the organizational structure of the neural network. We can choose how many neurons there are in the input and output layers as well as the number of hidden layers and the number of neurons in a given hidden layer. 

\mpar
We make three simplifying assumptions. First, layers are ordered from left to right with the input layer on the left, the output layer on the right, and the hidden layers in between. Second, there are the same number of neurons in each hidden layer and third, that each neuron in a given layer is connected to every neurons in the layer to its immediate right. Other choices can be made but these are the simplest.

\mpar
Thus, network architecture is specified by integers: \texttt{InputNeuron, OutputNeuron, \\HiddenLayer,} and \texttt{HiddenNeuron}. For example, if we wish to create a neural network with 5 hidden layers of 11 nodes each, with an input layer of 26 nodes and an output layer of 2 nodes then we would declare that,

$$
\texttt{InputNeuron} =26,\, \texttt{OutputNeuron}=2,\, \texttt{HiddenLayer}=5,\, \texttt{HiddenNeuron}=11.
$$

\mpar
We think of neurons as containing a number called its \textit{activation} and we write, following Nielsen,  $a^l_j$ for the activation of the $j^{th}$ neuron in the $l^{th}$ layer. The superscript $l=0$ refers to the input layer and $l=L$ refers to the output layer, so there are $L+1$ layers altogether. Thus, $a^0$ is a vector of dimension \texttt{InputNeuron}, $a^L$ is a vector of dimension \texttt{OutputNeuron}, and for each $1\leq j\leq\texttt{HiddenLayer},\, a^l$ is a vector of dimension \texttt{HiddenNeuron}.

\mpar
We think of each edge as having a \textit{weight} and we write $w^l_{jk}$ for the weight of the arrow connecting the $k^{th}$ neuron in the $(l-1)^{st}$ layer to the $j^{th}$ neuron in the $l^{th}$ layer. Thus each $w^l$ is a matrix. Observe that $w^1$ connects the input layer to the first hidden layer and $w^L$ connects the last hidden layer to the output layer. Thus $w^1$ is a matrix of dimension \texttt{HiddenNeuron} $\negthickspace\times\negthickspace$ \texttt{InputNeuron}, $w^L$ is a matrix of dimension \texttt{OutputNeuron} $\negthickspace\times\negthickspace$ \texttt{HiddenNeuron}, and for all other values of $l,\, w^l$ is a matrix of dimension \texttt{HiddenNeuron} $\negthickspace\times\negthickspace$ \texttt{HiddenNeuron}. The collection of these matrices is denoted,
$$
W=(w^1, w^2, \,\dots,\,w^L).
$$

\mpar
We think of each neuron, except those in the input layer, as having it's own \textit{bias} and we write $b^l_j$ for the bias of the $j^{th}$ neuron in the $l^{th}$ layer. Thus  each $b^l$ is a vector. Observe that $b^L$ has dimension \texttt{OutputNeuron} while for each $1\leq l\leq L-1,\, b^l$ is a vector of dimension \texttt{HiddenNeuron}. The collection of these vectors is denoted,
$$
B=(b^1, b^2,\,\dots,\, b^L).
$$

\mpar
To simplify notation in the following subsections, let $d=\texttt{InputNeuron}, \,\, D=\texttt{OutputNeuron},$ and $H=\texttt{HiddenNeuron}.$

\mpar
The network has an \textit{activation function} $\sigma$ which is a function $\sigma\colon\mathbb{R}\to \mathbb{R}.$ There are many choices for $\sigma$ but we choose the sigmoidal function,
$$
\sigma(x)=\frac{1}{1+e^{-x}}
$$
for two reasons: (1) $0<\sigma(x)<1$ and (2) the graph of $\sigma$ is smooth and increasing. It turns out that the derivative of $\sigma$ occurs in the backpropagation subsection so it's helpful to know something about it.

\mpar
By a standard calculation the derivative is,
$$
\sigma^\prime (x)= \frac{e^{-x}}{(1+e^{-x})^2} =\frac{1}{\left(e^{\,x/2}+e^{-x/2}\right)^2}.
$$ 
Observe that $\sigma^\prime(x)=\sigma^\prime(-x)$ and $\lim_{x\to -\infty}\sigma^\prime(x) =\lim_{x\to \infty}\sigma^\prime(x) =0.$ It's not too hard to check that $\sigma^\prime(x)$ is increasing on $(-\infty,0)$ and decreasing on $(0,\infty)$ which implies $\sigma^\prime(x)$ attains it's largest value at $x=0.$ Thus, $0<\sigma^\prime(x)\leq\sigma^\prime(0)=\tfrac{1}{4}.$

\mpar
Training is an iterative process of modifying weights and biases based on network performance on the training data set. So, it's natural to ask how does the whole process get started? We have to have some initial values for weights and biases before we can update them so where do initial values for weights and biases come from in the first place? 

\mpar
The basic idea is to choose them randomly by saying that $w^l_{jk}$ and $b^l_j$ for all values of $j, k$ and $l$ are independent, identically distributed normal random variables of a chosen mean $\mu$ and a chosen variance $\sigma^2.$ Typical choices are $\mu=0$ and $\sigma=1.$ Note that there are two $\sigma$'s floating around here: an activation function and a standard deviation. This is ambiguous notation but it's usually clear from the context which $\sigma$ is which.

\mpar
The network is intended to learn by searching for patterns in data. This implies we have a set of data points of the form $(x,y)$ where $x$ is a vector thought of as the input and $y$ is another vector thought of as the output. Our attitude is that this pairing represents the ground truth. For example, ground truth could mean $x$ is the pixels of a digital image of human lungs and $y$ is code for a consensus diagnosis made by expert radiologists. Thus $x$ is a $d$-dimensional vector (where $d$ is the product of the number of pixels in a radiograph times the number of possible greyscales in a pixel) and $y$ is a $D$-dimensional vector (where $D$ is the number of possible diagnoses). These data are split in two subsets: $N$ will be used to train the network and $M$ will be used to evaluate the network's performance. 

\mpar
To summarize initialization,

\begin{itemize}
\item\texttt{TrainingData} and \texttt{EvaluationData} are lists of input/output tuples of the form $(x,y)$ where,
\begin{itemize}
\item $x$ is vector of dimension $d=\texttt{InputNeuron}$
\item $y$ is a vector of dimension $D=\texttt{OutputNeuron}$ 
\end{itemize}
\item  $N$ = \texttt{NumTrainingData} is the number of points in the training set
\item $M$ = \texttt{NumEvaluationData} is the number of points in the evaluation set
\item Choose a sigmoidal function $\sigma$
\item Randomly initialize all the weights and biases

\end{itemize}	

\mpar
\subsection{Feedforward}
\textit{This subsection uses addition and subtraction of vectors A(ii) and matrices and linear transformations A(v).}

\mpar
Activation flows from left to right through a neural network, beginning at the input layer and ending at the output layer. The activation $a^l_j$ is calculated in two steps. 

\mpar
First, the activations of neurons in the previous layer, layer $l-1,$ are multiplied by the weights of the edges connecting them to the $j^{th}$ neuron in layer $l,$ then added together along with the bias of the neuron's bias. This defines an affine function,
\begin{align*}
z_j^i(a^{l-1}) & =w^l a_j^{l-1}+b_j^l\\
& =\sum_{k=1}^{N_{l-1}} w^l_{jk} a^{l-1}_k +b^l_j
\end{align*}
which we call the \textit{score} of neuron $j.$ Here, $N_{l-1}$ is the number of neurons in layer $l-1, \, w^l$ is the $N_l\times N_{l-1}-$dimensional matrix,
$$
w^l =\left(
\begin{matrix}
w^l_{11} & w^l_{12} & \dots & w^l_{1\, N_{l-1}}\\
w^l_{21} & w^l_{j2} & \dots & w^l_{2 N_{l-1}}\\
\vdots & \vdots &\ddots &\vdots\\
w^l_{N_l\,1} & w^l_{N_l\, 2} & \dots & w^l_{N_l\, N_{l-1}}
\end{matrix}
\right)
$$
and $a^{l-1}$ is the $N_{l-1}-$dimensional vector of activations,
$$
a^{l-1} = (a^{l-1}_1, a^{l-1}_2, \cdots, a^{l-1}_{N_{l-1}}).
$$
To be specific $N_0=d=\texttt{InputNeuron},\, N_L=D=\texttt{OutputNeuron},$ and otherwise $N_l=H=\texttt{HiddenNeuron}.$

\mpar
We think of the weights $w^l_{jk}$ and biases $b^l_j$ as given and the activations $a^{l-1}_{k}$ as input variables to the function $z_j(\cdot).$ That's the first step. The second step is to define the activation $a^l_j$ as the composition of the activation function $\sigma$ with this affine function,
$$
a^l_j  =\sigma\circ z_j^l\, (a^{l-1}) = \sigma\left(\sum_{k=1}^{N_{l-1}} w^l_{jk} a^{l-1}_k +b^l_j\right).
$$

The shorthand notation for this is,
$$
a^l=\sigma\circ z^l(a^{l-1}) = \sigma(z^l(a^{l-1}))= \sigma\left(w^la^{l-1}+b^l\right).
$$
Note that this shorthand is technically incorrect. The input to $\sigma$ must be a number like $z^l_j$, not a vector like $z^l.$ By implication, the shorthand means: apply $\sigma$ to each of the coordinates of $z^l.$

\mpar
Thus, the activations of neurons in a given layer are functions of the activations of neurons in the previous layer which, in turn, are functions of the activations of neurons in the layer before that, etc., all the way back to the input layer. 

\mpar
Turning this around, the numerical values provided by training data to neurons in the input layer \textit{feed forward} through all the hidden layers to produce numerical values for the activations in the output layer. What's more, the numerical values in the output layer are a function not only of the input values but also of all the weights and biases of the network. This point of view is crucial to understanding how neural networks are trained.

\mpar
To summarize feedforward,
\begin{itemize}
\item{$d=\texttt{InputNeuron},\, H=\texttt{HiddenNeuron},\,D=\texttt{OutputNeuron},\\ L=\texttt{HiddenLayer}+1$}
\item{The coordinates of the $d$-dimensional vector $x$ of training data are the activations of the $d$ neurons of the input layer}
\item{The activations of neurons in the $l^{th}$ hidden layer are calculated from activations of the $(l-1)^{st}$ hidden layer by the formula,}
$$
a^l_j=\sigma(z^l_j(a^{l-1})) = \sigma\left(\sum_{k=1}^Hw^l_{jk}a^{l-1}_k +b^l_j\right)
$$
\item{The predicted putput, or predicted value, $\hat{y}(x)$ of input $x$ is the $D$-dimensional vector of activations in the output layer: $\hat{y}(x)=a^L$}
\end{itemize}

\subsection{Backpropagation}
\textit{This subsection uses the distance and angle between vectors A(iii), the chain rule B(vi), directional derivatives and the gradient B(viii), the method of gradient descent B(ix).}

\mpar
To fix ideas let's look at the squared distance between the predicted output $\hat{y}(x)$ and the the actual output $y(x)$ for a single input data point $(x,y),$  namely,
$$
\delta^2(x) = || \hat{y}(x)-y(x)||^2 =\sum_{j=1}^D |\hat{y}_j(x) -y_j(x)|^2. 
$$
Recall that $\hat{y}(x)$ is a function of all the weights and biases in $W$ and $B.$ The essential step in applying the method of gradient descent to the mean squared error $\Delta$ is to calculate the gradient of $|| \hat{y}(x)-y(x)||^2 $ whose components are the partial derivatives,
$$
\frac{\partial}{\partial\, w^l_{j,k}}\, || \hat{y}(x)-y(x)||^2\quad\text{and}\quad \frac{\partial}{\partial\, b^l_j }   || \hat{y}(x)-y(x)||^2 
$$
for all possible values of $j,k,$ and $l.$ Since, by the chain rule,
$$
\frac{\partial}{\partial\, w^l_{j,k}}\, || \hat{y}(x)-y(x)||^2 =2\sum_{j=1}^D\left(|\hat{y}_j(x) -y_j(x)|\times\frac{\partial}{\partial w^l_{jk}}\,\hat{y}_j(x)\right)
$$
this comes down to calculating the partial derivative of $\hat{y}_j(x) = \hat{y}_j(x,W,B)$ with respect to $w^l_{jk}.$ Likewise for the  partial derivative of $\hat{y}_j(x) = \hat{y}_j(x,W,B)$ with respect to $b^l_j.$
\mpar

OK, we can do that in principle, but calculating these partial derivatives in an effective way depends on understanding how the predicted value $\hat{y}_i$ depends on all the network's weights and biases. And we have to understand these dependencies well enough to apply the standard rules of calculus to get a formula we can evaluate on a computer. 

\mpar
That's the bottom line. If this can't be done then we can't train the network. If this can be done but it's too expensive or takes too much time, then training the network is a theoretical possibility but a practical impossibility. In historical terms, this was a choke point for the field of neural networks until a new computational method came along called \textit{backpropagation} that made the computational cost of taking one step in the gradient descent algorithm proportional to the computational cost of calculating the predicted value in the first place. 

\begin{remark}
\textit{The paper that brought backpropagation into prominence is Rummelhart, Hinton and Williams \cite{RHW}. The underlying idea had been discussed previously by several authors, but the use of backpropagation in the context of neural networks was firmly established in this publication. Geoffery Hinton was a co-recipient of the 2024 Nobel Prize in physics for his pioneering work in neural networks and machine learning.}
\end{remark}

The underlying idea of  backpropagation is to take advantage of the special structure of the \textit{feedforward} function $x\to\hat{y}(x,W,B)$ that computes the predicted value given the input value. Let's take some time to understand this special structure. It isn't all that hard but it is notationally tricky and potentially confusing since there are a lot of moving parts. Some patience and diligence are required of the reader to make sense of it all.

\mpar
The key idea is to work backwards through the feedforward step. At the last step, all the activations in the last hidden layer are known and the predicted value $\hat{y}$ has the abstract form,
$$
\hat{y} =\sigma\circ z^L=\sigma \left(w^La^{L-1} +b^L\right)
$$
where $a^{L-1}$ is the $H$-dimensional vector of activations in the last hidden layer, 
$$
w^L = (w^L_{jk}), \,\,1\leq j\leq D,\,\, 1\leq k\leq H,
$$ 
is the matrix of weights associated with the edges connecting the $H$ neurons in the last hidden layer to the $D$ neurons in the output layer and, 
$$
b^L=(b^L_1, b^L_2,\,\dots,\, b^L_D)
$$ 
is the $D$-dimensional vector of biases of the output layer neurons. 

\mpar
It's helpful to introduce the so-called \textit{score}, $z^L=w^La^{L-1}+b^L,$ which is the $D$-dimensional vector whose $j^{th}$ component is,
$$
z^L_j = \sum_{k=1}^M w^L_{jk} a^{L-1}_k + b^L_j.
$$
Evidently,  the $j^{th}$ component of the vector $\hat{y}$ is $a^L_j=\sigma(z^L_j).$ In this notation, $z^L\colon\mathbb{R}^H\to\mathbb{R}^D$ is an affine function  $z^L(v)=w^Lv+b^L$  and it's important to think of $z^L$ as a function of the weights and biases: $z^L(v)=z^L(v, w^L, b^L).$

\mpar
For clarity's sake let's restate the formula for $\hat{y}$ as,

\begin{align*}
\hat{y} & = \sigma(z^L)\\
& = (\sigma(z^L_1), \sigma(z^L_2),\cdots, \sigma(z^L_D))\\
& = (\thickspace\sigma\left(\sum_{k=1}^M w_{1k} a^{L-1}_k + b_1\right)\negthickspace, \sigma\left(\sum_{k=1}^M w^L_{2k} a^{L-1}_k + b^L_2\right)\negthickspace, \cdots, \sigma\left(\sum_{k=1}^M w^L_{Dk} a^{L-1}_k + b^L_D\right)).
\end{align*}

\mpar\mpar
There is an important insight here. Because the network feeds forward, the activations in the last hidden layer  don't involve the weights and biases of the output layer. This implies that varying $w^L_{jk}$ and $b^L_j$ leaves each $a^{L-1}_k$ unchanged, so we can think of $a^{L-1}_k$ as a constant when we compute partial derivatives with respect $w^L_{jk}$ and $b^L_l$. That's super important!

\mpar
Also the formula for $\sigma(z^L_1)$ only involves the variables $w^L_{1k}$ and $b^L_1,$ so varying the other weights and biases leaves $\sigma(z^L_1)$ unchanged. This means we can think of,
$\sigma(z^L_1)$ as constant when we compute partial derivatives with respect to $w^L_{jk}$ and $b^L_j$ for $j=2,3,\cdots, D.$ That's just as important because it drastically reduces the number of partial derivatives we have to compute. 

\mpar
Of course this observation applies just as well to $\sigma(z^L_2), \sigma(z^L_3),$ etc. So, the only non-zero partial derivatives of $\sigma(z^L_j)$ are the partial derivatives with respect to $w^L_{jk}$ and $b^L_j$ and these are easily calculated by the chain rule,
$$
\frac{\partial \,\sigma(z^L_j)}{\partial\, w^L_{jk}} = \sigma^\prime (z^L_j)\,\frac{\partial \,z^L_j }{\partial\, w^L_{jk}} =  \sigma^\prime (z^L_j)\, a^{L-1}_k,
$$
and,
$$
\frac{\partial \,\sigma(z^L_j)}{\partial\, b_j} = \sigma^\prime (z_j)\,\frac{\partial \,z^L_j }{\partial\, b^L_j} =  \sigma^\prime (z_j),
$$
since $\partial/\,\partial w^L_{jk}\ z_j^L=a^{L-1}_k$ and $\partial/\partial\, b^L_j\,z_j^L=1.$

\mpar
OK. So far, so good. Now what about taking partial derivatives with respect to the weights and biases in the last hidden layer or, equivalently, the next-to-last layer of the network? There are more layers involved now, so let's label the parameters associated with the output layer with the superscript $L$, for ``last,'' parameters associated with the last hidden layer with the superscript $L-1$ for ``next-to-last,'' the parameters associated with the next-to-last hidden layer with the superscript $L-2$ for ``second-to-last,'' etc. 

\mpar
This notational sleight-of-hand allows us to write a recursive formula for the activations,
$$
a^L=\sigma\circ z^L(a^{L-1})\quad\text{and}\quad a^{L-1}=\sigma\circ z^{L-1}(a^{L-2})
$$
and therefore, 
\begin{align*}
\hat{y} & = a^L\\
& = \sigma\circ z^L (a^{L-1})\\
& = \left(\sigma\circ z^L\right)\circ\left(\sigma\circ z^{L-1}\right)(a^{L-2}).
\end{align*}

\mpar
If you think about it for a minute you'll see this makes calculating partial derivatives with respect to $W^{L-1}$ and $b^{L-1}$ pretty straightforward. The activations $a^{L-2}$ don't involve these parameters because they haven't been fed forward to layer $L-1$ yet. We could say they are ``upstream'' from layer $L-1.$  And the function $z^L$ doesn't involve them either because it's defined only in terms of the ``downstream'' parameters associated with layer $L.$ All the action involves the ``current'' parameters in the term $\sigma\circ\ z^{L-1}.$

\mpar
Once you see this you realize it can be repeated - that's exactly what a recursive formula means! Taking it one step further we'd have,
\begin{align*}
\hat{y} & = a^L\\
& = \left(\sigma\circ z^L\right)\circ\left(\sigma\circ z^{L-1}\right)(a^{L-2})\\
& = \left(\sigma\circ z^L\right)\circ\left(\sigma\circ z^{L-1}\right)\circ\left(\sigma\circ z^{L-2}\right)(a^{L-3}),
\end{align*}
and so on, all the way back to the input data,
$$
\hat{y}(x)=\left(\sigma\circ z^L\right)\circ\left(\sigma\circ z^{L-1}\right)\circ\cdots\circ\left(\sigma\circ z^1\right)(x),
$$
where the superscript 1 refers to the first hidden layer. 

\mpar
This is such a pretty formula that it must be powerful in some way, if only we could see it. How can we use the repeated composition of relatively simples function to our advantage? How can we turn this advantage into an efficient method to calculate partial derivatives?

\mpar
It all turns on thinking about the chain rule in a new way. Recall derivative of the composition $f\circ g(x) = f(g(x))$ is,
$$
\frac{d}{dx} f\circ g(x) = f^\prime(g(x))\,g^\prime(x),
$$
and it's proved using an algebraic trick,
\begin{align*}
\frac{d}{dx} f\circ g(x) = & \lim_{h\to 0} \frac{f(g(x+h))- f(g(x))}{h}\\
& = \lim_{h\to 0} \frac{f(g(x+h))- f(g(x))}{h}\times\frac{g(x+h)-g(x)}{g(x+h)-g(x)}\\
& =\lim_{h\to 0} \frac{f(g(x+h))- f(g(x))}{g(x+h)-g(x)}\times\frac{g(x+h)-g(x)}{h}\\
& \lim_{h\to 0} \frac{f(g(x+h))- f(g(x))}{g(x+h)-g(x)}\times\lim_{h\to 0}\frac{g(x+h)-g(x)}{h}.
\end{align*}

Let's write this as,
$$
\frac{df}{dx} = \frac{df}{dg}\frac{dg}{dx},
$$
which just says,  

\mpar
\textit{The rate at which the function f(x) is changing with respect to x is the rate at which f(x) is changing with respect to g(x) times the rate at which g(x) is changing with repect to x.}

\mpar\mpar
The fact that $dg$ appears in numerator and denominator is a reminder of the algebraic trick we used in the proof.

\mpar
Let's apply this mnemonic device on a known formula to see if it's reliable. It all depends on seeing $z^L$ in two ways: first as an independent variable in its own right and second and a function of $w^L$ and $b^L$ with the activations $a^{L-1}$ kept fixed. For example, the activation $a^L$ is a function the variable $z^L$ and, switching perspectives,  $z^L$ is a function of the weights $w^L.$ Then we should have,
$$
\frac{d a^L}{d w^L}= \frac{d a^L}{d z_L}\,\frac{d z^L}{d w^L}.
$$
Does this make sense? Well, $ da^L/d w^L$ stands for the rate of change of activations $a^L_j$ with respect to weights $w^L_{jk}.$ Similarly $d a^L/ d z^L$ stands for the rate of change of activations $a^L_j$ with respect to the scores $z^L_j,$ and $dz^L/dw^L$ stands for the rate of change of the scores $z_j^L$ with respect to the weights $w^L_{jk}.$ Now we already calculated that,
$$
\frac{d a_j^L}{d w^L_{jk}} = \sigma^\prime (z^L_j)\, a^{L-1}_k.
$$

Does the mnemonic predict this formula? Let's see. We know $a^L_j=\sigma(z^L_j)$ and therefore,
$$
\frac{\partial a^L_j}{\partial z^L_j}  = \frac{\partial \sigma (z^L_j)}{\partial z^L_j}  = \sigma^\prime(z^L_j).
$$
And we know $z^L_j=\sum_k w^L_{jk} a_k^{L-1} + b^L_j$ and therefore,
$$
\frac{\partial z^L_j}{\partial w^L_{jk}} = a^{L-1}_k.
$$
Multiplying the two terms together we find,
$$
\frac{\partial\, a^L_j}{\partial\, w^L_{jk}} = \frac{d a^L_j}{d z^L}_j\,\frac{d z^L_j}{d w^L_{jk}} = \sigma^\prime (z^L_j)\, a^{L-1}_k,
$$ 
as it should. So, yes! The mnemonic device correctly predicts the formula for the partial derivatives of activations with respect to weights and it's easy to check it also correctly predicts the partial derivatives of activations with respect to biases.

\mpar
Just to be careful, observe that $\partial\, a_j^L/\partial\, w^L$ is not a number but rather a $D\times H$ dimensional matrix,
$$
\frac{\partial\, a^L}{\partial\, w^L} =\left(
\begin{matrix}
\sigma^\prime(z^L_1)\,a^{L-1}_1 & \sigma^\prime(z^L_1)\,a^{L-1}_2 & \dots & \sigma^\prime(z^L_1)\,a^{L-1}_H\\
\sigma^\prime(z^L_2)\,a^{L-1}_1 & \sigma^\prime(z^L_2)\,a^{L-1}_2& \dots & \sigma^\prime(z^L_2)\,a^{L-1}_H\\
\vdots & \vdots &\ddots &\vdots\\
\sigma^\prime(z^L_D)\,a^{L-1}_1 & \sigma^\prime(z^L_D)\,a^{L-1}_2& \dots & \sigma^\prime(z^L_D)\,a^{L-1}_H
\end{matrix}
\right).
$$

\mpar
Similarly, $\partial\, z^L/\partial w^L$ is a $D\times H$ diagonal dimensional matrix. Thus,

$$
\frac{\partial\, z^L}{\partial\, w^L} =\left(
\begin{matrix}
a^{L-1}_1 & a^{L-1}_2 & \dots & a^{L-1}_H\\
a^{L-1}_1 & a^{L-1}_2& \dots & a^{L-1}_H\\
\vdots & \vdots &\ddots &\vdots\\
a^{L-1}_1 & a^{L-1}_2& \dots & a^{L-1}_H
\end{matrix}
\right).
$$

\mpar
Finally, $\partial\, a^L/\partial\,z^L$ is an $D\times D$ dimensional matrix but it is a diagonal matrix since $a^L_j$ only depends on $z^L_j,$ not on $z^L_k$ for $k\neq j.$ Thus,

$$
\frac{\partial\, a^L}{\partial\, z^L} =\left(
\begin{matrix}
\sigma^\prime(z^L_1) & 0 & \dots & 0\\
0 & \sigma^\prime(z^L_2)& \dots & 0\\
\vdots & \vdots &\ddots &\vdots\\
0 & 0& \dots & \sigma^\prime(z^L_D)
\end{matrix}
\right).
$$

It's easy to check that,
$$
\frac{\partial\, a^L}{\partial\, w^L} = \frac{\partial\, a^L}{\partial\,z^L}\,\frac{\partial\, z^L}{\partial w^L}
$$

in the sense of matrix multiplication and it's a useful exercise to do it by hand. In particular, it is dimensionally correct since the product of a $D\times D$ matrix with a $D\times H$ matrix is a $D\times H$ matrix.

\mpar
These observations also hold for the partial derivatives with respect to the biases. In particular, $\partial\,a^L/\,\partial b^L$ is the $D\times D$ matrix,
$$
\frac{\partial\, a^L}{\partial\,b^L} = \left(
\begin{matrix}
\sigma^\prime(z^L_1) & \sigma^\prime(z^L_2) & \dots & \sigma^\prime(z^L_D)\\
\sigma^\prime(z^L_1) & \sigma^\prime(z^L_2) & \dots & \sigma^\prime(z^L_D)\\
\vdots & \vdots &\ddots &\vdots\\
\sigma^\prime(z^L_1) & \sigma^\prime(z^L_2) & \dots & \sigma^\prime(z^L_D)
\end{matrix}
\right),
$$

and $\partial\, a^L/\partial z^L$ is the $D\times D$ matrix,
$$
\frac{\partial\, a^L}{\partial\,z^L} = \left(
\begin{matrix}
\sigma^\prime(z^L_1) & \sigma^\prime(z^L_2) & \dots & \sigma^\prime(z^L_D)\\
\sigma^\prime(z^L_1) & \sigma^\prime(z^L_2) & \dots & \sigma^\prime(z^L_D)\\
\vdots & \vdots &\ddots &\vdots\\
\sigma^\prime(z^L_1) & \sigma^\prime(z^L_2) & \dots & \sigma^\prime(z^L_D)
\end{matrix}
\right),
$$

and $\partial\, z^L/\partial\, b^L$ is the $D\times D$ just  identity matrix,

$$
\frac{\partial\, z^L}{\partial\,b^L} = \left(
\begin{matrix}
1 & 0 & \dots & 0\\
0 & 1 & \dots & 0\\
\vdots & \vdots &\ddots &\vdots\\
0 & 0 & \dots & 1
\end{matrix}
\right).
$$
 
Thus, in terms of matrix multiplication, 

$$
\frac{\partial\, a^L}{\partial\,b^L} = \frac{\partial\, a^L}{\partial\,z^L}\, \frac{\partial\, z^L}{\partial\, b^L}.
$$

\mpar
Our point of view will be to have confidence in the mnemonic device and use it to predict formulas for the components of the gradient of $\delta^2(x)$ with respect to the weights and biases of the network. Done carefully, this leads to formulas that can be checked with real mathematical rigor and can be easily coded into a computer. 

\mpar
To gain a bit more experience, let's calculate the partial derivatives of $a^L$ with respect to the weights and biases of the last hidden layer, $w^{L-1}$ and $b^{L-1}.$  Using the mnemonic device twice we have,
$$
\frac{\partial\,a^L}{\partial\, w^{L-1}} = \frac{\partial\, a^L}{\partial\, z^L}\,\frac{\partial\, z^L}{\partial\, a^{L-1}}\,\frac{\partial\, a^{L-1}}{\partial\, w^{L-1}}
$$
and,
$$
\frac{\partial\, a^L}{\partial\, b^{L-1}} =  \frac{\partial\, a^L}{\partial\, z^L}\frac{\partial\, z^L}{\partial\, a^{L-1}}\frac{\partial\, a^{L-1}}{\partial\, b^{L-1}}.
$$

\mpar\mpar
In each case, the right hand side is the product of three matrices of the appropriate dimensions.

\mpar
Here's why backpropagation is so powerful. We already calculated $\partial a^L/\partial z^L.$ And the matrices $\partial a^{L-1}/\partial w^{L-1}$ and $\partial a^{L-1}/\partial b^{L-1}$ have the same functional form as before,
$$
\frac{\partial\, a^{L-1}}{\partial\, w^{L-1}} =\sigma^\prime (z^{L-1}) a^{L-2}\quad\text{and}\quad \frac{\partial\, a^{L-1}}{\partial\, b^{L-1}}=\sigma^\prime(z^{L-1}),
$$ 
so most of the work has already been done. Even better, the remaining derivatives of the $D\times H$ matrix $\partial z^L/\partial\, a^{L-1}$ are easy to calculate from the definition. In fact,

$$
\frac{\partial\, z^L}{\partial\,a^{L-1}} = \left(
\begin{matrix}
w^L_{11} & w^L_{12} & \dots & w^L_{1H}\\
w^L_{21} & w^L_{22} & \dots & w^L_{2H}\\
\vdots & \vdots &\ddots &\vdots\\
w^L_{D1} & w^L_{D1} & \dots & w^L_{DH}
\end{matrix}
\right) =w^L.
$$

\mpar
So, we're on a roll! Given that we've already calculated all the partial derivatives with respect to weights and biases from layer $l+1$ forward to layer $L,$ we can calculate the partial derivatives with respect to parameters in the previous layer upstream in terms of derivatives we've already calculated plus the matrices,
$$
\frac{\partial a^{l}}{\partial w^{l}} =\sigma^\prime (z^{l}) a^{l-1}\quad\quad \frac{\partial a^{l}}{\partial b^{l}}=\sigma^\prime(z^{l}) \quad\quad \frac{\partial\, z^{l}}{\partial\,a^{l-1}} = w^{l}.
$$

\mpar
Let's make these calculations explicit by finding recursive formulas for $\partial a^L/\partial w^l$ and $\partial a^L/\partial b^l$ for all $1\leq l\leq L.$ We have,

\begin{align*}
\frac{\partial a^L}{\partial w^l} & =\left(\frac{\partial a^L}{\partial z^L}\,\frac{\partial z^L}{\partial a^{L-1}}\right)\,\frac{\partial a^{L-1}}{\partial w^l}\\
& = \left(\frac{\partial a^L}{\partial z^L}\,\frac{\partial z^L}{\partial a^{L-1}}\right)\,\left(\frac{\partial a^{L-1}}{\partial z^{L-1}}\,\frac{\partial z^{L-1}}{\partial a^{L-2}}\right)\,\frac{\partial a^{L-2}}{\partial w^l}\\
& = \quad \cdots \\
& = \left(\frac{\partial a^L}{\partial z^L}\,\frac{\partial z^L}{\partial a^{L-1}}\right)\,\left(\frac{\partial a^{L-1}}{\partial z^{L-1}}\,\frac{\partial z^{L-1}}{\partial a^{L-2}}\right)\,\cdots \left(\frac{\partial a^{l+1}}{\partial z^{l+1}}\,\frac{\partial z^{l+1}}{\partial a^l}\right) \frac{\partial a^l}{\partial w^l}.\\
\end{align*}
Similarly,
\begin{align*}
\frac{\partial a^L}{\partial b^l} & =\left(\frac{\partial a^L}{\partial z^L}\,\frac{\partial z^L}{\partial a^{L-1}}\right)\,\frac{\partial a^{L-1}}{\partial b^l}\\
& = \left(\frac{\partial a^L}{\partial z^L}\,\frac{\partial z^L}{\partial a^{L-1}}\right)\,\left(\frac{\partial a^{L-1}}{\partial z^{L-1}}\,\frac{\partial z^{L-1}}{\partial a^{L-2}}\right)\,\frac{\partial a^{L-2}}{\partial b^l}\\
& = \quad \cdots \\
& = \left(\frac{\partial a^L}{\partial z^L}\,\frac{\partial z^L}{\partial a^{L-1}}\right)\,\left(\frac{\partial a^{L-1}}{\partial z^{L-1}}\,\frac{\partial z^{L-1}}{\partial a^{L-2}}\right)\,\cdots \left(\frac{\partial a^{l+1}}{\partial z^{l+1}}\,\frac{\partial z^{l+1}}{\partial a^l}\right) \frac{\partial a^l}{\partial b^l}.\\
\end{align*}
Notice the terms in parentheses all have the same functional form - only the indices change from step to step. In fact, it's easy to calculate the product in parenthesis and we find the $(j,k)$ entry of the product matrix is 
$$
\left(\frac{\partial a^l}{\partial z^l}\,\frac{\partial z^l}{\partial a^{l-1}}\right)_{j,k} = \sigma^\prime(z^l_j)w^l_{jk}.
$$

Also notice terms from previous steps get reused in subsequent steps. The iterative nature of calculations that start from the output layer and propagate backwards to the first hidden layer makes the cost of computing $\nabla\hat{y}(x)$ comparable to the cost of computing $\hat{y}(x).$ That's the take home message of backpropagation. 

\mpar
To summarize backpropagation,
\begin{itemize}
\item{The predicted output can be written as the composition of functions,}
$$
\hat{y}(x)=\left(\sigma\circ z^L\right)\circ\left(\sigma\circ z^{L-1}\right)\circ\cdots\circ\left(\sigma\circ z^1\right)(x),
$$
\item{This functional form makes it relatively easy to calculate partial derivatives of $\hat{y}(x)=a^L$ by repeated application of the chain rule}
\item{For each $2\leq l\leq L-1,$}
\begin{align*}
\frac{\partial a^l}{\partial z^l}\cdot\frac{\partial z^l}{\partial a^{l-1}} & = \left(
\begin{matrix}
\sigma^\prime(z^l_1)w^l_{11} & \sigma^\prime(z^l_2)w^l_{12} & \dots & \sigma^\prime(z^l_H)w^l_{1H}\\
\sigma^\prime(z^l_1)w^l_{21} & \sigma^\prime(z^l_2)w^l_{22} & \dots & \sigma^\prime(z^l_H)w^l_{2H}\\
\vdots & \vdots &\ddots &\vdots\\
\sigma^\prime(z^l_1)w^l_{H1} & \sigma^\prime(z^l_2)w^l_{H2} & \dots & \sigma^\prime(z^l_H)w^l_{HH}
\end{matrix}
\right)\\ \\
\frac{\partial a^l}{\partial w^l} & =\left(
\begin{matrix}
\sigma^\prime(z^l_1)\,a^{l-1}_1 & \sigma^\prime(z^l_1)\,a^{l-1}_2 & \dots & \sigma^\prime(z^l_1)\,a^{l-1}_H\\
\sigma^\prime(z^l_2)\,a^{l-1}_1 & \sigma^\prime(z^l_2)\,a^{l-1}_2& \dots & \sigma^\prime(z^l_2)\,a^{l-1}_H\\
\vdots & \vdots &\ddots &\vdots\\
\sigma^\prime(z^l_H)\,a^{l-1}_1 & \sigma^\prime(z^l_H)\,a^{L-1}_2& \dots & \sigma^\prime(z^l_H)\,a^{l-1}_H
\end{matrix}
\right)\\ \\
\frac{\partial a^l}{\partial b^l} & = \left(
\begin{matrix}
\sigma^\prime(z^l_1) & \sigma^\prime(z^l_2) & \dots & \sigma^\prime(z^l_H)\\
\sigma^\prime(z^l_1) & \sigma^\prime(z^l_2) & \dots & \sigma^\prime(z^l_H)\\
\vdots & \vdots &\ddots &\vdots\\
\sigma^\prime(z^l_1) & \sigma^\prime(z^l_2) & \dots & \sigma^\prime(z^l_H)
\end{matrix}
\right)
\end{align*}

\item{For $l=1, L$ these matrices have a similar form but different dimensions. They are not $H\times H$ dimensional matrices like those above but rather $H\times d$ dimensional matrices if $l=1$ and $D\times H$ dimensional matrices if $l=L.$}
\item{For all $1\leq l\leq L$,}

\begin{align*}
\frac{\partial a^L}{\partial w^l} & =\left(\frac{\partial a^L}{\partial z^L}\,\frac{\partial z^L}{\partial a^{L-1}}\right)\,\left(\frac{\partial a^{L-1}}{\partial z^{L-1}}\,\frac{\partial z^{L-1}}{\partial a^{L-2}}\right)\,\cdots \left(\frac{\partial a^{l+1}}{\partial z^{l+1}}\,\frac{\partial z^{l+1}}{\partial a^l}\right) \frac{\partial a^l}{\partial w^l}\\ \\
\frac{\partial a^L}{\partial b^l} & =\left(\frac{\partial a^L}{\partial z^L}\,\frac{\partial z^L}{\partial a^{L-1}}\right)\,\left(\frac{\partial a^{L-1}}{\partial z^{L-1}}\,\frac{\partial z^{L-1}}{\partial a^{L-2}}\right)\,\cdots \left(\frac{\partial a^{l+1}}{\partial z^{l+1}}\,\frac{\partial z^{l+1}}{\partial a^l}\right) \frac{\partial a^l}{\partial b^l}
\end{align*}
\end{itemize}

\mpar\mpar
\subsection{Mini-batch}
\textit{This subsection uses the method of gradient descent B(ix), Bernoulli random variables C(i), the law of large numbers C(ii),  and stochastic approximation to the gradient C(vi).}

\mpar\mpar
Recall the the predicted value $\hat{y}(x)=\hat{y}(x, W, B)$ is a joint function of the input $x$, weights $W$, and biases $B,$ and $\delta^2(x)= ||\hat{y}(x) -y(x)||^2$ is the squared length of the difference between the predicted and actual values. Define the mean squared error $\Delta^2$ by the formula,
$$
\Delta^2(W,B) =\frac{1}{N}\sum_{i=1}^N \delta^2(x_n)
$$
where $N=\texttt{NumTrainingData}.$ Since $\Delta^2$ is a function of the weights and biases, it makes sense to calculate partial derivatives like,
$$
\frac{\partial}{\partial \,w^l_{jk}} \Delta^2 (W,B) \quad\text{and}\quad \frac{\partial}{\partial \,b^l_j} \Delta^2(W,B),
$$
for all possible values of the indices $j, k, l.$ These partial derivatives are the coordinates of the gradient vector $\nabla \Delta(W,B)$ where $-\nabla \Delta (W,B)$ is the direction in which the mean square error decreases most rapidly when the weights are $W$ and the biases are $B.$

\mpar
True gradient descent is the following algorithm for updating the weights and biases,
$$
w^l_{jk}\to w^l_{jk} -\eta\times\frac{\partial}{\partial\, w^l_{jk}} \Delta^2(W,B) \quad\text{and}\quad b^l_j\to b^l_j-\eta\times\frac{\partial}{\partial \,b^l_j} \Delta^2(W,B),
$$
where $\eta$ is a small positive parameter called the \textit{learning rate}. This sets up the gradient descent iteration, 

\begin{itemize}
\item{$W_0$ and $B_0$ are randomly chosen weights and biases}
\item {Repeat until a pre-determined stop rule occurs:}
\begin{itemize}
\item{$\hat{y}_{n+1}(x)$ is the feedforward of $x$ using $W_n, B_n$ for all $x$ in the training set}
\item{$\Delta^2(W_n, B_n)$ is the corresponding mean squared error}
\item{$(W_{n+1}, B_{n+1}) =(W_n, B_n) -\eta\times\nabla\, \Delta^2(W_n, B_n)$}
\end{itemize}
\end{itemize}

\mpar\mpar
Theoretically speaking, the stop rule is when the length of the gradient of $\Delta^2(W_n, B_n)$ drops below a specified level. Practically the stop rule is when a specified number of steps have been completed.

\mpar
Notice that $\Delta^2$ is a sum of the terms $\delta^2(x_n)$ corresponding to the  contribution of the $n^{th}$ data point $(x_n, y_n)$ to the total mean square error. \textit{Stochastic gradient descent} is based on the idea of choosing a training data point $(x_r, y_r)$ at random, which is akin to throwing a die with $N$ sides, and updating weights and biases via this rule,
$$
w^l_{jk}\to w^l_{jk} -\eta\times\frac{\partial}{\partial\, w^l_{jk}} \delta^2(x_r)\quad\text{and}\quad b^l_j\to b^l_j-\eta\times\frac{\partial}{\partial \,b^l_j} \delta^2(x_r).
$$ 

Note that the derivatives of $\delta^2(x_r)$ can be calculated efficiently via backpropagation. 

\mpar
Obviously, since $N$ is large, the computational cost of a single step of stochastic gradient descent is dramatically less than the cost a single step of true gradient descent, but the direction chosen by this method, $\nabla\delta(x_r)$ is random and quite unlikely to be the direction of steepest descent. However, since the data point $x_r$ was chosen uniformly at random, it is true that the average of $\nabla\delta(x_r)$ is,
$$
\nabla\, \Delta^2 = E[\,\nabla\delta^2(x_r)\,],
$$
so $-\nabla\delta^2(x_r)$ points in the direction of steepest descent \textit{on average}. To repeat: $\nabla\delta^2(x_r)$ is random and $\nabla\, \Delta^2$ is its mean value. But the variance of $\nabla\delta^2(x_r),$ being the sum of $N$ independent random variables, is proportional to $N$ which implies the signal-to-noise ratio will be pretty small, basically proportional to $1/\sqrt{N}.$  So, there is a lot of noise in stochastic gradient descent and not much signal. It's the classic tradeoff of cost versus accuracy.

\mpar
Mini-batch gradient descent is a compromise between stochastic and true gradient descent. Instead of choosing a single data point at random, let's choose a random subset of $R$ data points $(x_{r_1}, y_{r_1}), (x_{r_2}, y_{r_2}), \cdots, (x_{r_R}, y_{r_R})$
where $1< R< N$ is a suitably chosen number. In principle we can do this by rolling an $\binom{N}{R}$ sided die having one face for each subset of size $R.$ Then we form the mean squared error for this subset,
$$
\Delta^2_{MB}(W,B) = \frac{1}{R}\sum_{k=1}^R \delta^2(x_{r_k})
$$
and use the algorithm,
$$
w^l_{jk}\to w^l_{jk} -\eta\times\frac{\partial}{\partial\, w^l_{jk}} \Delta^2_{MB}(W,b)\quad\text{and}\quad b^l_j\to b^l_j-\eta\times\frac{\partial}{\partial \,b^l_j} \Delta^2_{MB}(W,b).
$$
The computational cost goes up again with mini-batch but the update direction, $-\nabla \Delta^2_{MB},$ may be appreciably closer to the direction of steepest descent, $-\nabla \Delta^2,$ than the update direction for stochastic gradient descent. Specifically, 
$$
\nabla \Delta^2=E[\,\nabla \Delta^2_{MB}\,]
$$ 
as before but the variance of $\nabla\Delta^2_{MB}$ deceases with increasing $R$  or, equivalently, the signal-to-noise ratio is increasing with increasing $R.$ That is just to say that larger values for $R$ mean higher computational costs but they also mean an update direction that is statistically closer to the optimal direction. Once again, it's a tradeoff. But mini-batch descent, for an appropriately chosen value of $R$ compared to $N,$ has become a standard choice. Notice that the case $R=1$ is stochastic gradient descent and the case $R=N$ is classic gradient descent, so mini-batch interpolates between the two extremes.

\mpar
In practice mini-batch is done in a completely different way. The issue is that there are multiple steps in classic gradient descent as applied to $\Delta^2_{MB}$ for one randomly chosen subset of $R$ training data points.  But for the law of large numbers to help, we need to combine multiple, independent choices of random subsets of $R$ data points. How do we take both these things into account in a computationally efficient way?

\mpar
Here's the trick. Randomly shuffle the $N$ training data points, which means apply a randomly chosen permutation $\pi$ of the set $\{1, 2, \,\dots,\, N\}.$ (This just means a uniformly distributed random variable on a set of $N!$ points). Then the first mini-batch $MB_1$ corresponds to the data points with indices $\pi_1, \pi_2,\,\dots,\, \pi_R,$ the second mini-batch $MB_2$ corresponds to the data points with indices $\pi_{R+1}, \pi_{R+2},\,\dots,\, \pi_{2R},$ etc. Note that the last mini-batch may be smaller in size than $R.$ In fact if $N=kR + l$ where $0\leq l<R$ then there are $k$ mini-batches of size $R$ and the remaining mini-batch $MB_{k+1}$ has size $l$. 

\mpar
We say an \textit{epoch} of training corresponds to one gradient descent update of the weights and biases using $\nabla \Delta_{MB_1}$ followed a second gradient descent update using $\nabla \Delta_{MB_2},$ etc., until the $k^{th}$ gradient descent update is made using $\nabla \Delta_{MB_k}.$ There is more than one way to deal with the remainder mini-batch $MB_{k+1},$ if there happens to be one, but for the sake of simplicity we're just going to ignore it. Experience shows we should reshuffle the data using an independent random permutation and do $k$ more gradient descent steps, one for each new mini-batch. That's a second epoch and there can be additional epochs, as needed. The mini-batch method requires use to choose a number $T$ of epochs after which network training is declared to be over.

\mpar
To summarize, mini-batch is a nested pair of loops. The outer loop consists of $T$ sweeps through the training data. Each sweep, or epoch, consists of randomly shuffling the training data set and  breaking it up into consecutive blocks of size $R.$ The $j^{th}$ mini-batch of the $t^{th}$ epoch is denoted $MB_{t,j}$ and a gradient descent update of the weights and biases is made using $\nabla\,\Delta_{MB_{t,j}}$ for $j=1, 2,\,\dots,\, k.$ 

\mpar
We need some notation to keep track of the weights and biases used in each update. Specifically, $W_0$ and $B_0$ are randomly initialized and serve as the weights and biases of the first update on the first mini-batch in the first epoch. The weight and biases of the $j^{th}$ update in the $t^{th}$ epoch are denoted $W_{t,j},\, B_{t,j}$ where $1\leq t\leq T, \, 1\leq j\leq k$ and $N=kR+l.$ Then mini-batch descent consists of the following steps,

\begin{itemize}
\item{\texttt{Eta, MinibatchSize}, and \texttt{Epochs} are adjustable parameters}
\item{$\eta$ = \texttt{Eta} is the learning rate}
\item{$R=\texttt{MinibatchSize}$ is the number of training data points in each mini-batch}
\item{$T=\texttt{Epochs}$ is the number of times the algorithm passes through the training set}
\item{$(W_{1,1}, B_{1,1}) = (W_0, B_0),\,t=1,$ and $j=1$}
\item{While $1\leq t < T,$ do}
\begin{itemize}\item{Shuffle training data and split it into $k$ consecutive blocks of size $R,$ where $N=kR+l$ and $0\leq l <R.$}
\item{While $1\leq j< k,$ do}
\item{For all $(x,y)\in MB_{t,j},\, \hat{y}(x)$ is the feedforward of $x$ using $(W_{t,j}, B_{t,j})$}
\item{$\Delta_{MB_{t,j}}(W_{t,j}, B_{t,j})$ is the corresponding mean squared error}
\item{$(W_{t,j+1}, B_{t,j+1}) =(W_{t,j}, B_{t,j}) -\eta\times\nabla\, \Delta(W_{t,j}, B_{t,j})$}
\item{$(W_{t+1,1}, B_{t+1,1})=(W_{t,k}, B_{t,k})$ and $j=1$}
\end{itemize}
\item{At the end of the training period the weights and biases are $(W_{T,k}, B_{T,k})$}
\end{itemize}

\mpar
The remaining question is: how do you choose $\eta,\, R$ and $T$? This question doesn't have an easy answer except, perhaps, \textit{it depends.} It depends on the size of the network and the amount of data, that much is pretty clear. But exactly - that is to say quantitatively  - how it depends is actually a current research question. At present it's a practitioner's art that is guided partly by theory and partly by practice. At any rate it lies well outside our meager expertise so we'll set the question aside, except for a few remarks in the next subsection.

\mpar\mpar
\subsection{Evaluation}
Implicitly, we've made a choice that our network is a \textit{classifier,} meaning its job is to choose among several competing alternatives. When we choose \texttt{OutputNeuron} $=D$ we are saying the neurons in the output layer correspond to $D$ possible  alternatives. For example, in our cartoon network $D=3$ corresponds to the outcomes \textit{Dog, Bear} and \textit{Human}. 

\mpar
Each output neuron contains an activation, which is a number between 0 and 1. One way to interpret these activations is, 
$$
\textit{The larger the activation, the liklier the classification}.
$$

This suggests we adopt the rule that the network's classification is the one with the highest activation. So, when the training's done we can feed the training data back into the network and see how well it's classified the data it was trained on. This step is called \textit{self-evaluation}. It has the form, 
$$
\textit{Classification $j$ was correct $X_j$\% of the time, for $1\leq j\leq D$,}
$$

and is an important tool in optimizing training. 

\mpar
Here's why. Among other things, training the network depends on choosing the number $T$ of epochs, a mini-batch size $R,$  and  the learning rate $\eta$. In the absence of firm guidance we could choose reasonable values for $T$ and $\eta$, train the network, then check the self-evaluation. If it's not good then maybe we can consider using more epochs and a different learning rate. Let's focus on that. 

\mpar
If you think about it for a minute, a relatively large learning rate makes relatively large changes in the weights and biases. This may prevent the mini-batch descent algorithm from converging to an approximate minimum of the mean squared error because the weights and biases are  wildly jumping all over the place. On the other hand, a relatively small learning rate may be leading us appropriately toward an approximate minimum of the mean squared error, but doing it so slowly that we need to increase the number of epochs which drives up the computational cost. 

\mpar
Turning this around, if the classification is not good then the culprit might be $\eta$ is too large, so reduce $\eta$ and see what happens. If the self-evaluation improves then you may be on the right track. If it doesn't improve then maybe the culptit is the mini-batch size $R$ is too small or there aren't enough epochs. So, increase $R$ and/or $T$ and see if that helps. It's a process.

\mpar
Setting aside the art and science of training, let's suppose the self-evaluation percentages are pretty good -  let's say all of them are over 90\%. That means the network correctly classifies the data it was trained on at least 9 times out of 10. OK, that's good as far as it goes. But what about data it hasn't seen yet? That's where the evaluation data set comes in. Feed those data into the trained network and check the evaluation percentages. What do you expect to see?

\mpar
It would be miraculous if the percentages went up, right? That's too much to expect. But what if they don't go down very much, say all of them over 80\%? Getting the classification right at least 4 times out of 5? Outstanding.\mpar
To summarize evaluation,

\begin{itemize}
\item We think of the neural network as a classifier, meaning a method to sort input data into one of $D$ possible alternatives
\item The coordinate of the predicted output vector $\hat{y}(x)$ that has the largest activation is network's classification of input $x$
\item The fraction, expresses as a percentage, of training data that are correctly classified is the \textit{self-evaluation}
\item The fraction, expresses as a percentage, of evaluation data that are correctly classified is the \textit{evaluation}
\end{itemize}
\mpar\mpar

\section{Pseudocode}
It's convenient to break the pseudocode into three pieces: Initialization, Training, and Evaluation. Initialization contains the Initialization step discussed above as well as methods that will be used repeatedly in other pieces of pseudocode. Evaluation consists only of the Evaluation step. Effectively we're combining the Mini-batch, Feedforward, and Backpropagation steps into the Training piece. It develops that Feedforward and Backpropagation are essential methods used to implement Mini-batch.  

\mpar
Note that we have not tried to optimize the pseudocode for computational efficiency but rather written it instead for ease of understanding. This makes our pseudocode longer, clumsier, and more pedantic than Nielsen's Python code. While this may frustrate or even irritate a knowledgeable reader, we do so in the hope that our pseudocode more transparent to the uninitiated.

\mpar
// \textbf{Initialization.} \textit{This piece of pseudocode asks the user for input for numerical values of certain variables. It also defines methods that will be used repeatedly.}\\

\texttt{InputNeuron} = \textit{*user input*\quad//Number of neurons in the input layer.}\\
\texttt{OutputNeuron} = \textit{*user input*\quad//Number of neurons in the output layer.}\\
\texttt{HiddenLayer} = \textit{*user input*\quad//Number of hidden layers.}\\
\texttt{HiddenNeuron} = \textit{*user input*\quad//Number of neurons in each hidden layer.}\\

$d=\texttt{InputNeuron}$\\
$D=\texttt{OutputNeuron}$\\
$L=\texttt{HiddenLayer} +1$\\
$H=\texttt{HiddenNeuron}$\\

//\textit{This section of the pseudocode defines three methods; that is, self-contained pieces of pseudocode that execute a specific task. We will use standard mathematical symbols for sums, inequalities, and familiar functions.} \\

\texttt{Length} \textit{accepts a vector as input and returns square root of the the sum of the squares of the coordinates.}\\

\texttt{FF} \textit{accepts an input vector, a set of weights, and a set of biases as input and returns the activations of every neuron, as calculated by the feedforward mechanism. The input vector $x$ is the vector $a^0$ of activations in the input layer and the output vector $\hat{y}(x)$ is the vector $a^L$ of activations in the output layer.}\\
 
\texttt{BP} \textit{accepts activations of all neurons, along with the weights and biases that were used to calculate them, and returns the gradient of $a^L$ with respect the the weights and biases in every layer.}\\
 
//\textit{The method} $\texttt{Length} = \texttt{Length}(x).$ \textit{The input variable is a vector understood as a list.}\\
$c=\textit{length-of-list}(x)$\\
$\texttt{LengthSquared}(x)=\sum_{i=1}^c x_i^2$\\
$\texttt{Length}(x) =\sqrt{\texttt{LengthSquared}(x)}$

\mpar\mpar
//\textit{The method} $\texttt{FF}=\texttt{FF}(W,B,x).$ \textit{Input variables are weights $W=w[l][j][k],$ biases $B=b[l][j],$ and a $d$-dimensional input vector $x$ whose coordinates are the activations in the input layer neurons: $a[0][j]=x_j.$ Output variables are the score functions $\texttt{FF}(W,B,x)[1][l][k]=z[l][j]$ and activations $\texttt{FF}(W,B,x)[0][l][k]=a[l][k]$ of all remaining neurons in the network. In particular, $\hat{y}(x)_k=\texttt{FF}(W, B, x)[0][L][k]$ for all $1\leq k\leq D.$ Sometimes it's convenient to write $Z=\texttt{FF}(W,B,x)[1][l][k]$ and $A=\texttt{FF}(W,B,x)[0][l][k]$ to denote the set of activations and score functions.}\\

//\textit{In principle the sigmoid function $\sigma$ is a user choice. For simplicity we've made a specific choice.}\\
$\texttt{Sigma}(x)= 1/(1+e^{-x})$\\
$\texttt{SigmaPrime}(x) = 1/(e^{-x/2}+e^{x/2})^2$\\

//\textit{Calculate activations in terms of the score function z.}\\
//\textit{In the first hidden layer:}\\
\textit{While $j\leq H$}\\
$\phantom{mm} z[1][j] =\sum_{k=1}^d w[1][j][k]\cdot x_k + b[1][j]$\\
\phantom{mm} $a[1][j]=\texttt{Sigma}(z[1][j])$\\
\textit{end}\\

\mpar
//\textit{In the remaining hidden layers:}\\
\textit{While $l\geq 2$ and $l\leq L-1$}\\
\phantom{mm} \textit{While $j\leq H$}\\
\phantom{mmmm} $z[l][j]=\sum_{k=1}^H w[L][j][k]\cdot a[l-1][k] + b[l][j]$\\
\phantom{mmmm} $a[l][j]=\texttt{Sigma}(z[l]j])$\\
\phantom{mm} \textit{end}\\
\textit{end}\\

//\textit{In the output layer:}\\
\textit{While $j\leq D$}\\
\phantom{mm} $z[L][j] =\sum_{k=1}^H w[L][j][k]\cdot a[L-1][k] + b[L][j]$\\
\phantom{mm} $a[L][j]=\texttt{Sigma}(z[L][j])$\\
\textit{end}\\
//\textit{This concludes the method \texttt{FF}.}\\ 

//\textit{The method} $\texttt{BP} = \texttt{BP}(W,B,Z, A).$ \textit{Input variables are weights $W=w[l][j][k],$ biases $B=b[l][j],$ and the score functions  $Z=\texttt{FF}(W,B,x)[1][l][k]$ and activations $A=\texttt{FF}(W,B,x)[0][l][k]$ that were calculated by \texttt{FF} using $W$ and $B.$ Output variables are the partial derivatives of components of the output vector $a[L][j]$ with respect to weights $w[l][j][k]$ and biases $b[l][j].$}\\

//\textit{We write $ALW[l][j][k] =\partial\, a^L/\partial\, w^l_{jk}$ and $ALB[l][j]=\partial\, a^L/\, \partial\, b^l_j$ for the partial derivatives of activations with respect to the weights and biases. Then,
\begin{align*}
\texttt{BP}(W,B, Z, A)[0][l][j][k] & = ALW[l][j][k]\\
\texttt{BP}(W,B, Z, A)[1][l][j]&  = ALB[l][j]
\end{align*}
where $Z=\texttt{FF}(W,B,x)[1][l][k]$ and $A=\texttt{FF}(W,B,x)[0][l][k]$}\\

//\textit{Define matrices which are building blocks in the calculation the partial derivatives of $a^L.$ Specifically, $AA[l]=\partial\, a^l/\partial\, z^l\cdot\partial\, z^l/\partial\, a^{l-1}, AW[l]=\partial\, a^l/\,\partial\, w^l,$ and $AB[l] = \partial\, a^l/ \partial\, b^l.$ See the summary of Section 3(iii) \textbf{Backpropagation} to review the mathematical definitions of these matrices.}\\

//\textit{Defining auxiliary matrices.}\\
//\textit{In the first hidden layer:}\\
\textit{While $j\leq H$}\\
\phantom{mm} \textit{While $k\leq d$}\\
\phantom{mmmm} $AA[1][j][k] =\texttt{SigmaPrime}(z[1][j])\cdot w[1][j][k]$\\
\phantom{mmmm} $AW[1][j][k] = \texttt{SigmaPrime}(z[1][j])\, a[0][k]$\\
\phantom{mmmm} $AB[1][j][k] =\texttt{SigmaPrime}(z[1][j])$\\
\phantom{mm}\textit{end}\\
\textit{end}\\

//\textit{In the remaining hidden layers:}\\
\textit{While $2\leq l$ and $l\leq L-1$}\\
\phantom{mm}\textit{While $j\leq H$}\\
\phantom{mmmm} \textit{While $k\leq H$}\\
\phantom{mmmmmm} $AA[l][j][k] =\texttt{SigmaPrime}(z[l][j])\cdot w[l][j][k]$\\
\phantom{mmmmmm} $AW[l][j][k] = \texttt{SigmaPrime}(z[l][j])\, a[l-1][k]$\\
\phantom{mmmmmm} $AB[l][j][k] =\texttt{SigmaPrime}(z[l][j])$\\
\phantom{mmmm}\textit{end}\\
\phantom{mm}\textit{end}\\
\textit{end}\\

//\textit{In the output layer:}\\
\textit{While $j\leq D$}\\
\phantom{mm} \textit{While $k\leq H$}\\
\phantom{mmmm} $AA[L][j][k] =\texttt{SigmaPrime}(z[L][j])\cdot w[L][j][k]$\\
\phantom{mmmm} $AW[L][j][k] = \texttt{SigmaPrime}(z[L][j])\, a[L-1][k]$\\
\phantom{mmmm} $AB[1][j][k] =\texttt{SigmaPrime}(z[L][j])$\\
\phantom{mm}\textit{end}\\
\textit{end}\\

//\textit{Calculate partial derivatives of $a^L$ iteratively starting with the weights and biases of the output layer, or $L^{th}$ layer.}\\

//\textit{Note that partial derivatives are calculated in terms of descending partial products,
\begin{align*} 
DPP[L] & = AA[L]\\
DPP[L-1] &=DPP[L]\cdot AA[L-1]\\
& = AA[L]\cdot AA[L-1]\\
DPP[L-2]& =DPP[L-1]\cdot AA[L-2]\\
& = AA[L]\cdot AA[L-1]\cdot AA[L-2],
\end{align*}
etc. Again, see the summary of Section 3(iii) \textbf{Backpropagation} to review the mathematical background.}\\

//\textit{Calculate auxiliary matrices DPP.}\\
//\textit{In the output layer:}\\
\textit{While $j\leq D$}\\
\phantom{mm}\textit{While $k\leq H$}\\
\phantom{mmmm}$DPP[L][j][k]=AA[L][j][k] =\texttt{SigmaPrime}(z[L][j])\cdot w[L][j][k]$\\
\phantom{mm}\textit{end}\\
\textit{end}\\

//\textit{From the last hidden layer to second hidden layer:}\\
\textit{While $l\geq 2$ and $l\leq L-1$}\\
\phantom{mm}\textit{While $j\leq H$}\\
\phantom{mmmm}\textit{While $k\leq H$}\\
\phantom{mmmmmm}$DPP[L-l][j][k]=\sum_{i=1}^H DPP[L-l+1][j][i]\cdot AA[L-l][i][k]$\\
\phantom{mmmm}\textit{end}\\
\phantom{mm}\textit{end}\\
\textit{end}\\

//\textit{Calculate partial derivatives of $a^L$ in the output layer.}\\
\textit{While $j\leq D$}\\
\phantom{mm} $ALB[L][j] = \texttt{SigmaPrime}(z[L][j])$\\
\phantom{mm}\textit{While $k\leq H$}\\
\phantom{mmmmmm}$ALW[L][j][k] = \texttt{SigmaPrime}(z[L][j]) \cdot a[L-1][k]$\\
\phantom{mm}\textit{end}\\
\textit{end}\\

//\textit{Calculate partial derivatives of $a^L$ from the last hidden layer to the second hidden layer.}\\
\textit{While $l\geq 2$ and $l\leq L-1$}\\
\phantom{mm}\textit{While $j\leq H$}\\
\phantom{mmmm}$ALB[l][j] =\sum_{i=1}^H DPP[L-l+1][j][i]\cdot AB[l][i][k]$\\
\phantom{mmmm}\textit{While $ k\leq H$}\\
\phantom{mmmmmm}$ALW[l][j][k] = \sum_{i=1}^H DPP[L-l+1][j][i]\cdot AW[l][i][k]$\\
\phantom{mmmm}\textit{end}\\
\phantom{mm}\textit{end}\\
\textit{end}\\

//\textit{Calculate partial derivatives of $a^L$ in the first hidden layer.}\\
\textit{While $j\leq H$}\\
\phantom{mm} $ALB[1][j] = \texttt{SigmaPrime}(z[1][j])$\\
\phantom{mm}\textit{While $k\leq d$}\\
\phantom{mmmmmm}$ALW[1][j][k] =\sum_{i=1}^d DDP[2][j][i]\cdot AW[1][i][k]$\\
\phantom{mm}\textit{end}\\
\textit{end}\\
//\textit{This concludes the method \texttt{BP}.}\\

//\textit{\textbf{Training.} This section of pseudocode executes the mini-batch step. This is done in\, $T$ epochs. In the $t^{th}$ epoch the training data are randomly shuffled and partitioned into $k$ disjoint subsets of size $R.$ The $n^{th}$ random subset of the $t^{th}$ epoch defines a mini-batch $MB_{t,n}$ of training data and a $\Delta^2_{MB_{t,n}}$ denotes the corresponding partial mean squared error.  A single gradient descent step is performed using this function to update the weights and biases. The updated weights and biases are denoted $W[t][n]=w[t][n][l][j][k]$ and $B[t][n]=b[t][n][l][j],$ respectively. These calculations required multiple calls to the methods \texttt{Length}, \texttt{FF}, and \texttt{BP} that we formulate as a new method \texttt{MB}.}\\

\texttt{TrainingData} =\textit{*user input*}\\
\texttt{EvaluationData} =\textit{*user input*}\\
//\textit{These are lists of pairs $(x,y)$ where $x$ is a $d$-dimensional vector and $y$ is a $D$-dimensional vector.}\\\\
\texttt{NTD} = \textit{length-of-list}(\texttt{TrainingData})\\
\texttt{NED} =\textit{ length-of-list}(\texttt{EvaluationData})\\
//\textit{These are the numbers of pairs in each list.}\\

\texttt{Eta} = \textit{*user input*\quad//This is the learning rate.}\\
$R=$ \textit{*user input*\quad//This is the size of each mini-batch. $R$ must be less than \texttt{NTD}.}\\
$k=floor(\texttt{NTD}\div R)\quad$//\textit{This is the number of mini-batches in each epoch. Recall $floor(x)$ is the largest integer less than or equal to $x.$}\\
$T=$ \textit{*user input*}\quad//\textit{This is the number of epochs.}\\

//\textit{The method $\texttt{MB} = \texttt{MB}(W,B, MB).$ The input variables are a set of weights $W$, biases $B,$ and a mini-batch MB of R training data pairs $(x[1], y[1]), (x[2], y[2]), \cdots, (x[R], y[R]).$ The output variables are updated weights $W_*=\texttt{MB}(W,B, MB)[1]$ and updated biases $B_*=\texttt{MB}(W,B, MB)[0].$ In this formulation, 
\begin{align*}
W[t][n] &=\texttt{MB}(W[t][n-1], B[t][n-1], MB_{t, n-1})[1]\\
B[t][n] & =\texttt{MB}(W[t][n-1], B[t][n-1], MB_{t, n-1})[0].
\end{align*}}

//\textit{Calculate score functions $Z$ and activations $A$ for the R input vectors of the mini-batch $MB$. Note that in the loop below, $Z[m]$ is a matrix with components $Z[m]=Z[m][l][j]$ and $A[m]$ is a matrix with components $A[m][l][j]$ for every $1\leq m\leq R.$}\\

\textit{While $m\leq R$}\\
\phantom{mm}$Z[m]=\texttt{FF}(W,B, x[m])[1]$\\
\phantom{mm}$A[m]=\texttt{FF}(W,B, x[m])[0]$\\
\textit{end}\\

//\textit{Calculate the gradient of $\Delta^2_{MB}.$ This uses the fact, for example, that 
$$
\partial\, ||a^L -y||^2/\partial\, w^l_{jk} = 2 ||a^L -y||\times \partial\, a^L/\partial\, w^l_{jk}.
$$}

\textit{Note that for $1\leq m\leq R,$ $a[m][L]$ is the vector of activations in the output layer. We denote gradients with respect to biases by $Grad[0]$ and gradients with respect to weights by $Grad[1].$}\\

//\textit{Gradients in the first hidden layer.}\\
\textit{While $k\leq H$}\\
\phantom{mm}$Grad[0][1][k]=\frac{2}{R}\sum_{m=1}^R \texttt{Length}(a[m][L] -y[m])\cdot \texttt{BP}(W, B, Z, A, x[m])[0][1][k]$\\
\phantom{mmmm}\textit{While $j\leq d$}\\
\phantom{mmmm}$Grad[1][1][j][k]=\frac{2}{R}\sum_{m=1}^R \texttt{Length}(a[m][L] -y[m])\cdot \texttt{BP}(W, B, Z, A, x[m])[0][1][j][k]$\\
\phantom{mm}\textit{end}\\
\textit{end}\\

//\textit{Gradients in the remaining hidden layers.}\\
\textit{While $l\leq L-1$}\\
\phantom{mm}\textit{While $k\leq H$}\\
\phantom{mmmm}$Grad[0][l][j]=\frac{2}{R}\sum_{m=1}^R \texttt{Length}(a[m][L] -y[1][m])\cdot \texttt{BP}(W,B,Z,A, x[m])[l][k]$\\
\phantom{mmmm}\textit{While $j\leq H$}\\
\phantom{mmmmmm}$Grad[1][l][j][k] = \frac{2}{R}\sum_{m=1}^R \texttt{Length}([m][L] -y[1][m])\cdot \texttt{BP}(W,B,Z,A, x[m])[l][j][k]$\\
\phantom{mmmm}\textit{end}\\
\phantom{mm}\textit{end}\\
\textit{end}\\

//\textit{Gradients in the output layer.}\\
\textit{While $k\leq D$}\\
\phantom{mm}$Grad[0][L][k]=\frac{2}{R}\sum_{m=1}^R \texttt{Length}(a[m][L] -y[m])\cdot \texttt{BP}(W, B, Z, A, x[m])[0][L][k]$\\
\phantom{mmmm}\textit{While $j\leq H$}\\
\phantom{mmmm}$Grad[1][L][j][k]=\frac{2}{R}\sum_{m=1}^R \texttt{Length}(a[m][L] -y[m])\cdot \texttt{BP}(W, B, Z, A, x[m])[0][L][j][k]$\\
\phantom{mm}\textit{end}\\
\textit{end}\\

//\textit{Update the weights and biases from $W, B$ to $W_*, B_*.$}\\
//\textit{In the initial layer.}\\
\textit{While $k\leq H$}\\
\phantom{mm}$b_*[1][k] =b[1][k] -\texttt{Eta}\times Grad[0][1][k]$\\
\phantom{mm}\textit{While $\leq d$}\\
\phantom{mmmm}$w_*[1][j][k] = w[1][j][k] -\texttt{Eta}\times Grad[1][1][j][k]$\\
\phantom{mm}\textit{end}\\
\textit{end}\\

//\textit{In the hidden layers.}\\
\textit{While $l\leq L-1$}\\
\phantom{mm}\textit{While $k\leq H$}\\
\phantom{mmmm}$b_*[l][k] =b[l][k] -\texttt{Eta}\times Grad[0][l][k]$\\
\phantom{mmmm}\textit{While $j\leq H$}\\
\phantom{mmmmmm}$w_*[l][j][k] = w[l][j][k] -\texttt{Eta}\times Grad[1][l][j][k]$\\
\phantom{mmmm}\textit{end}\\
\phantom{mm}\textit{end}\\
\textit{end}\\

//\textit{In the output layer.}\\
\textit{While $k\leq D$}\\
\phantom{mm}$b_*[L][k] =b[L][k] -\texttt{Eta}\times Grad[0][L][k]$\\
\phantom{mm}\textit{While $\leq H$}\\
\phantom{mmmm}$w_*[L][j][k] = w[L][j][k] -\texttt{Eta}\times Grad[1][L][j][k]$\\
\phantom{mm}\textit{end}\\
\textit{end}\\
//\textit{This concludes the method \texttt{BP}.}\\

//\textit{Initialize the weights $W$ and biases $B.$ This entails repeated calls to a random number generator, a standard feature of programming languages like Python. We use \texttt{RNG} to denote a single call to the random number generator that returns a sample of a standard normal random variable. We assume repeated calls to \texttt{RNG} are statistically independent.} 

\mpar
//\textit{Initialize weights and biases from the input layer to the first hidden layer.}\\
\textit{While $k\leq \texttt{HiddenNeuron}$\\
\phantom{mm} $b[1][k] = \texttt{RNG}$\\
\phantom{mm}While $j\leq \texttt{InputNeuron}$\\
\phantom{mmmm} $w[1][j][k] =\texttt{RNG}$\\
\phantom{mm}\textit{end}\\
\textit{end}}

\mpar
//\textit{Initialize weights and biases between hidden layers.}\\
\textit{While $2\leq l$ and $l\leq L-1$\\
\phantom{mm}\textit{While $k\leq \texttt{HiddenNeuron}$\\
\phantom{mmmm} $b[l][j] = \texttt{RNG}$\\
\phantom{mmmm}While $j\leq \texttt{HiddenNeuron}$\\
\phantom{mmmmmmm} $w[l][j][k] =\texttt{RNG}$\\
\phantom{mmmm}\textit{end}\\
\phantom{mm}\textit{end}}\\
end}

\mpar
//\textit{Initialize weights and biases from the last hidden layer to the output layer.}\\
\textit{While $k\leq \texttt{OutputNeuron}$\\
\phantom{mm} $b[L][j] = \texttt{RNG}$\\
\phantom{mm}While $k\leq \texttt{HiddenNeuron}$\\
\phantom{mmmm} $w[L][j][k] =\texttt{RNG}$\\
\phantom{mm}\textit{end}\\
\textit{end}}\\

//\textit{A standard feature of programming languages is a command which randomly shuffles the entries in a list. If the list has $N$ entries the means choosing a random permutation $\sigma$ of the list or, equivalently, sampling a uniform random variable on a set of size $N!$ We use \texttt{Shuffle} to denote this feature and assume repeated calls to $\texttt{Shuffle}$ are statistically independent. So, if $x=(x_1, x_2,\cdots, x_N)$ is a list then $\texttt{Shuffle}(x) = (x_{\sigma(1)}, x_{\sigma(2)}, \cdots, x_{\sigma(N)})$ is a list of the elements of $x$ in random order. We write $\texttt{Shuffle}(x)_n = x_{\sigma(n)}$ to denote the $n^{th}$ element in the shuffled list.}\\

$W[1][1]=W$\\
$B[1][1] = B$\\
\textit{While $t\leq T$}\\
\phantom{mm}\textit{While $n\leq k$}\\
\phantom{mmmm}\textit{While $m\leq R$}\\
\phantom{mmmmmm}$MB[t][n][m]=\texttt{Shuffle}(\texttt{TrainingData})_{(n-1)R+m}$\\
\phantom{mmmm}\textit{end}\\
\phantom{mmmm}$W[t][n+1]=\texttt{MB}(W[t][n], B[t][n], MB[t][n])(1)$\\
\phantom{mmmm}$B[t][n+1]=\texttt{MB}(W[t][n], B[t][n], MB[t][n])(0)$\\
\phantom{mm}\textit{end}\\
\phantom{mm} $W[t+1][1]=W[t][k]$\\
\phantom{mm} $B[t+1][1]=B[t][k]$\\
\textit{end}\\

//\textit{\textbf{Evaluation.} The training pseudocode produces ``optimal'' weights $W_*=W[T+1][1]$ and biases $B_*=B[T+1][1].$ The evaluation pseudocode measures how well the neural network performs using them. We think of the neural network as a classifier, meaning that for each training datum $(x,y)$ the vector $y$ has exactly one coordinate equal to 1 and all the remaining coordinates equal 0. This signifies that fact that each vector $x$ represents one of $D$ mutually exclusive alternatives. Recall that $\texttt{FF}(W_*, B_*, x)[0][L]$ is the vector of activations in the output layer using the ``optimal'' weights and biases. Evaluation consists of counting the number of times the coordinate where the largest value of $\texttt{FF}(W_*, B_*, x)[0][L][j]\, 1\leq j\leq D$ occurs equals the coordinate of $y$ where 1 occurs.}\\

//\textit{The method \texttt{ArgMax}=\texttt{ArgMax}(x). The input variable is a vector thought of as a list. The output variable is an integer equal to the coordinate where the largest value of $x[n],\, 1\leq n\leq length-of-list(x)$ occurs. If the largest value of $x[n]$ occurs in multiple coordinates, \texttt{ArgMax} chooses the last coordinate where the maximum value occurs. Since it is unlikely that the largest value of $\texttt{FF}(W_*, B_*, x)[0][L]$ occurs in multiple coordinates, this artifact of \texttt{ArgMax} is unlikely to significantly affect the evaluation process.}\\

$Max(x) = x[1]$\\
$\texttt{ArgMax}(x)=1$\\
\textit{While $n\leq length-of-list(x)$}\\
\phantom{mm}\textit{if $x[n]> Max$}\\
\phantom{mmmm}\textit{then $Max=x[n], \texttt{ArgMax}=n$}\\
\phantom{mm}\textit{end}\\
\textit{end}\\

//\textit{Self evaluation \texttt{SuccessPercentage(TrainingData)}, meaning the percentage of successful classifications of training data $(x,y)\in\texttt{TrainingData}.$}\\
\textit{Count=0}\\
\textit{While $n\leq \texttt{NTD}$}\\
\phantom{mm}$\textit{if \texttt{Argmax}(\texttt{FF})}(W_*, B_*, x[n])=\texttt{ArgMax}(y[n])$\\
\phantom{mmmm}\textit{Count = Count+1}\\
\phantom{mm}\textit{end}\\
\textit{end}\\
\texttt{SuccessPercentage(TrainingData)}=$100\times \textit{Count}\div\texttt{NTD}$\\

//\textit{Evaluation \texttt{SuccessPercentage(TrainingData)}, meaning the percentage of successful classifications of evaluation data $(x,y)\in\texttt{EvaluationData}.$}\\
\textit{Count=0}\\
\textit{While $n\leq \texttt{NED}$}\\
\phantom{mm}$\textit{if \texttt{Argmax}(\texttt{FF})}(W_*, B_*, x[n])=\texttt{ArgMax}(y[n])$\\
\phantom{mmmm}\textit{Count = Count+1}\\
\phantom{mm}\textit{end}\\
\textit{end}\\
\texttt{SuccessPercentage(EvaluationData)}=$100\times \textit{Count}\div\texttt{NED}$\\

//\textit{This ends the neural network pseudocode.}\\

\mpar\mpar
\appendix 
\section{Linear Algebra}
Linear algebra is all about vectors and relationships between vectors. This appendix  focuses and cartesian and polar coordinates, length and scalar multiples of a vector, sum, difference and angle between two vectors, and linear transformation of vectors. For simplicity's sake, we'll start by talking about two dimensional vectors because they're easy to visualize.  At the end of the appendix, we'll show how it all works for three dimensions, four dimensions, or in fact any number of dimensions. This is necessary because in practice, neural networks require using vectors whose dimensions are thousands or even tens of thousands of dimensions.
\mpar
\subsection{Vectors: Length, Direction, and Coordinates}The first thing to consider is the idea of a vector. You can think of it as a direction and a positive number which we think of as a length, speed, or rate. For example, the velocity of a truck from south Jersey heading up I-95 toward Newark Liberty International Airport is a vector: the direction is North North East and the number is $65$ (that is $65$ miles per hour, if the driver is sticking to the legal limit). 

\mpar
Here's another example. If you're standing on a street corner and someone asks you for directions to a store, you're likely to say something like, 

\mpar
\textit{Just go up three blocks, turn right, then go over a block and a half and you'll see the store on the left side of the street.'}
\mpar

In math terms, you just provided the coordinates of the store relative to the spot where you're standing, and we'll write $(\text{over}, \text{up})=(1.5,3)$ as shorthand for these coordinates. We also use it as shorthand for the directed line segment - or the \textit{vector} - from $(0,0),$ the coordinates of the place you're standing, to $(1.5, 3),$ the coordinates of the store. (See Figure 8).

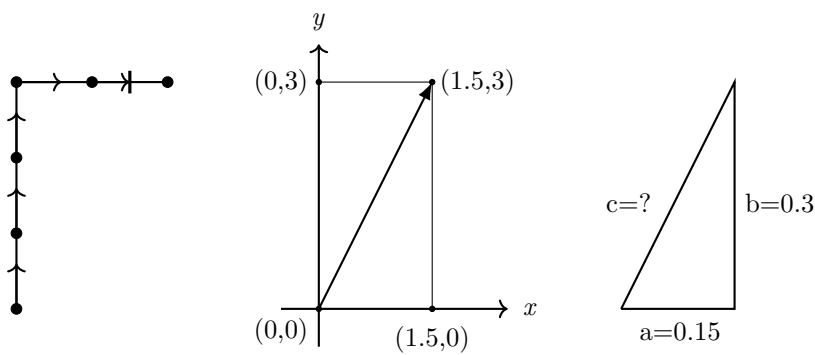
\begin{figure} [h]
\quad
\begin{tikzpicture}
\draw[fill=black] (-0.5,0) circle (2pt);
\draw[fill=black] (-0.5,1) circle (2pt);
\draw[fill=black] (-0.5,2) circle (2pt);
\draw[fill=black] (-0.5,3) circle (2pt);
\draw[fill=black] (0.5, 3) circle (2pt);
\draw[fill=black] (1.5,3) circle (2pt);

\draw[thick] (-0.5,0) -- (-0.5,3) -- (1.5, 3);
\draw[very thick] (1, 2.85) -- (1, 3.15);
\draw[thick,->] (-0.5,0) -- (-0.5,0.6);
\draw[thick,->] (-0.5,1) -- (-0.5,1.6);
\draw[thick,->] (-0.5,2) -- (-0.5,2.6);
\draw[thick,->] (-0.5,3) -- (0.1,3);
\draw[thick,->] (0.5, 3) -- (1,3);

\draw[thick,->] (3,0) -- (6, 0);
\draw[thick,->] (3.5,-0.5) -- (3.5, 3.5);

\node at (6.3,0) {\textit{x}};
\node at (3.5,3.8) {\textit{y}};

\draw[fill=black] (3.5, 3) circle (1pt);
\draw[fill=black] (5,0) circle (1pt);
\draw[fill=black] (5, 3) circle (1pt);
\draw[fill=black] (3.5, 0) circle (1pt);

\node at (3,3) {(0,3)};
\node at (5,-0.4) {(1.5,0)};
\node at (5.6,3) {(1.5,3)};
\node at (3,-0.3) {(0,0)};

\draw[thin] (3.5, 3) -- (5, 3);
\draw[thin] (5,0) -- (5, 3);
\draw[thick,-Latex] (3.5, 0) -- (5, 3);

\draw[thick] (7.5, 0) --(9,3) -- (9,0) -- (7.5,0); 

\node at (8.28,-0.3) {\text{a=0.15}};
\node at (9.6, 1.4) {\text{b=0.3}};
\node at (7.6, 1.4) {\text{c=?}};

\end{tikzpicture}

\caption{On the left is a cartoon of a city street with directions from where you're standing to the store (vertical bar). In the middle is the vector $(1.5, 3)$ written in terms of a rectangular $x-y$ grid. On the right is the corresponding right triangle with sidelengths labelled $a, b$ and $c.$ }
\end{figure}

\mpar
If you're on the street, you'll walk four and a half blocks to get to the store. If each block is a tenth of a mile long, that means you'll walk $0.45$ miles altogether. But the store is closer ``as the crow flies,'' meaning in a straight line, so how far away is it in that sense? You can see from Figure 1 that the vector $(1.5, 3)$ defines a right triangle with side lengths $a=0.15$ and $b=(0.3)$ (don't forget that city blocks are a tenth of a mile long!). If you remember Pythagoras' theorem, you'll know that the length of the hypotenuse of a right triangle is related to its side lengths by the formula $c^2=a^2+b^2.$ Let's use this formula to calculate the straight line distance to the store by plugging in numbers for $a$ and $b$ then taking the square root,
$$
c=\sqrt{a^2+b^2} = \sqrt{(0.15)^2 + (0.3)^2} =\sqrt{0.1125} \cong 0.335
$$

So, the store is almost half a mile on foot but only about a third of a mile as the crow flies. 
\mpar
In general, we can use a pair of numbers $(x,y)$ to represent a point in the plane which is $x$ units of distance from the origin of coordinates in the horizontal direction and $y$ units of distance from the origin in the vertical direction. We can also use the same pair $(x,y)$ to represent the vector corresponding to this point, defined as the directed line segment going from $(0,0)$ - the \textit{base point} - toward $(x,y)$ - the \textit{end point}.

\mpar
Let's call this vector $v,$ or sometimes $v=(x,y)$ when we want to emphasize the coordinates of the vector. The exact value of $x$ and $y$ don't really matter, you can calculate the length of the vector $v,$ which we denote $||v||,$ in the same way as before by observing that $v$ is the hypotenuse of a right triangle with side lengths $x$ and $y$. So, by Pythagoras' theorem, the length of $v$ is,
$$
||v||=\sqrt{x^2+y^2}.
$$
Think of it as a formula: if you know $x$ as a number and $y$ as a number then you can calculate the length of $v$ as a number. You can use a calculator, if you want, but the main thing is just to realize that it can be calculated, whenever you want.

\mpar
If you think about it for a second, it makes sense to talk about doubling a vector or cutting it in half, or multiplying it by any number you please, like $-17.$ For example, if $v=(x,y)$ then twice $v$ is $2v=(2x,2y)$, half $v$ is $\tfrac{1}{2}v = (0.5x, 0.5 y)$ and $-17 v=(-17x, -17y).$ From now on, if $\lambda$ is some number (like $\lambda = 2$ or $\tfrac{1}{2}$ or $-17$ or anything else) then $\lambda v$ is the vector with coordinates $(\lambda x, \lambda y).$

\mpar
A clever choice of $\lambda$ is the reciprocal of the length of $v$ itself: 
$$
\lambda =\frac{1}{||v||} = \frac{1}{\sqrt{x^2+y^2}}.
$$ 
It seems obvious that the vector $u=\lambda v$ has unit length, meaning $||u||=1.$ But there is a lot going on here so let's try to prove it. 

\mpar
First of all, what are the coordinates of $u$? Well, $u=\lambda v$ so, 
$$
u=(\lambda x, \lambda y) =(\frac{x}{\sqrt{x^2+y^2}}, \frac{y}{\sqrt{x^2+y^2}}).
$$
OK, so what is the length of $u$? Well, it's the square root of the sum of the squares of its coordinates. Let's write out $||u||^2$ instead of $||u||$, just to eliminate the pesky square root:
$$
||u||^2 =(\frac{x}{\sqrt{x^2+y^2}})^2 +(\frac{y}{\sqrt{x^2+y^2}})^2  = \frac{x^2}{x^2+y^2} + \frac{y^2}{x^2+y^2}=1.
$$
But if $||u||^2=1$ then $||u||=1,$ and this is what we wanted to prove. 

\begin{figure} [h]
\begin{tikzpicture}
\draw[thick,->] (5,-2) -- (5,2);
\draw[thick,->] (3,0) -- (7,0);

\draw (5,0) circle (1);
\node at (6.4, 0.3) {(1,0)};
\draw[fill=black] (6, 0) circle (1.5pt);

\draw[thin,->] (5,0) -- (5.7, 1.7);
\draw[thick,-Latex] (5,0) -- (5.38, 0.92);

\draw[thin,->] (5,0) -- (2.75, 0.825);
\draw[thick,-Latex] (5,0) -- (4.06, 0.365);

\draw[thin,->] (5,0) -- (5.75, -1.3);
\draw[thick,-Latex] (5,0) -- (5.502, -0.858);
\end{tikzpicture}
\caption{The unit circle is drawn in the $x-y$ plane along with three vectors of varying lengths shown in grey. Their directions are unit vectors which are shown in black. In this diagram, all the lengths are bigger than 1 but they don't have to be. The length of a vector can be any number bigger than or equal to 0.}
\end{figure}
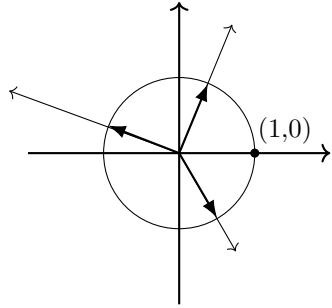

\mpar
The number $\lambda = ||v||$ is called the \textit{length} of $v$ and the vector $u=v/||v||$ is called the \textit{direction} of $v$. Both $\lambda$ and $u$ are functions of the coordinates $(x,y).$ 

\mpar
If you think about it for a minute, though, you'll realize all this can be turned around. This means if you choose a vector of unit length $u$ and choose a number $\lambda$ any way you like, then there is exactly one pair of coordinates $(x,y)$, equivalently, exactly one vector $v=(x,y)$ such that the length of $v$ is $\lambda$ and the direction of $v$ is $u.$ The upshot is there are two ways of representing any point in the plane using either \textit{Cartesian} coordinates $(x,y)$ or \textit{polar} coordinates $(\lambda, u).$ Figure 9 gives some examples to illustrate the point.

\mpar
\subsection{Addition and Subtraction}
We've seen that you can multiply a vector by a real number. It turns out you can add vectors just like you can add numbers: all you do is just add their corresponding coordinates together. For example if $v=(a,b)$ and $w=(x,y)$ are vectors specified by their Cartesian coordinates then their sum is the vector with coordinates,
$$
v+w=(a+x, b+y).
$$
You can visualize this bit of algebra in the following way. (See Figure 10). Translate $v$ along $w$ from $(0,0),$ the base point of $w,$ to $(x,y),$ the endpoint of $w$ making sure that the translation of $v$ remains parallel to $v$ and doesn't change its length in the process - no twisting or stretching allowed! Let $v^\prime$ be the translation of $v$ along $w.$ The endpoint of $v^\prime$ is a point in the plane with coordinates  $(a+x, b+y).$ The vector you get by drawing a line from $(0,0)$ to $(a+x, b+y)$ is the vector $v+w.$ 

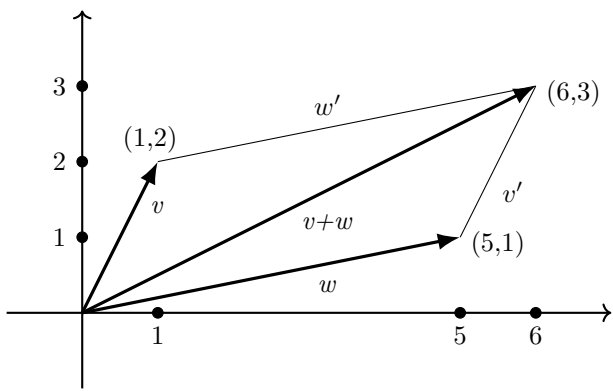
\begin{figure} [h]
\begin{tikzpicture}
\draw[thick,->] (0,0) -- (8,0);
\draw[thick,->] (1,-1) -- (1,4);

\draw[very thick,-Latex] (1,0) -- (2,2);
\draw[very thick,-Latex] (1,0) -- (6,1);

\draw[thin] (2,2) -- (7,3);
\draw[thin] (6,1) -- (7,3);

\draw[very thick, -Latex] (1,0) -- (7,3);

\draw[fill=black] (2,0) circle (2pt);
\draw[fill=black] (6,0) circle (2pt);
\draw[fill=black] (7,0) circle (2pt);

\node at (2,-0.3) {1};
\node at (6,-0.3) {5};
\node at (7,-0.3) {6};

\draw[fill=black] (1,1) circle (2pt);
\draw[fill=black] (1,2) circle (2pt);
\draw[fill=black] (1,3) circle (2pt);

\node at (0.7, 1) {1};
\node at (0.7, 2) {2};
\node at (0.7, 3) {3};

\node at (1.9,2.3) {(1,2)};
\node at (6.5,0.9) {(5,1)};
\node at (7.5,2.9) {(6,3)};

\node at (4.25, 0.35) {\textit{w}};
\node at (4.25, 1.2) {\textit{v+w}};
\node at (4.25, 2.75) {\textit{$w^\prime$}};

\node at (2,1.4) {\textit{v}};
\node at (6.7,1.6) {\textit{$v^\prime$}};
\end{tikzpicture}
\caption{The parallelogram determined by $v=(1,2)$ and $w=(5,1).$ The top and bottom of the rectangle are formed by $w^\prime$ and $w$ while the left and right sides of the rectangle are formed by $v$ and $v^\prime.$ The diagonal of the rectangle is the sum $v+w.$}
\end{figure}

\mpar
You can do this the other way around, if you like. Translate $w$ along $v$ from $(0,0),$ the base point of $v,$ to $(a,b),$ the endpoint of $v,$ making sure that the translation of $w,$ call it $w^\prime,$ remains parallel to $w$ and doesn't change its length in the process. The endpoint of $w^\prime$ is a point in the plane with coordinates  $(a+x, b+y),$ just like $v^\prime.$ The vector you get by drawing a line from $(0,0)$ to $(a+x, b+y)$ is the same vector as before. This just means that $v+w$ is the diagonal of the parallelogram formed by $v, w, v^\prime,$ and $w^\prime.$ (Again, see Figure 10).

\mpar
To subtract $w$ from $v$ just add $v$ and $-w,$ meaning,
$$
v-w= v+(-w) = (a,b)+(-x, -y) = (a-x, b-y).
$$
Makes sense, right?

\mpar
\subsection{Distance and Angle Between Vectors}
We can use subtraction of vectors to define a notion of distance between them. If $v$ and $w$ are two vectors then $v-w$ is a third vector and it seems pretty clear that $v=w$ if and only if $v-w=0.$ This leads to the idea that we can define the distance between $v$ and $w$ as the length of $v-w.$ Let's modify Figure 10 in order to visualize the difference. The key observation is that $v-w$ is the vector you must add to $w$ in order to get $v.$ This just means $v$ is the diagonal of the rectangle formed by $w$ and $v-w.$ (See Figure 11).

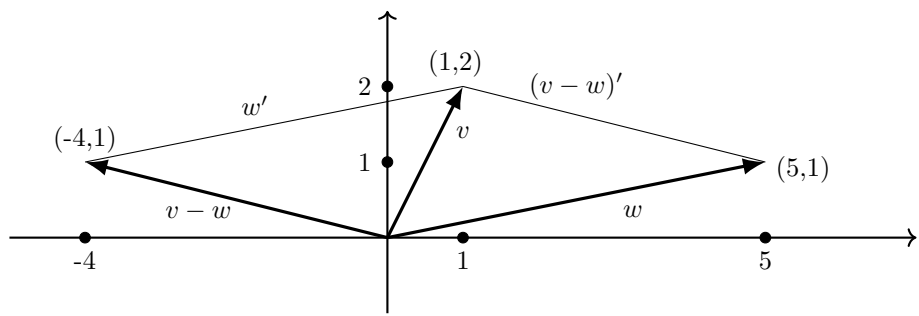
\begin{figure} [h]
\begin{tikzpicture}
\draw[thick,->] (-4,0) -- (8,0);
\draw[thick,->] (1,-1) -- (1,3);

\draw[very thick,-Latex] (1,0) -- (2,2);
\draw[very thick,-Latex] (1,0) -- (6,1);

\draw[thin] (-3, 1) -- (2,2);
\draw[thin] (6,1) -- (2,2);

\draw[very thick, -Latex] (1,0) -- (-3,1);

\draw[fill=black] (2,0) circle (2pt);
\draw[fill=black] (-3,0) circle (2pt);
\draw[fill=black] (6,0) circle (2pt);

\node at (2,-0.3) {1};
\node at (6,-0.3) {5};
\node at (-3,-0.3) {-4};

\draw[fill=black] (1,1) circle (2pt);
\draw[fill=black] (1,2) circle (2pt);

\node at (0.7, 1) {1};
\node at (0.7, 2) {2};

\node at (1.9,2.3) {(1,2)};
\node at (6.5,0.9) {(5,1)};
\node at (-3,1.3) {(-4,1)};

\node at (4.25, 0.35) {\textit{$w$}};
\node at (-1.5, 0.35) {\textit{$v-w$}};
\node at (-0.75, 1.75) {\textit{$w^\prime$}};
\node at (3.5, 2) {\textit{$(v-w)^\prime$}};
\node at (2,1.4) {\textit{$v$}};
\end{tikzpicture}
\caption{The parallelogram determined by $v=(1,2)$ and the difference $v-w=(1-5, 2-1) =(-4,1).$ The top and bottom of the rectangle are formed by $w^\prime$ and $w$ while the left and right sides of the rectangle are formed by $v-w$ and $(v-w)^\prime.$ The diagonal of the rectangle is the sum $v= w+(v-w).$}
\end{figure}

\mpar
It's easy to calculate distances based on what we already know about vectors. In fact,
$$
||v-w|| = ||(a-x, b-y)|| =\sqrt{(a-x)^2 +(b-y)^2}.
$$

\mpar
This formula just says that if you know the coordinates of $v$ and $w$ as numbers then you can calculate the distance between them as a number, just by plugging the coordinates into the formula above. Let's use the vectors in Figure 11 as an example. Since $v=(1,2)$ and $w=(5,1),$ the distance between them is,
$$
||v-w|| =\sqrt{(1-5)^2 +(2-1)^2} =\sqrt{4^2+1^2} =\sqrt{17} \cong 4.123.
$$
It may seem silly at first, but it's important to calculate the distance between $v$ and the zero vector, $w=(0,0).$ In this case, it's easy to see that $||v-w|| =\sqrt{a^2+b^2}$ which we also denote by $||v||.$ We call $||v||$ the \textit{length} of $v.$
\mpar
It turns out you can visualize distances between vectors in terms of the angle between them. It's easiest to start by looking at two vectors of unit length; meaning two vectors whose end point lies on the circle of radius 1, like those in Figure 9. In this case, distance clearly has something to do with the angle between the vectors: the smaller the angle is, the closer the coordinates must be to one another, and therefore the smaller is the distance between the two vectors. A diagram makes this clear. (See Figure 12).

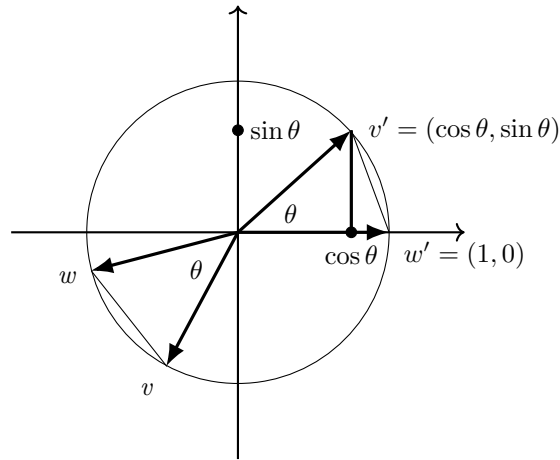
\begin{figure} [h]
\begin{tikzpicture}
\draw (0,0) circle (2);
\draw[thick, ->] (0,-3) -- (0,3);
\draw[thick, ->] (-3,0) -- (3,0);

\draw[very thick, -Latex] (0,0) -- (-1.95, -0.51);
\draw[very thick, -Latex] (0,0) -- (-0.95, -1.775);
\draw[thin] (-1.94, -0.51) -- (-0.94, -1.775);
\draw[very thick, -Latex] (0,0) -- (2, 0);
\draw[very thick, -Latex] (0,0) -- (1.5, 1.35);
\draw[very thick] (1.5,0) -- (1.5, 1.35);
\draw[thin] (2,0) -- (1.5, 1.35);
\node at (-0.55, -0.5) {$\theta$};
\node at (0.7, 0.25) {$\theta$};

\draw[fill=black] (1.5, 0) circle (2pt);
\draw[fill=black] (0, 1.35) circle (2pt);

\node at (0.5, 1.35) {$\sin\theta$};
\node at (1.5, -0.3) {$\cos\theta$};
\node at (3, 1.35 ) {$v^\prime =(\cos\theta, \sin\theta)$};
\node at (3, -0.3 ) {$w^\prime =(1,0)$};
\node at (-1.2, -2.075) {\textit{v}};
\node at (-2.25, -0.6) {\textit{w}};
\end{tikzpicture}
\caption{Two vectors $v$ and $w$ on the unit circle are pictured in the third quadrant. The length between them is length of the line segment joining their endpoints. The angle between them is $\theta.$ Clearly, the smaller the angle $\theta$, the smaller the distance $||v-w||.$ Rotate $w$ so that it lies on the horizontal axis and then rotate $v$ by the same amount. The vectors $v^\prime$ and $w^\prime$ are the rotations of $v$ and $w,$ respectively. Rotation doesn't change the angle between vectors so the angle between $v^\prime$ and $w^\prime$ is also $\theta.$ The coordinates of $w^\prime$ are $(1,0)$ because $w^\prime$ lies on the unit circle. The coordinates of $v^\prime$ are $(\cos\theta, \sin\theta)$ by the definition of $\cos$ and $\sin$. We'll use this fact to find a formula for $||v-w||$ in terms of $\theta.$}
\end{figure}

\mpar
In case you've forgotten your trigonometry (or never had the chance to learn it!), the \textit{sine} and \textit{cosine} of an angle are defined as certain ratios of side lengths in a right triangle. The cosine is the ratio of the length of the side adjacent to the angle to the length of the hypotenuse. And the sine is the ratio of the length of the side opposite from the angle to the length of the hypnotenuse. We used this definition in Figure 13 to figure out the coordinates of $v^\prime.$ In the right triangle formed by $v^\prime$ along with the horizontal and vertical axes, the hypotenuse has length one because $v^\prime$ lies on the  unit circle. Thus $\cos\theta$ is the length of the adjacent side, meaning the horizontal coordinate of $v^\prime,$ and $\sin\theta$ is the length of the opposite side, meaning the vertical coordinate of $v^\prime.$

\mpar
Now we can actually calculate the distance between $v$ and $w$ in terms of $\theta.$ This is the same as the distance between $v^\prime$ and $w^\prime$ because rotating two vectors by the same amount doesn't change the distance between them. Actually, it's easier to compute the square of the distance, because it get's rid of the square root sign, so let's go that route. We have,
\begin{align*}
||v-w||^2 & = ||v^\prime - w^\prime||^2\\
& =(\cos\theta -1)^2+(\sin\theta -0)^2\\
& = \cos^2\theta -2\cos\theta +1+\sin^2\theta\\
& = (\sin^2\theta +\cos^2\theta)  + 1 -2\cos\theta\\
& =2-2\cos\theta.
\end{align*}

In the last line we used the fact that $\sin^2\theta + \cos\theta =1.$ Why is this true? Well, it follows from Pythagoras theorem applied to the right triangle depicted in Figure 12 because the hypotenuse of that triangle has length 1. Thus $||v-w||=\sqrt{2(1-\cos\theta)}$. It's clear from Figure 12 that the closer $\theta$ is to zero, the closer $\cos\theta$ is to 1, so this gives a pretty good quantitative understanding of the distance between unit vectors in terms of the angle between them. 

\mpar
But what about the case that $v$ and $w$ are any two vectors, regardless of their length? Let's figure it out in two steps. The first step is to write out the formula for the squared length of $v-w$ which is,
\begin{align*}
||v-w||^2 & = (a-x)^2+(b-y)^2\\
& = (a^2-2ax +x^2) + (b^2-2by +y^2)\\
& = (a^2+b^2) +(x^2+y^2) -2(ax+by).
\end{align*}

If you think about it for a minute, you'll recognize that the first two terms above are the squared lengths of $v$ and $w$; that is,
$$
(a^2+b^2) = ||v||^2\quad\text{and}\quad (x^2+y^2)=||w||^2.
$$
The third term above involves something new that we're going to call the \textit{inner product} or sometimes the \textit{dot product} of $v$ and $w,$
$$
v\cdot w = (a,b)\negthinspace\cdot\negthinspace (x,y) = ax+by.
$$
Notice that the formula for the dot product just tells you to multiply corresponding coordinates together and add them up. You can also see that if $\lambda$ and $\mu$ are numbers then,
$$
(\lambda\, v)\negthinspace\cdot\negthinspace(\mu\, w) = \lambda\,\mu\, (v\negthinspace\cdot\negthinspace w).
$$

Given the definition of dot product, we can write the squared distance between two vectors as,
$$
||v-w||^2 = ||v||^2+||w||^2 -2\,\,v\negthinspace\cdot\negthinspace w.
$$

\mpar
Now, the second step is to see that the dot product of two vectors is related to the angle between them. If $v$ and $w$ happen to have length 1 we already know that,
$$
||v-w||^2 = 2-2\cos\theta,
$$
and we just figured out that,
$$
||v-w||^2 = ||v||^2+||w||^2 -2\,\,v\negthinspace\cdot\negthinspace w = 2-2v\negthinspace\cdot\negthinspace w,
$$
so if you equate the two formulas for $||v-w||$ you find that,
$$
v\negthinspace\cdot\negthinspace w=\cos\theta.
$$

This formula just says the dot product between two vectors of length 1 equals their dot product. OK, so what if $v$ and $w$ don't have unit length; what then? There's a simple trick for that - just divide $v$ and $w$ by their length,
$$
v^\prime =\frac{v}{||v||}\quad\text{and}\quad w^\prime=\frac{w}{||w||}
$$
and observe that,
$$
v=||v||v^\prime\quad\text{and}\quad w=||w|| w^\prime.
$$
But $v^\prime$ and $w^\prime$ have unit length so,
$$
v\negthinspace\cdot\negthinspace w = ||v||\, ||w||\, v^\prime\negthinspace\cdot\negthinspace w^\prime = ||v||\, ||w||\, \cos\theta.
$$

\subsection{Higher Dimensions}
We have been working with two dimensional vectors because $2D$ geometry is pretty easy to visualize. But vectors come in many dimensions like the three dimensional space we live in, where $v=(a, b, c),$ or the four dimensions of spacetime where $v=(a,b, c, t).$  Here, the first three coordinates refer to the three familiar directions of space: \textit{North/South, East/West, Up/Down} and the fourth coordinate refers to the time dimension: \textit{Past/Future}.

\mpar
We say $v$ is an \textit{n-dimensional} vector if it has $n$ coordinates and we write,
$$
v=(v_1, v_2, \cdots, v_n).
$$
where $v_1, v_2, \cdots, v_n$ are arbitrary numbers. Once the dimension is $5$ or more, it's hard at first to see how $n$-dimensional vectors are relevant, but it turns out that they are really all around you if you know where to look.

\mpar
For example, every image stored on your phone is represented by a very high dimensional vector. A black-and-white digital photo is basically a rectangle made up of little squares called \textit{pixels} and each pixel contains a number between $0$ and $255$ representing the \textit{greyscale}. Here $0$ denotes black, $255$ denotes white, and the numbers in between represent varying shades of grey. 

\mpar
A typical image taken with a $3$-megapixel camera has something like $1200\times 2100$ pixels which we can think of as a vector where the coordinates are the greyscales. If you do the math, each picture on your phone is a vector in $2,940,000$-dimensional space. This is the typical dimension of vectors processed by commercial neural networks like the one under the hood of ChatGPT.

\mpar
We want to extend a couple of formulas from two dimensional vectors to all higher dimensions: the formula for length of a vector and the formula for the angle between two vectors. We'll need them when we discuss multivariable calculus. A picture of the three dimensional case illustrates the general argument. (See Figure 13).

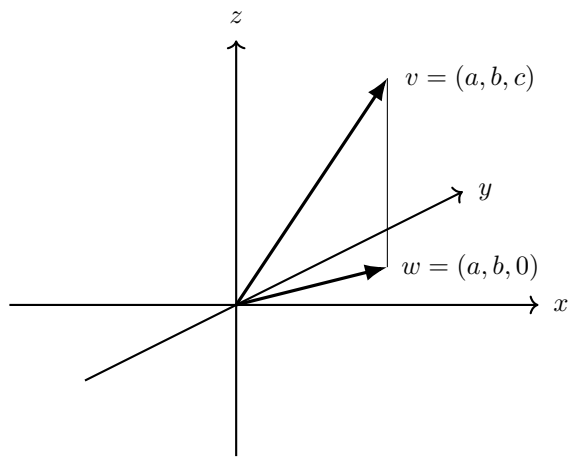
\begin{figure} [h]
\begin{tikzpicture}
\draw[thick, ->] (0,-2) -- (0,3.5);
\draw[thick, ->] (-3,0) -- (4,0);
\draw[thick, ->] (-2, -1) -- (3,1.5);
\draw[very thick, -Latex] (0,0) -- (2,3);
\draw[thin] (2,3) -- (2, 0.5);
\draw[very thick, -Latex] (0,0) -- (2, 0.5);

\node at (0, 3.8) {$z$};
\node at (4.3, 0) {$x$};
\node at (3.3, 1.5) {$y$};
\node at (3.1,3) {$v=(a,b,c)$};
\node at (3.1,0.5) {$w=(a,b,0)$};
\end{tikzpicture}
\caption{A right triangle formed by vectors $v$ and $w,$ where $v$ is the hypotenuse, shown in the $x-y-z$ coordinate system. This means the length of $v$ can be calculated from Pythagoras theorem, once you know the lengths of the other two sides.} 
\end{figure}

In Figure 13, $v$ is a typical $3D$ vector and $w$ is the vector you get by setting the $z$ coordinate of $v$ equal to $0.$ Observe that $v$, $w$ and the line segment connecting the end points of $v$ and $w$ form a right triangle. Because its $z$ coordinate is zero, $w$ lies in the $x-y$ plane so we can think of it as a $2D$ vector with coordinates $(a,b)$. Thus, the length of the hypotenuse of this triangle is $||v||,$ the length of $v,$ and the side lengths of the triangle are $||(a,b)|| = \sqrt{a^2+b^2}$ and $|c|,$ the absolute value of $c$. (Notice that $c$ could be a negative number although its a positive number as drawn in Figure 13). But by Pythagoras theorem,
\begin{align*}
||v||^2 & = ||(a,b)||^2 + |c|^2\\
& = (\sqrt{a^2+b^2}\,)^2 + |c|^2\\
& = a^2+b^2+c^2,
\end{align*}
and therefore,
$$
||v|| =\sqrt{||v||^2} = \sqrt {a^2+b^2+c^2}.
$$
\mpar
This formula says that the length of a $3D$ vector is just the square root of the sum of the squares of its coordinates. And when you say it this way, you realize you didn't need to know exactly how many coordinates there are, so the statement ought to hold in any number of dimensions. Like a lot of things in mathematics, once you see it, it's easy to prove. Here's how.

\mpar
Let $v=(v_1, v_2,\cdots , v_{n-1}, v_n)$ be an $n$ dimensional vector and let,
$$
w=(v_1, v_2,\cdots, v_{n-1}, 0)
$$ 
be the vector obtained by setting the $n^{th}$ coordinate of $v$ equal to zero. We're going to think of $w$ as an $n-1$ dimensional vector, just like we did in the $3D$ case depicted in Figure 13, and we're going to suppose we already proved that,
$$
||w|| = \sqrt{v_1^2+v_2^2 +\cdots +v_{n-1}^2}.
$$
Then $v$ and $w$ are two sides of a right triangle whose third side has length $|v_n|.$ Therefore, by Pythagoras theorem,
\begin{align*}
||v||^2 & = ||w||^2 + |v_n|^2\\
& = (\sqrt{v_1^2+v_2^2 +\cdots +v_{n-1}^2}\,)^2 + |v_n|^2\\
& =v_1^2+v_2^2 +\cdots + v_{n-1}^2+v_n^2,
\end{align*}
and it follows that,
$$
||v|| = \sqrt{||v||^2} = \sqrt{v_1^2+v_2^2 +\cdots + v_{n-1}^2+v_n^2}.
$$

So, what's happening here? What we've shown is a conditional statement: \textit{if} the desired formula is true in a given dimension \textit{then} it is true in the next higher dimension. But it is true in $2D,$ so it's also true in $3D$. But since its true in $3D$ it must also be true in $4D,$  and so on, ad infinitum. This kind of two-step argument is called \textit{mathematical induction}: you can prove infinitely many statements at a stroke by proving a conditional (the truth of one statement implies the truth of the next statement) and then proving an absolute (the first statement is true). Mathematical induction is like setting up a logical chain reaction then lighting the fuse!

\mpar
The final thing we need to know is the formula connecting the dot product of two vectors with the angle between them is true in any dimension. This amounts to revisiting the calculations we did in the $2D$ case. Suppose $v$ and $w$ are $n$ dimensional vectors. Then,
\begin{align*}
& ||v-w||^2  = (v_1-w_1)^2 +(v_2-w_2)^2+\cdots + (v_n-w_n)^2\\
& = (v_1^2-2v_1w_1 +w_1^2) +(v_2^2-2v_2w_2 +w_2^2) + \cdots + (v_n^2-2v_nw_1 +w_n^2)\\
& = (v_1^2+v_2^2+\cdots + v_n^2)+(w_1^2+w_2^2+\cdots + w_n^2) -2(v_1w_1+v_2w_2+\cdots +v_nw_n)\\
& = ||v||^2+||w||^2-2\,\,v\negthinspace\cdot\negthinspace w.
\end{align*}

Here we've used the notation $v\negthinspace\cdot\negthinspace w$ for the dot product,
$$
v\negthinspace\cdot\negthinspace w = v_1w_1+v_2w_2+\cdots +v_nw_n.
$$

If $v$ and $w$ have unit length then this formula becomes,
$$
||v-w||^2=2-2\,\,v\negthinspace\cdot\negthinspace w
$$
Now, here's the key insight. Two distinct vectors lie in a plane just like two distinct point lie on a line. So the geometry of $v$ and $w$, thought of a $n$ dimensional vectors, is identical to the geometry of the two dimensional vectors $v^\prime$ and $w^\prime$ in Figure 12 and it follows that,
$$
||v-w||^2=2-2\cos\theta
$$
where $\theta$ is the angle between them. On comparing these two expressions for $||v-w||^2$ we find,
$$
v\negthinspace\cdot\negthinspace w = \cos\theta.
$$
This gives the the formula,
$$
||v-w||^2 = ||v||^2 + ||w||^2 - 2\,||v||\,||w|| \cos\theta
$$
using the same trick as before - multiplying and dividing $v$ and $w$ by their lengths.

\mpar
\subsection{Matrices and Linear Transformations} 
When we look at the feed forward and backpropagation steps of neural networks we will need to know something about matrices acting on vectors and products of matrices creating new matrices. This subsection provides a quick and dirty introduction to these ideas.

\mpar
Other than adding vectors together, the simplest way to change one vector $x$ into another vector $y$ is to make the coordinates of $y$ be functions of the coordinates of $x$. To keep things a simple as possible let's make them linear functions like this,
$$
x=(x_1, x_2)\to y=(y_1, y_2) = (2x_1-3x_2, -4 x_1+ 5 x_2).
$$
This is equivalent to saying,
\begin{align*}
y_1 & = \phantom{-}2x_1 -3x_2\\
y_2 & = -4 x_1+5 x_2
\end{align*}
which we represent graphically as,
$$
\begin{pmatrix}
y_1\\
y_2
\end{pmatrix}
=
\begin{bmatrix}
\phantom{-}2 & -3\\
-4 & \phantom{-}5
\end{bmatrix}
\begin{pmatrix}
x_1\\
x_2
\end{pmatrix}
$$
\mpar
The square of numbers above is called a \textit{matrix} or more specifically a $2\times 2$ dimensional matrix which stands for 2 rows and 2 columns. The graphical representation invites you to make $2D$ vectors from the rows of the matrix, namely $(2, -3)$ and $(-4,5),$ and take their inner product with the vector $(x_1, x_2)$ to get the vector $(y_1, y_2).$ Specifically,
\begin{align*}
y_1 & = 
\begin{pmatrix}
\phantom{-}2\\
-3
\end{pmatrix}
\cdot
\begin{pmatrix}
x_1\\
x_2
\end{pmatrix}
=\phantom{-}2x_1-3x_2\\
y_2 & =\begin{pmatrix}
-4\\
\phantom{-}5
\end{pmatrix}
\cdot
\begin{pmatrix}
x_1\\
x_2
\end{pmatrix}
=-4x_1+5x_2
\end{align*}

Of course there's noting special about the numbers $2,-3, -4, 5$ of the matrix, they could just as well be any numbers whatsoever. Let $A$ be the $2\times 2$ matrix,
$$
A = \begin{bmatrix}
a_{11} & a_{12}\\
a_{21} & a_{22}
\end{bmatrix}
$$
where $a_{ij}$ are arbirtary. Then the equation $y=Ax$ stands for
\begin{align*}
y_1 & = 
\begin{pmatrix}
a_{11}\\
a_{12}
\end{pmatrix}
\cdot
\begin{pmatrix}
x_1\\
x_2
\end{pmatrix}
=a_{11}x_1+a_{12}x_2\\
y_2 & =\begin{pmatrix}
a_{21}\\
a_{22}
\end{pmatrix}
\cdot
\begin{pmatrix}
x_1\\
x_2
\end{pmatrix}
=a_{21} x_1+a_{22}x_2
\end{align*}

\mpar
Let's step away from the actual formulas and look at the bigger picture. $Ax$ is the vector you get by taking the inner product of the rows of the matrix $A$ with vector $x$. Keeping this in mind it's easy to see that if $x$ and $y$ are vectors then,
$$
A(x+y)= Ax+Ay.
$$
Also if $\lambda$ and $\mu$ are numbers then for the same reason,
$$
A(\lambda x+\mu y) = \lambda Ax +\mu Ay.
$$
Thus every $2\times 2$ dimensional matrix $A$ defines a function $A\colon \mathbb{R}^2\to\mathbb{R}^2$ having the property $A(\lambda x+\mu y) = \lambda Ax +\mu Ay.$ Functions with this property are called \textit{linear transformations}.

\mpar
It's extremely important to know, and not so hard to prove, that every linear transformation is defined by a matrix $A$. All you need to observe is that $A$ transforms the special vectors $e_1=(1,0)$ and $e_2=(0,1)$ into some vectors say,
$$
A e_1=(a,b)\quad\text{and}\quad A e_2=(c,d)
$$
and then the matrix,
$$
A = \begin{bmatrix}
a & b\\
c & d
\end{bmatrix}
$$
does the trick! \textit{Every linear transformation from $\mathbb{R}^2$ to $\mathbb{R}^2 $is represented by a $2\times 2$ matrix.}

\mpar
There's nothing special about $2\times 2$ matrices and there's a lot to be gained by considering $n\times m$ dimensional matrices. These are rectangles of numbers having $n$ rows and $m$ columns such as,
$$
A=
\begin{bmatrix}
a & b & c
\end{bmatrix}, \quad
B=
\begin{bmatrix}
d & e \\
f & g\\
h & i
\end{bmatrix}, \quad
C=
\begin{bmatrix}
j& k & l & m\\
n & o & p & q\\
r & s & t & u
\end{bmatrix}, \quad
D=
\begin{bmatrix}
v & w & x\\
y& z & \alpha
\end{bmatrix}.
$$
Observe that $A$ is $1\times 3$ dimensional matrix, $B$ is a $3\times 2,$ $C$ is a $3\times 4,$ and $D$ is a $2\times 3$ dimensional matrix. It turns out that if you get the dimensions right you can multiply two matrices together. Here's how.

\mpar
Suppose $E$ and $F$ are matrices and we want to know if the product $EF$ makes sense as a matrix. Taking a hint from the formula for $Ax$ above, if the number of columns of $E$ equals the number of rows of $F$ then we can think of a row of $E$ and a column of $F$ as vectors of the same dimension and calculate their inner product. This produces a number which depends on exactly which row and column you chose. Of course, you can do this for all possible choices of rows from $E$ and columns from $F$. 

\mpar
This may sounds confusing when you write it out in words, but it's actually pretty simple once you get the hang of it. Let's see how this works in the examples above.

\mpar
First of all, $AB$ makes sense because $A$ has 3 columns and $B$ has 3 rows. There are just two inner products,
$$
(a, b, c)\cdot (e,f,h)\quad\text{and}\quad (a,b,c)\cdot (e,g,i)
$$ 
and you can arrange them in a $1\times 2$ dimensional matrix,
$$
 AB=
\begin{bmatrix}
 ad+bf+ch & ae+ bg+ci
\end{bmatrix}.
$$
 
\mpar
There are two other dimensionally correct products: $AC, BD,$ and  $DB,$ meaning the number of columns of the first factor equals the number of rows of the second factor. For example $AC$ is a $1\times 3$ times a $3\times 4$ matrix, $BD$ is a $3\times 2$ times a $2\times 3$ matrix, and $DB$ is a $2\times 3$ times a $3\times 2$ matrix.

\mpar
Let's work out the product $DB$ by hand just to get some practice. There are two rows in $D$ and two columns in $B$ so there are four inner products,
$$
(v,w,x)\cdot (d,e,f)\quad (v,w,x)\cdot (e,g,i)\quad (y,z, \alpha)\cdot (d,f,h)\quad (y,z,\alpha)\cdot (e,g,i).
$$

You can arrange these numbers in a $2\times 2$ dimensional matrix this way,
$$
DB =
\begin{bmatrix}
vd+we+xf & ve+ wg+xi\\
yd+zf +\alpha h & ye+zg+\alpha i
\end{bmatrix}.
$$
In self-explanatory notation it has the pattern,
$$
DB =
\begin{bmatrix}
row_1(D)\cdot col_1(B) & row_1(D)\cdot col_2(B)\\
row_2(D)\cdot col_1(B) & row_2(D)\cdot col_2(B)
\end{bmatrix}.
$$
If you work out $AC$ you'll discover its a $1\times 4$ dimensional matrix and if you do the same for $BD$ you'll find it is a $3\times 3$ dimensional matrix. In each case the product matrix has entries of the form
$$
row_i(1^{st} factor)\cdot col_j(2^{nd} factor)
$$

\mpar
That's the general rule. If $S$ is an $n\times m$ matrix and $B$ is a $m\times q$ matrix then $A$ rectangle of numbers having $n$ rows and $m$ columns and $B$ is a rectangle of numbers having $m$ rows an $q$ columns which we write as
$$
A=
\begin{bmatrix}
a_{11} & a_{12} & \dots & a_{1m}\\
a_{21} & a_{22} & \dots & a_{2m}\\
\vdots & \vdots &\ddots &\vdots\\
a_{n1} & a_{n1} & \dots & a_{nm}
\end{bmatrix}
\quad
\quad
B=
\begin{bmatrix}
b_{11} & b_{12} & \dots & b_{1q}\\
b_{21} & b_{22} & \dots & b_{2q}\\
\vdots & \vdots &\ddots &\vdots\\
a_{m1} & a_{m1} & \dots & a_{mq}
\end{bmatrix}.
$$
The can be multiplied together in the correct order: $C=AB$ is a $n\times q$ matrix,
$$
C=
\begin{bmatrix}
c_{11} & c_{12} & \dots & c_{1q}\\
c_{21} & c_{22} & \dots & c_{2q}\\
\vdots & \vdots &\ddots &\vdots\\
c_{n1} & c_{n1} & \dots & c_{nq}
\end{bmatrix}
$$
where,
$$
c_{ij}=\sum_{k=1}^ma_{ik}b_{kj}.
$$
This formula just says that $c_{ij}$ is the inner product of the $i^{th}$ row of $A$ and the $j^{th}$ column of $B$.

\mpar
To sum it all up here are the main points:
\begin{itemize}
\item{A linear transformation $T\colon\mathbb{R}^m\to\mathbb{R}^n$ is a function such that 
$$
T(\lambda x+\mu y)=\lambda T(x)+\mu T(y)
$$
for all vectors $x,y\in\mathbb{R}^m$ and $\lambda,\mu\in\mathbb{R}.$}
\item{If $A$ is $n\times m$ matrix then $A$ defines a linear transformation from $\mathbb{R}^m$ to $\mathbb{R}^n$ by the rule,
$$
Ax=\left(\sum_{k=1}^m a_{1,k} x_k,\, \sum_{k=1}^m a_{2,k} x_k,\,\dots,\, \sum_{k=1}^m a_{n,k} x_k\right),
$$
any every linear transformation arises this way.}
\item{If $A$ is a $n\times m$ dimensional matrix and $B$ is a $m\times q$ dimensional matrix then the product $AB$ is an $n\times q$ dimensional matrix.}
\end{itemize}

\mpar
\mpar
\section{Calculus}
In this appendix we will discuss some of the most important ideas of calculus - functions and their graphs, linear functions as local approximations of more general functions, and how to calculate these local approximations by derivatives. We will focus first on functions of one variable and then consider functions of two or more variables. 

\mpar
The idea of the gradient of a function is the key concept. We'll see that the gradient is a vector whose \textit{direction} is the direction the function is changing most rapidly and whose \textit{length} is the rate at which it's changing. This fact is the basis of the \textit{method of steepest descent}, the workhorse algorithm for approximating the minimum of a function of many variables.

\mpar
\subsection{Functions and their Graphs}
One way to think of a function is as a rule, formula, or recipe for calculating the answer to a question, given some input. For example, if the question is \textit{``What is the length of a vector?''} then the answer is a non-negative real number, the input is the coordinates of the vector, and the formula is,
$$
\lambda = f(x,y)=\sqrt{x^2+y^2}.
$$
We use the notation $f(x,y)$ to emphasize the input variables and we write $\lambda = f(x,y)$ to emphasize that the length $\lambda$ can be calculated from the coordinates by a very specific formula.

\mpar
Here's some standard notation that is commonly used in calculus. We write $\mathbb{R}$ to denote the set of all real numbers. This includes whole numbers like $0, -1, 13,$ as well as fractions like $\tfrac{1}{2}, -\tfrac{4}{3}, \tfrac{219}{4396},$ and also irrational numbers like $\sqrt{2}, e,$ and $\pi.$ You can think of $\mathbb{R}$ as a straight line, infinite in both directions, and visualize numbers as points on the line. Every point on the line is a real number and every real number is a point on the line.

\mpar
Two dimensional vectors are pairs of points and we use the notation $\mathbb{R}^2$ to denote the set of all such pairs. For example $(-\tfrac{1}{2}, \pi)$ is a pair in $\mathbb{R}^2$ and we denote it by writing $(-\tfrac{1}{2}, \pi)\in\mathbb{R}^2.$ The notation suggests $\mathbb{R}^2=\mathbb{R}\times\mathbb{R},$ which is a good way to think of it because encourages you to think of $\mathbb{R}^2$ as a plane - a product to two perpendicular lines, infinite in all directions - with $(0,0)$ being the place the two lines intersect, the point we call the origin. There is $\mathbb{R}^3=\mathbb{R}^2\times\mathbb{R}$ and  $\mathbb{R}^4=\mathbb{R}^3\times\mathbb{R},$ etc., for any number of dimensions. $\mathbb{R}^n$ is called the Euclidean space of dimension $n$.

\mpar
Here a final bit of notation. We write $f\colon\mathbb{R}^2\to\mathbb{R}$ as shorthand for the fact that the function $f(x,y)=\sqrt{x^2+y^2}$ takes a pair of real numbers as input and gives a single real number as output. There are many other functions, for example $f(x,y)=y^3+x^2+1$ or $g(x,y)= e^x-e^y,$ and we use the letters $f$ and $g$ to stand for any function of two variables in the same way that we use the letters $x$ and $y$ to stand for any real number.

\mpar
The \textit{graph} of a function is the set of all pairs of the form $(\text{input}, \text{output}).$ In the example of the function calculating the length of a $2D$ vector $v$, the input/output pairs have the form, 
$$
(\text{input}, \text{output})=(v, ||v||) =(x,y,\sqrt{x^2+y^2}),
$$
so the graph of this function is a collection of $3D$ vectors. You can visualize the graph of this function as a surface in $\mathbb{R}^3$ lying above the $x-y$ plane: each input $v=(x,y)$ corresponds to the $3D$ vector $(x,y,0)$ and the the (input,output) vector  $(x,y,\sqrt{x^2+y^2})$ lies directly above it at the vertical height of $\sqrt{x^2+y^2}.$ As the input values of $x$ and $y$ vary, the output vector $(x,y,\sqrt{x^2+y^2})$ also varies and traces out a two dimensional surface as it does so.

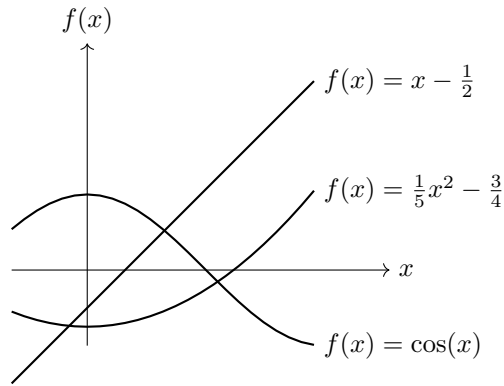
\begin{figure}[h]
\begin{tikzpicture} [domain =-1:3]
\draw[->](-1,0) -- (4,0) node[right] {$x$};
\draw[->](0,-1) -- (0,3) node[above] {$f(x)$};

\draw[thick] plot(\x,\x-0.5) node[right] {$f(x)=x-\tfrac{1}{2}$};
\draw[thick] plot(\x,{cos(\x r)}) node[right] {$f(x)= \cos(x)$};
\draw[thick] plot(\x,{0.2*(\x)^2-0.75}) node[right] {$f(x)= \frac{1}{5}x^2 -\frac{3}{4}$};
\end{tikzpicture}
\caption{The graphs of $f(x)=x+\tfrac{1}{2},\,\, f(x)=\tfrac{1}{5}x^2-\tfrac{3}{4},$ and $f(x)=\cos(x).$ Each graph is formed by the $2D$ vectors $(x, f(x))$ which lie directly over the the vectors $(x,0).$ As the values of $x$ vary, the vectors $(x,f(x))$ vary and so race out a curve in $\mathbb{R}^2.$}
\end{figure}

\mpar
It's helpful to looks at some examples where both input and output are one dimensional vectors, meaning they are just real numbers, so the (input,output) pairs on the graph trace out a one dimensional curve in $\mathbb{R}^2.$ See Figure 14.

\mpar
\subsection{Linear Functions as Local Approximations}
Linear functions are particularly important and understanding them is the key to understanding calculus. A function $f\colon\mathbb{R}\to\mathbb{R}$ is linear if its given by a formula of the form,
$$
f(x)=ax + b,
$$
and a function $g\colon\mathbb{R}^2\to\mathbb{R}$ is linear if it has the form,
$$
g(x,y)=ax +by +c.
$$
In these formulas, think of $a, b$ and $c$ as fixed real numbers like $a=7, b=\sqrt{3}$ and $c=-\tfrac{10}{7}$and think of $x$ and $y$ as variable real numbers in $\mathbb{R}.$ 

\mpar
Strictly speaking, $f$ and $g$ are \textit{affine} functions - they are linear functions if $b$ and $c$ equal zero. In either case, the graph of $f$ is a line in $\mathbb{R}^2$ and the graph of $g$ is a plane in $\mathbb{R}^3.$ If you think about it for a moment, the origin $(0,0)$ of $\mathbb{R}^2$ is in the graph of $f$ if and only if $b=0$ and the origin $(0,0,0)$ of $\mathbb{R}^3$  is in the graph of $g$ if and only if $c=0.$ So, the graphs of linear functions always pass through the origin but the graphs of affine functions do not. 
\mpar
Here is essential idea of calculus: \textit{it's useful to approximate the graph of a function $f$ at a point $(x,f(x))$ by the graph of an affine function that contains the point $(x,f(x)).$}

\mpar
This statement may seem a bit cryptic at first but a few pictures using a concrete function will help make sense of it. Lets take the function $f(x) =2-\tfrac{1}{2}x^2$ and try to approximate its graph at the point $(1,f(1))=(1,\tfrac{3}{2})$ and again at the point $(0,f(0))=(0,2).$

\begin{figure}[h]
\begin{tikzpicture} [domain =-1:3]
\draw[->](-1,0) -- (4,0) node[right] {$x$};
\draw[->](0,-1) -- (0,3) node[above] {$f(x)$};
\draw[thick] plot(\x,{2-0.5*(\x)^2}) node[right] {$f(x)= 2-\tfrac{1}{2}x^2$};
\draw[fill=black] (1,1.5) circle (2pt);
\draw[fill=black] (1,0) circle (2pt) node[below] {$1$};
\draw[fill=black] (0,1.5) circle (2pt) node[left] {$\tfrac{3}{2}$};
\draw[thin] (-1,-2.5) -- (2, 3.5);
\draw[thin] (-1,0.5) -- (3, 2.5);
\draw[thin] (-1,3.5) -- (3, -0.5);

\node at (2.5,3.5) {$s=2$};
\node at (3.5,2.5) {$s=\tfrac{1}{2}$};
\node at (3.7,-0.5) {$s=-1$};
\node at (1.75, 1.4) {$(1,\tfrac{3}{2})$};
\end{tikzpicture}
\caption{The graphs of three affine functions passing through the point $(1,\tfrac{3}{2})$ are labelled by their slope $s$. As the slope varies, the lines spin around the fixed point $(1,\tfrac{3}{2}).$ }
\end{figure}

\mpar
Of the three affine functions shown in Figure 15, the line with slope $s=-1$ seems intuitively like the best approximation. If the line were a little steeper it would cut through the graph of $f$  to the right of $(1,\tfrac{3}{2})$ and if it were a little flatter, it would cut through the graph of $f$ to the left of $(1,\tfrac{3}{2}).$  But when the slope is $-1$ the line seems balanced just right to exactly mirror the shape of the graph better than any other line. It is called the \textit{tangent line to $f$ at $x=1.$} The slope of the tangent line is called the \textit{derivative} of $f$ at $x=1$ and we denote it by writing $f^\prime(1)=-1.$

\mpar
You can find the equation of the tangent line using the facts that (1) the point $(1,\tfrac{3}{2})$ lies on the line and (2) the line has slope $-1.$ In fact, if $(x,y)$ is any other point on the tangent line then the \textit{rise} = $(y-\tfrac{3}{2})$ divided by the \textit{run} = $(x-1)$ equals the \textit{slope} = $rise\,/\, run =-1,$ and it follows that,
$$
(y-\tfrac{3}{2})=- (x-1),
$$ 
or, equivalently,
$$
y= -x +\tfrac{5}{2}.
$$
This is the equation of tangent line -  the line that best approximates $f$ at $x=1.$

\mpar
Just to get a bit more experience, let's find the affine function whose graph is the best approximation to the graph of $f$ at the point $(0,f(0))=(0,2).$ In other words, let find the slope of the tangent line to $f$ at $x=0.$

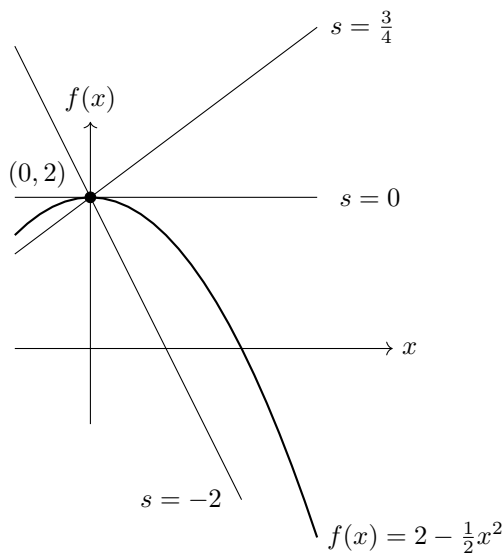
\begin{figure}[h]
\begin{tikzpicture} [domain =-1:3]
\draw[->](-1,0) -- (4,0) node[right] {$x$};
\draw[->](0,-1) -- (0,3) node[above] {$f(x)$};
\draw[thick] plot(\x,{2-0.5*(\x)^2}) node[right] {$f(x)= 2-\tfrac{1}{2}x^2$};
\draw[fill=black] (0,2) circle (2pt);
\node at (-0.7, 2.3) {$(0,2)$};
\draw[thin] (-1,2) -- (3,2);
\draw[thin] (-1,1.25) -- (3, 4.25);
\draw[thin] (-1,4) -- (2, -2);

\node at (3.6,4.25) {$s=\tfrac{3}{4}$};
\node at (3.7,2) {$s=0$};
\node at (1.2,-2) {$s=-2$};

\end{tikzpicture}
\caption{The graphs of three affine functions passing through the point $(0,2)$ are labelled by their slope $s$. As the slope varies, the lines spin around the fixed point $(0,2)$.}
\end{figure}

In Figure 16, the line with slope $0$ seems as if its the best approximation because lies flat at the top of the ``hill'' formed by the graph of $f.$ A line with negative slope would point downward to the right and a line with positive slope would point upward to the right. Only the line with slope $0$ exactly balances at the top of the ``hill.'' This just reflects the fact that as the input $x$ increases toward $0$ the output $f(x)$ increases toward $2$ and as the input $x$ increases away from $0$ the output decreases away from $2.$ 

\mpar
The line passing through the point $(0,2)$ with slope equal to zero is given by the equation $y=2.$ We can write it in the familiar form $y=ax+b$ by choosing $a=0$ and $b=2.$ This is the tangent line to $f$ at $x=0$ and it's slope is the derivative $f^\prime(0)=0.$ 

\mpar
These two examples show that finding the best linear approximation to $f$ at at point $x$ makes intuitive sense and in principle you could do it for any value of $x$, say $x=p,$not just $x=0$ or $x=1.$ A line passing through $(p, f(p))$ with slope $a$ has the form,
$$
(y-f(p))=a(x-p)
$$
which can be rewritten as,
$$
y=a(x-p)+f(p).
$$
So, the problem of finding the ``best'' linear approximation to $f$ at $x=p$ is the same as the problem of the ``correct'' choice of the slope $a$. If you can find the correct choice, we call it the derivative of $f$ and write it as $a=f^\prime(p).$ Then the equation for the tangent line becomes,

$$
y=f^\prime (p)(x-p)+f(p).
$$

\mpar
It's a pretty formula, right? It basically says if you give me a function $f$ then I'll give you it's tangent line at any point $p$ you like. All I need to do is find a practical method or an algorithm for computing the derivative $f^\prime (p).$

\mpar
\subsection{The Derivative}The hard part in finding such an algorithm is it depends of what the word ``best'' means when we ask for the best approximation to the graph of $f$ by the graph of an affine function. In the two examples above it was intuitively clear what best meant in geometric terms - something like ``balanced','' or ``mirrors'' the shape of the graph, or has the ``correct'' slope. But this just exchanges one unknown term like ``best'' for other unknown terms like ``balanced,'' ``mirrors,'' or ``correct''. It feels like we're stuck in a vicious circle.

\mpar
It took mathematicians of the late seventeenth and early eighteen centuries many years to finally break out of the vicious circle. To see how they did it let's go back the our example and figure out an algorithm for calculating $f^\prime(0),$ the derivative of $f$ at $x=0.$  If we can do it for this particular function at this particular point, we can do it for any function at any point. It takes two things: a key insight and a completely new idea. 

\mpar
The key insight is that two distinct points determine a line and we can calculate the slope of the line in terms of the coordinates of the two points by the rule,
$$
\textit{slope} = \frac{\textit{rise}}{\textit{run}}.
$$. 

Rise and run refer to a staircase: rise is the height of each step and run is its depth. 

\mpar
If $(a,b)$ and $(c,d)$ are the points then the rise is the difference in $y$ coordinates, $d-b,$ and the run is the difference in $x$ coordinates, $c-a$  so the slope is,
$$
s=\frac{d-b}{c-a}
$$
It doesn't matter which order you subtract them as long as its the same order for both rise and run. 

\mpar
Let's apply this insight to the function $f(x)=2-\tfrac{1}{2}x^2$ when $(0,f(0))=(0,2)$ is one of the points and we chose the other point to be $(p, f(p))$ where $p=1, \, p=\tfrac{1}{2},$ and $p=-\tfrac{1}{10}.$ The lines formed by $(0,f(0)$ and $(p,f(p))$ are called \textit{secant} lines. See Figure 17.

\begin{figure}[h]
\begin{tikzpicture} [domain =-1:3]
\draw[->](-1,0) -- (4,0) node[right] {$x$};
\draw[->](0,-1) -- (0,3) node[above] {$f(x)$};
\draw[thick] plot(\x,{2-0.5*(\x)^2}) node[right] {$f(x)= 2-\tfrac{1}{2}x^2$};

\draw[thin] (-1,2.5) -- (3,0.5);
\draw[thin] (-1,2.25) -- (3, 1.25);
\draw[thin] (-1,1.95) -- (3, 2.15);

\node at (3.7,0.5) {$s=-\tfrac{1}{2}$};
\node at (3.7, 1.25) {$s=-\tfrac{1}{4}$};
\node at (3.7, 2.15) {$s=\tfrac{1}{20}$};

\draw[fill=black] (1,0) circle (2pt);
\draw[fill=black] (1,1.5) circle (2pt);
\draw[fill=black] (0.5,0) circle (2pt);
\draw[fill=black] (0.5,1.875) circle (2pt);
\draw[fill=black] (-0.1,0) circle (2pt);
\draw[fill=black] (-0.1,1.995) circle (2pt);

\node at (1,-0.4) {1};
\node at (0.5,-0.4) {$\tfrac{1}{2}$};
\node at (-0.4,-0.4) {$-\tfrac{1}{10}$};
\end{tikzpicture}
\caption{Three secant lines determined by the points $(p, f(p))$ for $p= 1, \tfrac{1}{2},$ and $ -\tfrac{1}{10}.$ The lines get flatter the closer $p$ is to $0.$}
\end{figure}
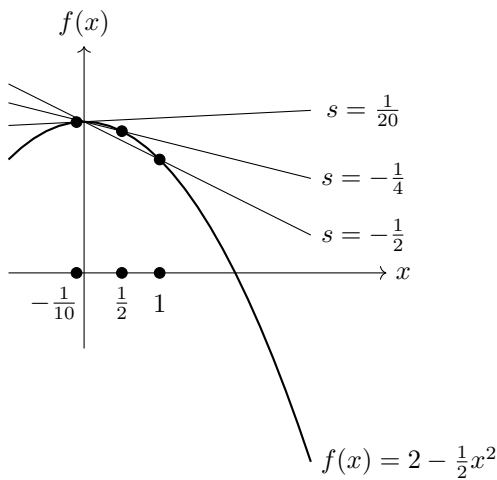

When $p=1,$ the other point is $(1,\tfrac{3}{2})$ and the slope is,
$$
s=\frac{f(1)-f(0)}{1-0}=\frac{\tfrac{3}{2}-2}{1-0} = -\tfrac{1}{2}.
$$
When $p=\tfrac{1}{2},$ the other point is $(\tfrac{1}{2},\tfrac{15}{8})$ and the slope is,
$$
s=\frac{f(\tfrac{1}{2})-f(0)}{\tfrac{1}{2}-0}=\frac{\tfrac{15}{8}-2}{\tfrac{1}{2}-0} = -\tfrac{1}{4}.
$$
When $p=-\tfrac{1}{10},$ the other point is $(-\tfrac{1}{10},\tfrac{399}{200})$ and the slope is,
$$
s=\frac{f(-\tfrac{1}{10})-f(0)}{-\tfrac{1}{10}-0}=\frac{\tfrac{399}{200}-2}{-\tfrac{1}{10}-0} = \tfrac{1}{20}.
$$
You can see what's happening: as you go from $p=1$ to $p=\tfrac{1}{2}$ to $p=-\tfrac{1}{10},$ the corresponding secant line gets flatter and flatter, which is the same as saying the slope of the secant line is getting closer to $0.$

\mpar
This is where the completely new idea comes in. Instead of doing the calculations for three secant lines, lets do it for infintely many secant lines determined by infinitely many points,
$$
 p_n, n\geq 1,\quad\text{meaning}\quad p_1, p_2, p_3, \cdots 
$$
supposing that the bigger the index $n$ is, the smaller the distance between $p_n$ and $0$ is. An example of such a sequence is,
$$
p_n =\frac{1}{n}, n\geq 1, \quad \text{meaning}\quad p_1=1,\, p_2=\tfrac{1}{2},\, p_3=\tfrac{1}{3}, \cdots
$$
and another one is,
$$
p_n =\frac{1}{10^n}, n\geq 1, \quad \text{meaning}\quad p_1=\tfrac{1}{10},\, p_2=\tfrac{1}{100},\, p_3=\tfrac{1}{1000}, \cdots
$$
There are lots of such sequences, once you start thinking about it.

\mpar
What does it mean to say ``$p_n$ gets closer to 0 as $n$ gets larger''? Turning this intuitive idea into solid mathematics is the completely new idea at the foundation of calculus, called \textit{limits}. It's a beautiful and subtle idea and it takes a few attempts to get the definition of limit right.

\mpar
To start, we say the the limit of $p_n$, as $n$ tends to infinity, equals $0$ if you can make $|p_n|$ is as small as you please just by taking $n$ large enough. If your idea of small is one one hundredth - meaning $10^{-2}$ - then you can find a number $N$ such that $|p_N|\leq 10^{-2}.$ But suppose your idea of small is now one one millionth, that is $10^{-7},$ what then? Well, then you can find another number $N^\prime$ bigger than $N$ such that $|p_{N^\prime}|\leq 10^{-7}.$ OK, so what if your idea of small is now an Angstrom ($10^{-10}$), Planck's length ($1.6\times 10^{-35}$), or even the reciprocal of a googolplex ($10^{-100}$)? The you can still find an even larger number $N^{\prime\prime}$ such that $|p_{N^{\prime\prime}}|\leq 10^{-100}.$ That's what it means to say you can make $|p_n|$ is as small as you please just by taking $n$ large enough. If you tell me your definition of small, call it $\epsilon,$ then I can find a big enough number, call it $N(\epsilon),$ such that, 
$$
|p_{N(\epsilon)}|\leq \epsilon,
$$
and I can do it for any $\epsilon>0.$ It's like a game. If I can always find $N(\epsilon)$ that does the trick no matter which positive number $\epsilon$ you choose, then I win the game and get to say that the limit of $p_n$ is $0.$

\mpar
Or do I? We want $p_N$ to be sufficiently small if $N$ is sufficiently large, but is that really enough? Don't we want $p_N$ to get small and then \textit{stay} small? That would mean,
$$
|p_{N(\epsilon)}|\leq \epsilon,\, |p_{N(\epsilon) +1}|\leq \epsilon, \, |p_{N(\epsilon)+2}|\leq \epsilon,\,\,\,\cdots \quad\text{etc.,}
$$
or more succinctly,
$$
|p_{n}|\leq \epsilon\quad\text{for all}\quad n\geq N(\epsilon).
$$
That turns out to be the correct definition of limit. You may have to think this through a couple of times before you get it but don't worry, everyone has to. It's one of those things that eventually become obvious once you've thought about it enough but doesn't seem obvious at all the first time you see it. 

\mpar
To summarize, we say the limit of $p_n$ is zero as $n$ tends to infinity and write,
$$
\lim_{n\to\infty}p_n=0
$$
provided for every $\epsilon>0$ there exists a natural number $N$ such that $|p_n|\leq\epsilon$ for all $n\geq N.$ 

\mpar
It's important to realize that the number $0$ didn't really play any role in the argument except as the proposed limit and neither did the actual value of the numbers $p_n,$ or the fact that we used $n$ to denote the index of the approximating numbers. The only thing that does matter is the relationship of the approximating numbers to the proposed limit. Therefore an expression like,
$$
\lim_{k\to\infty}q_k=\alpha
$$
means exactly what you think it does: for every $\epsilon>0$ there exists a natural number $N$ such that $|q_k-\alpha|\leq\epsilon$ for every $k\geq N.$  In the first case we say $p_n$ \textit{converges} to $0$ and in the second case we say $q_k$ \textit{converges} to $\alpha.$
\mpar
Let's go back to the example in Figures 16 and 17 and apply the idea of limits to the slope of the secant line. Using geometric arguments we convinced ourselves that,
$$
\text{if}\quad \lim_{n\to\infty}p_n=1 \quad\text{then}\quad \lim_{n\to\infty}\frac{f(p_n)-f(1)}{p_n-1} = 1
$$
and
$$
\text{if}\quad \lim_{n\to\infty}p_n=0 \quad\text{then} \quad\lim_{n\to\infty}\frac{f(p_n)-f(0)}{p_n-0} = 0.
$$

It's going to help us later on if we learn how to express limits in a more compact and succinct way. Notice the statements above are conditional: \textit{if} something happens \textit{then} something else happens. The condition on $p_n$ refers to any sequence whose limit is $0$ or $1,$ but the actual numerical values of the sequences are not really relevant. The only relevant thing is that they converge to $0$ or $1,$ as the case may be. Calling out a sequence of numbers whose specific values are actually irrelevant seems like overkill so we use the shorthand notation,
$$
f^\prime(p) = \lim_{x\to p}\frac{f(x)-f(p)}{x-p} 
$$
to mean the longer conditional statements where either $p=1$ and $f^\prime(1)=1$ or $p=0$ and $f^\prime(0)=0.$ We're just using the continuous variable $x$ to stand for any discrete sequence of numbers $p_n$ converging to $p.$

\mpar
Now that we've written the formula for the slope of the tangent line  this way, you realize that $p$ could be any number you please, not just $0$ or $1$ and, in fact, $f$ could be any function you please, not just $f(x)=2-x^2/2.$ 

\mpar
This is the definition of the derivative. The idea is to \textit{define} the slope of the best linear approximation to $f$ as the limit of the slope of a sequence of secant lines. More precisely, if the limit,
$$
f^\prime(p) = \lim_{x\to p}\frac{f(x)-f(p)}{x-p} 
$$
exists at a point $x=p$ then we declare that,
$$
y-f(p)=f^\prime(p)(x-p)
$$
is the equation of the line that is best approximates the graph of $f$ at $x=p.$ The derivative $f^\prime(p)$ is nothing else but the slope of the tangent line to the graph of $f$ at the point $(p, f(p)).$

\mpar
That's the basic idea of calculus in a nutshell: \textit{the slope of the tangent line is the limit of the slopes of a sequence of secant lines.}

\mpar
However, it's a entirely separate question whether you can calculate the derivative of important functions and turn those calculations to useful ends. The answer is, of course, yes you can but it takes a lot of work. In fact, it's the largest share of the work you do in a typical calculus class. 

\mpar
\subsection{Examples of Derivatives}
It turns out, for example, if $f(x)=1,\, x,\, x^2,\, x^3,$ etc. then,
$$
f^\prime(x) =0,\, 1,\, 2x,\, 3x^2,\,\text{etc.}
$$
It's also true that if $f(x) =\cos(x)$ then,
$$
f^\prime(x) = -\sin(x),
$$
and if $f(x)=\text{log}(x)$ is the natural logarithm then,
$$
f^\prime(x) =\frac{1}{x}.
$$
It's convenient to have alternative notation for the derivative of a function. This seems redundant at first but it turns out to be very useful. We write,
$$
\frac{df}{dx}(x)=\tfrac{d}{dx}f(x)= f^\prime(x).
$$
to mean the same thing and use whichever notation is most convenient in any given context. So, for example we can write,
$$
\frac{dx^2}{dx} =2x,\quad\tfrac{d}{dx}\cos(x)=-\sin(x),\quad \text{and}\quad \text{log}^\prime(x)=\frac{1}{x}.
$$
A more elaborate example is the derivative of  polynomial function,
\begin{align*}
(3x^2-5x+17)^\prime & =\tfrac{d}{dx}3x^2 +\tfrac{d}{dx}(-5x) +\tfrac{d}{dx}17\\
& = 3\tfrac{d}{dx}x^2-5\tfrac{d}{dx}x+17\tfrac{d}{dx} 1\\
& = 3\times 2x -5\times 1 +17\times 0\\
& = 6x-5.
\end{align*}
We slipped a few things under the rug here such as the  derivative of a sum of two functions is the sum of their derivatives and that the derivative of a constant multiple of a function is a constant multiple of the derivative. Again, these are the kinds of things you learn in first semester calculus, we're just making a note of them here.

\mpar
\subsection{Rates of Change}
We've used geometric ideas to formulate the definition of derivative in terms of slopes of approximating secant lines but there is a second, and very important, interpretation of the derivative as a \textit{rate of change}. Think of the function $f$ as describing the position of a car on a straight section of highway. At time $0$ the car starts to move down the highway and at time $t$ its position is $f(t).$ The ratio,
$$
r=\frac{f(t)-f(s)}{t-s}
$$
is the car's change of position over the time interval $[s,t], s<t,$ divided by the elapsed time. It has units of ``length divided by time'' so it is the average speed of the car over the interval $[s,t].$ Therefore the derivative,
$$
f^\prime
(t)=\lim_{s\to t}\frac{f(t)-f(s)}{t-s}
$$
\mpar
is the car's \textit{instantaneous} velocity at time $t$. This interpretation of the derivative as the instantaneous rate of change plays an important role in neural networks.

\mpar
It turns out to be useful to think of the instantaneous rate of change of a function as an approximation to the difference of the function at two nearby points. This just involves rewriting out the definition of derivative and interpreting the result in a new way. The definition of  $f^\prime(x)$ is just the statement that,
$$
\lim_{h\to 0}\left[ f^\prime(x)-\frac{f(x+h)-f(x)}{h}\right]=0
$$
which just means,
$$
f(x+h)\cong f(x)+f^\prime(x)\times h
$$
when $h$ is close to zero and the approximation gets better and better as $h$ gets closer to zero. This is just to say once again that the tangent line is the best linear approximation to the function but said in an effective way.
\mpar
\subsection{\textbf{The Chain Rule}}
The activations in a neural network are sums of compositions of functions and we need to be able to differentiate them and keep careful track of the terms. It turns out that the success of neural networks turns on a clever use of the \textit{chain rule} which is a formula for calculating the derivative of the composition of two functions. The chain rule states that the derivative of the composition $f\circ g$ is,
$$
\tfrac{d}{dx} f\circ g(x) = f^\prime(g(x))\,g^\prime(x),
$$
and it's proved this way,
\begin{align*}
\frac{d}{dx} f\circ g(x) = & \lim_{h\to 0} \frac{f(g(x+h))- f(g(x))}{h}\\
& = \lim_{h\to 0} \frac{f(g(x+h))- f(g(x))}{h}\times\frac{g(x+h)-g(x)}{g(x+h)-g(x)}\\
& =\lim_{h\to 0} \frac{f(g(x+h))- f(g(x))}{g(x+h)-g(x)}\times\frac{g(x+h)-g(x)}{h}\\
& \lim_{h\to 0} \frac{f(g(x+h))- f(g(x)))}{g(x+h)-g(x)}\times\lim_{h\to 0}\frac{g(x+h)-g(x)}{h}
\end{align*}
 By definition, the second limit equals $g^\prime(x).$ To handle the first limit note that,
 
 \begin{align*}
 \lim_{h\to 0} \frac{f(g(x+h))- f(g(x))}{g(x+h)-g(x)} & = \lim_{h\to 0}\frac{f(g(x)+g^\prime(x)\times h))-f(g(x))}{g^\prime(x)\times h}\\
 & = f^\prime(g(x)),
 \end{align*}
 which yields the proof. Actually, we've been a little cavalier in using the linear approximation $g(x+h)\cong g(x)+g^\prime(x)\times h,$ but a more careful argument would obscure the main idea without adding any real insight.
 \mpar
A simple example illustrates the rule. Suppose $f(x) = x^2$ and $g(x) = 3x+5.$ Then $f^\prime(x)=2x,\, g^\prime(x)=3,$ and so the chain rule says, 
 $$
\tfrac{d}{dx} f\circ g(x)=f^\prime(g(x))\,g^\prime(x)=2(3x+5)\times 3= 18x+30.
$$

On the other hand, we can multiply out $f\circ g(x)= (3x+5)^2 = 9x^2+30x +25$ and calculate the derivative as,
$$
\tfrac{d}{dx}f\circ g(x) = 9\tfrac{d}{dx} x^2+ 30 \tfrac{d}{dx} x + \frac{d}{dx} 25 = 9\times 2x+30 +0 = 18x+30,
$$
as it should be.

\mpar
\subsection{The Fundamental Theorem of Calculus}
Given a smooth positive function $f\colon\mathbb{R}\to\mathbb{R}$ and two points $a,b\in\mathbb{R}$ we can think of the set of points 
$$
B(f)_a^b=\{(t, f(t))\mid a\leq t\leq b\}
$$ 
as a box with a flat bottom and straight sides but a curvy top given by the graph of $f.$ The box has a definite area which we denote by the symbol,

$$
\text{area of $B(f)_a^b$}  =\int_a^b f(t)\, dt.
$$
If $f$ isn't positive but changes sign occasionally then we can think of $B(f)_a^b$ as a set of adjacent boxes, some with flat bottoms and curvy tops above the $x$-axis if $f>0$ and some with flat tops and curvy bottoms below the $x$-axis if $f<0$. The points $t$ where a curvy top box touches a curvy bottom box is just a point where $f$ changes sign from positive to negative. The \textit{signed area} of $B(f)_a^b$ is the sum of the areas of curvy top boxes minus the sum of the areas of the curvy bottom boxes.

\begin{figure}[h]
\begin{tikzpicture} [domain =-1:3.1]
\draw[->](-1,0) -- (4,0) node[right] {$x$};
\draw[->](0,-1) -- (0,3) node[above] {$f(x)$};
\draw[thick] plot(\x,{2-0.5*(\x)^2});
\fill[gray!30, opacity=0.5, domain = 1:2] plot(\x,{2-0.5*(\x)^2}) --(2,0) -- (1,0) -- cycle;
\fill[gray!30, opacity=0.5, domain = 2:3] plot(\x,{2-0.5*(\x)^2}) --(3,0) -- (2,0) -- cycle;

\draw[fill=black] (1,0) circle (2pt);
\draw[fill=black] (2,0) circle (2pt);
\draw[fill=black] (3,0) circle (2pt);
\draw[thin] (3,0)-- (3,-2.5);
\draw[thin] (1,0)-- (1,1.5);

\node at (1.34, 0.5) {\textbf{+}};
\node at (2.65, -0.65) {\textbf{--}};
\node at (1, -0.3) {1};
\node at (1.8, -0.3) {2};
\node at (3.2, -0.3) {3};
\end{tikzpicture}
\caption{The shaded region is $B(f)_1^3$, the space between the $x$-axis and the graph of $f(x)=2-x^2/2$ bounded by the vertical lines  $x=1$ and $x=3.$ The portion bounded by $x=1$ and $x=2$ has positive area and the portion bounded by $x=2$ and $x=3$ has negative area. The region $B(f)_1^3$ has signed area $\int_1^3f(x)dx=\int_1^2f(x)dx+\int_2^3f(x)dx.$}
\end{figure}
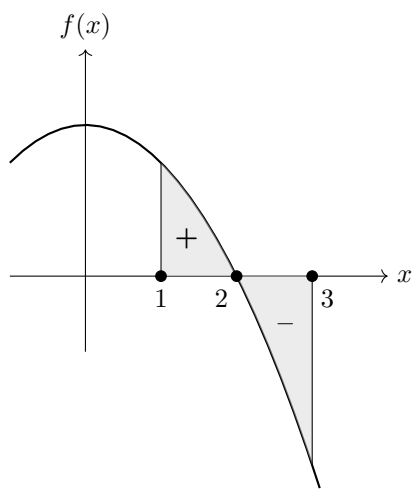
\mpar
The \textit{definite integral} of $f$ is the function $F$ you get by fixing the value of $a$ and letting the value of $b$ vary in the signed area of $B(f)_a^x,$ namely,
$$
F(x)=\int_a^xf(t)dt.
$$

If you think about it for a second you'll see that $F(x+h)-F(x)$ is the area of a box of width exactly $h$ and height approximately $f(x).$ This implies, 

$$
F(x+h)-F(x)\cong f(x)\times h \quad\text {or equivalently}\quad\frac{F(x+h)-F(x)}{h}\cong f(x)
$$

and therefore, 
$$
F^\prime(x)=\lim_{h\to 0}\frac{F(x+h)-F(x)}{h} =f(x).
$$

\mpar
This simple geometric picture is the basic idea behind \textit{integral calculus}. It takes some real effort to make solid mathematics out of the picture - that the typical subject matter of second semester calculus -  but the picture of integrals as signed areas under the graph of functions does faithfully represent the mathematics. We will not discuss integrals in this appendix, except to say that derivatives and integrals are inverse operations in much the same sense that squares and square roots or exponentials and logarithms and are inverses. That's the meaning of the \textit{fundamental theorem of calculus} which says if $f(x)$ is a smooth function and $f^\prime (x)$ is its derivative then,
$$
f(b) = f(a) +\int_a^bf^\prime (x) dx
$$
This formula says that if you know the derivative of a function at every point $x,$ and you know the value of $f$ at a single point $a,$ then you can figure out the value of $f$ and any other point $b$ via integration. It confirms that our definition the slope of the tangent line is derivative is the is the correct definition. Here's why. We've seen that if you know $f$ then you can calculate the derivative $f^\prime$. The fundamental theorem of calculus says the reverse: if you know the derivative $f^\prime$ then you can calculate $f$, at least if you know the value of $f$ at some point. This is the sense in which derivatives and integrals are inverse operations.

\mpar
The technicalities of defining integrals properly and then proving the fundamental theorem of calculus are a bit involved. In fact, the proof is often presented as the major highlight at the end of your first calculus course. But the basic underlying idea is pretty easy to grasp so it's worth knowing. 

\mpar
Instead of considering a function $f(x)$ of a continuous variable $x$ let's look at something simpler - a function $f(n)$ of a discrete variable $n$ where $n$ can be any whole number $n=0,\, \pm 1,\, \pm 2 \pm 3,$ etc. Let's say the \textit{discrete derivative} of $f$ at $n$ is,
$$
f^\prime (n) = \frac{f(n)- f(n-1)}{n-(n-1)} = f(n) - f(n-1).
$$
When we add the derivatives $f^\prime (n)$ together for all $a\leq n\leq b$ then we find,
\begin{align*}
\sum_{n=a}^bf^\prime(n) & = \sum_{n=a}^b \left(f(n)-f(n-1)\right)\\
& = f(b)-f(a-1).
\end{align*}

The sum is \textit{telescoping} in the sense that every term $f(n)$ appears twice - once with a plus sign and once with a minus sign - except for the endpoints $f(a-1)$ and $f(b).$ Thus,
$$
f(b) = f(a-1) + \sum_{n=a}^bf^\prime(n)
$$

If $g(n)$ is a function of the discrete variable $n$ let's call $\sum_{n=a}^bg(n)$ the \textit{discrete integral} of $g.$ Then what we've proved is the fundamental theorem of \textit{discrete calculus}: the discrete derivative and the discrete integral are inverse operations.

\mpar
\subsection{Directional Derivatives and the Gradient}
We'll see that calculating the derivative of multidimensional versions of the simple quadratic function $f(x)=x^2$ is a key step in training neural networks.  We know how to define the derivative of a function of one variable so it's natural to wonder if the ideas behind the definition extend to functions of many variables. 

\mpar
Let's start in $2D$ with a quadratic function,
\begin{align*}
g(x,y) & =2(x-1)^2 - 3(y+2)^2 -(x-y)^2\\
& = 2(x^2-2x+1) -3(y^2+4y+4) -(x^2-2xy+y^2)\\
& = (x^2+2xy-4y^2) -(4x+12y) -10.\\
\end{align*}

Here, the function is expressed as a sum of squares. Then we multiplied out the squares and grouped in parentheses the terms which are quadratic or linear in the  in the $x,y$ variables. The basic idea is two reduce from two dimensions to one dimension by freezing one of the variables so that it has a fixed value, and letting the other variable change. There are two ways to do this. First, suppose $y$ has a fixed value, say $y=-\tfrac{1}{2},$ regardless of the value of the $x$ variable. Plugging $y=-\tfrac{1}{2}$ into the formula for $g$ and multiplying everything out defines a new function,
\begin{align*}
f_1(x) & =g(x,-\tfrac{1}{2})\\
& = [x^2+2x\times(-\tfrac{1}{2}) -4\times(-\tfrac{1}{2})^2] -[4x+12\times(-\tfrac{1}{2})]-10\\
& =  x^2-5x+5,
\end{align*}
and its derivative is,
\begin{align*}
f_1^\prime(x) & = (x^2-5x +5)^\prime\\
& =\tfrac{d}{dx}x^2 -5\tfrac{d}{dx}x +5\tfrac{d}{dx} 1\\
& =2x-5\times 1 -5\times 0\\
& = 2x-5.
\end{align*}
Now suppose we freeze $x,$ say $x=2,$ and let $y$ vary. This defines a new function of one variable,
\begin{align*}
f_2(y) &= g(2,y)\\
& = [2^2+2\times 2y-4y^2] -[4\times 2+12y] -10\\
& = -4y^2-8y -14
\end{align*}
whose derivative is,
\begin{align*}
f_2^\prime(y) & = (-4y^2-8y-14)^\prime\\
& =-4\tfrac{d}{dy}y^2 -8\tfrac{d}{dy}y -14\tfrac{d}{dy} 1\\
& = -4\times 2y -8\times 1 -14\times 0\\
& = -8y-8.
\end{align*}

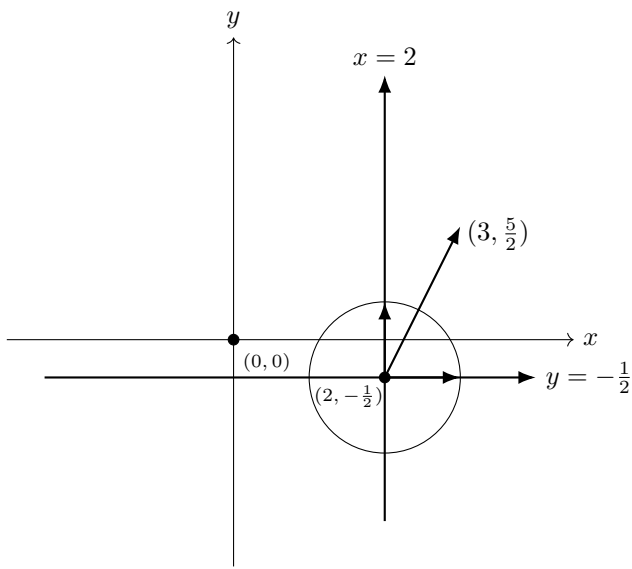
\begin{figure}
\begin{tikzpicture} 
\draw[thin, ->](-3,0) -- (4.5,0) node[right] {$x$};
\draw[thin, ->](0,-3) -- (0,4) node[above] {$y$};
\draw[thick, -Latex] (-2.5,-0.5) -- (4,-0.5) node[right] {$y=-\tfrac{1}{2}$};
\draw[thick, -Latex] (2,-2.4) -- (2,3.5) node[above] {$x=2$};
\draw[fill=black] (2, -0.5) circle (2pt);

\node at (1.525, -0.75) {\scriptsize $(2,-\tfrac{1}{2})$};
\draw[fill=black] (0,0) circle (2pt);
\node at (0.44, -0.3) {\scriptsize$(0,0)$};
\node at (3.5,1.4) {$(3,\tfrac{5}{2})$};
\draw[thick, -Latex] (2,-0.5) -- (2,0.5);
\draw[thick, -Latex] (2,-0.5) -- (3,1.5);
\draw[thick, -Latex] (2,-0.5) -- (3,-0.5);
\draw (2,-0.5) circle (1);

\end{tikzpicture}
\caption{The lines $x=2$ and $ y=-\tfrac{1}{2},$  intersect at the point  $(2,-\tfrac{1}{2}).$ We can calculate the derivative of $g$ at the point $(x,y)=(2,\tfrac{1}{2})$ in the unit directions $(0,1), (1,0)$  by restricting $g$ to these lines and differentiating the resulting functions of one variable $f_1(x)$ and $f_2(y)$. For reference these vectors are drawn as vectors on the unit circle centered at $(2,-\tfrac{1}{2})$ along with a third vector in the direction $(1,3).$ The coordinates of the endpoint of this vector, with respect to the origin $(0,0)$, are given by the vector sum $(2,-\tfrac{1}{2}) +(1,3) = (3,\tfrac{5}{2}).$}
\end{figure}

Notice that the lines $x=2$ and $y=-\tfrac{1}{2}$ intersect at the point $(2,-\tfrac{1}{2}),$ so that's a natural place to evaluate the derivatives of $f_1(x)$ and $f_2(y).$ See Figure 19. First we have,
\begin{align*}
\lim_{x\to 2}\frac{g(x,-\tfrac{1}{2})-g(2, -\tfrac{1}{2})}{x-2} & =\lim_{x\to 2}\frac{f_1(x)-f_1(2)}{x-2}\\
& = f_1^\prime(2)\\
& = 2\times 2 -5 \\
& =-1.
\end{align*}
On the other hand we have,
\begin{align*}
\lim_{y\to -\tfrac{1}{2}}\frac{g(2,y)-g(2,-\tfrac{1}{2})}{y-(-\tfrac{1}{2})} & =\lim_{y\to -\tfrac{1}{2}}\frac{f_2(y)-f_2(-\tfrac{1}{2})}{y-(-\tfrac{1}{2})}\\
& = f_2^\prime(-\tfrac{1}{2})\\
& = -8\times(-\tfrac{1}{2})-8 \\
& =-4.
\end{align*}
We call these the \textit{partial derivatives} of $g$ at the point $p=(1,-\tfrac{1}{2})$ and we denote them by,
$$
\frac{\partial}{\partial x}\, g(1,-\tfrac{1}{2})=-1\quad\text{and}\quad \frac{\partial}{\partial y}\, g(1,-\tfrac{1}{2})=-4.
$$
(These derivatives are called \textit{partial} because one variable is fixed and the other is changing).
\mpar\mpar
Let's step back from all these calculations and ask what do they really mean?

\mpar\mpar
Well first, if you recall the interpretation of the derivative as a rate of change, the calculations say the rate of change of $g$ at the point $p=(2,-\tfrac{1}{2})$ in the direction $u=(1,0)$ (that is, in the $x$ direction) is $-1$ while the rate of change of $g$ at the point $p=(2,-\tfrac{1}{2})$ in the direction $u=(0,1)$ (that is, in  the $y$ direction) is $-4.$ It seems pretty clear in geometric terms that if the slope of the tangent line is negative then it is ``pointing downwards''. Therefore, the graph of the function must also be pointing downwards which means it's decreasing. What we just learned is that $g$ is decreasing four times faster in the $y$ direction than it is in the $x$ direction at the reference point $(2, -\frac{1}{2})$ .
 
 \mpar\mpar
And second, observe that as the parameter $h$ varies in $\mathbb{R}$  from $-\infty$ to $+\infty$ the vectors,
 $$
 (2,-\tfrac{1}{2}) +h(1,0)\quad\text{and}\quad (2,-\tfrac{1}{2}) + h(0,1),
 $$
 trace out the horizontal and vertical lines in the plane that pass through our reference point $(2,-\tfrac{1}{2})$ when the parameter $h$ equals $0.$ Again, see Figure 19. Since $g((2,-\tfrac{1}{2}) +h(1,0)) = f_1(2+h, -\tfrac{1}{2}),$ we have,
 \mpar
 $$
 \lim_{h\to 0}\frac{g((2,-\tfrac{1}{2}) +h(1,0))-g((2,-\tfrac{1}{2}))}{h}=f_1^\prime(2) = -1,
 $$
 
 \mpar\mpar
 and since $g((2,-\tfrac{1}{2}) +h(0,1))=f_2(-\tfrac{1}{2}+h),$ we have,
 
 \mpar
 $$
 \lim_{h\to 0}\frac{g((2,-\tfrac{1}{2}) +h(0,1))-g((2,-\tfrac{1}{2})}{h}=f_2^\prime(-\tfrac{1}{2}) = -4.
 $$
 
\mpar\mpar
Now there are lots of directions in the plane other than the horizontal direction $(1,0)$ and the vertical direction $(0,1),$ right? If you think about it for a minute, you'll realize that you can calculate the rate of change of $g$ at the point $p=(2,-\tfrac{1}{2})$ in any direction you please, say direction $u=(a,b),$ by the formula,
 $$
 \lim_{h\to 0}\frac{g((1,-\tfrac{1}{2}) +h(a,b))-g((1,-\tfrac{1}{2}))}{h}.
 $$
The basic idea is to view,
$$ 
f_3(h) = g((2,-\tfrac{1}{2}) +h(a,b)) = g(2+ah, -\tfrac{1}{2} +bh)
$$
as a function of a single variable $h$ and then calculate $f_3^\prime(0)$, remembering that $a$ and $b$ are definite but as yet unspecified numbers.

\mpar
If you think about it for another minute, you'll also realize that the reference point $p=(2, -\tfrac{1}{2})$ plays a minor role in the calculation and that really you could repeat them at a general point $p=(x,y)$, meaning calculating the limit,
 $$
 \nabla_{(a,b)} g(x,y) = \lim_{h\to 0}\frac{g((x,y) + h(a,b))-g((x,y))}{h},
 $$
 or more succinctly, the limit,
 $$
 \nabla_u g(p) = \lim_{h\to 0}\frac{g(p + hu)-g(p)}{h}.
 $$

\mpar
The quantity $\nabla_ug(p)$ is  called the \textit{directional derivative of $g$ at the point $p$ in the direction $u$.} 

\mpar
This is a concise formula for a very powerful idea, so there's a lot to unpack. To grasp what's being said here, let's go back to the reference point $p=(2,-\tfrac{1}{2})$, pick a concrete direction like $u=(1,3),$  and calculate the corresponding directional derivative,
\begin{align*}
\nabla_{(1,3)}g(2,-\tfrac{1}{2}) & = \lim_{h\to 0}\frac{g\left((2,-\tfrac{1}{2}) + h(1,3)\right)-g(2,-\tfrac{1}{2})}{h}\\
& = \lim_{h\to 0}\frac{g(2+h,-\tfrac{1}{2}+3h) -g(2,-\tfrac{1}{2})}{h}.
\end{align*}

Let's simplify this a bit by introducing a new function,
$$
f_3(h) = g(2+h,-\tfrac{1}{2}+3h)
$$
and observing $f_3(0) = g(1,-\tfrac{1}{2})$ so that,
$$
\nabla_{(1,3)}g(1,-\tfrac{1}{2}) =\lim_{h\to 0}\frac{g(1+h,-\tfrac{1}{2}+3h) -g(1,-\tfrac{1}{2})}{h} = \lim_{h\to 0}\frac{f_3(h)-f_3(0)}{h}
$$
and therefore,
$$
\nabla_{(1,3)}\,g(2,-\tfrac{1}{2}) = f_3^\prime(0).
$$

\mpar
Apparently, we can calculate the directional derivative of $g$ by multiplying out the formula for $f_3$ and calculating derivative of $f_3^\prime(0).$ 

\mpar
This requires some work but it's not difficult work, just a bit lengthy. Referring back to the definition of the quadratic function,
$$
g(x,y)=(x^2+2xy-4y^2) -(4x+12y) -10,
$$
we have $x=2+h$ and $y=-\tfrac{1}{2}+3h,$ and so,
\begin{align*}
f_3(h) &= g(2+h, -\tfrac{1}{2}+3h)\\
& = [(2+h)^2 +2(2+h)(-\tfrac{1}{2}+3h) -4(-\tfrac{1}{2}+3h)^2] \\
& \quad\quad - [4(2+h) + 12(-\tfrac{1}{2}+3h)] -10\\
& = [(4+4h+h^2) +2(-1 +6h -\tfrac{1}{2}h+3h^2) - 4(\tfrac{1}{4} -3h+9h^2)]\\
& \quad\quad -[8+4h-6+36h] -10\\
& = -29h^2-13h -1.
\end{align*}

\mpar
So, when the dust settles, $f_3(h)$ is just a quadratic function of $h$ and we can calculate its derivative as,
\begin{align*}
f_3^\prime(h) & =-29\tfrac{d}{dh} h^2 -13\tfrac{d}{dh} h -{d}{dh} 1\\ 
& = -29\times 2h -13\times 1 -\tfrac{21}{2}\times 0\\
& = -29h - 13,
\end{align*}
hence,
$$
\nabla_{(1,3)}\,g(1,-\tfrac{1}{2}) = f_3^\prime(0) = -13.
$$

What comes next is a key insight in multivariable calculus. There's a simple relationship between the directional derivative $\nabla_{(1,3)}\,g(1,-\tfrac{1}{2}) = -13,$ the direction $u=(1,3),$ and the partial derivatives
$$
\frac{\partial}{\partial x}\, g(1,-\tfrac{1}{2})=-1\quad\text{and}\quad \frac{\partial}{\partial y}\, g(1,-\tfrac{1}{2})=-4.
$$
Can you see it? If you play around with the numbers long enough you'll realize that,
$$
-13 = 1\times (-1) + 3\times (-4)
$$
or, in other words,
\begin{align*}
\nabla_{(1,3)}\,g(2,-\tfrac{1}{2}) & = (1,3)\cdot\left(\frac{\partial}{\partial x}g(2,-\tfrac{1}{2}), \frac{\partial}{\partial y}g(2,-\tfrac{1}{2})\right)\\
& = (1,3)\cdot \nabla g(2,-\tfrac{1}{2}).
\end{align*}

\mpar
Here we introduced the notation $\nabla g(2,-\tfrac{1}{2})$ for the vector whose components are the two partial derivatives of $g$ at the point $(2,-\tfrac{1}{2}).$ 

\mpar
This dot product formula is not a coincidence. If you perform the same calculations at any point $p=(x,y)$ in any direction $u=(a,b)$
you always find that,
$$
\nabla_{(a,b)}\,g(x,y) = (a,b)\cdot\left(\frac{\partial}{\partial x}g(x,y), \frac{\partial}{\partial y}g(x,y)\right) = (a,b)\cdot\nabla g(x,y)
$$
This formula doesn't just hold for our particular quadratic function either, it holds for any function whose graph is a smooth surface. The vector $\nabla g(x,y)$ is called the \textit{gradient} of $g$ at the point $(x,y)$ so the dot product formula says,

\mpar
\textit{the directional derivative of a smooth function $g$ at a point $p$ in the direction $u$ is the dot product of $u$ and the gradient of $g.$}
\mpar
It's an amazingly powerful formula because it says you can calculate any directional derivative (and there are infinitely many of them!) by calculating just two, namely the two partial derivatives.

\mpar
The gradient of $g$ has a very important geometric meaning because it specifies the slope of a pair of tangent lines to the graph of $g$ at $p.$ Remember, the graph of $g$ is a smooth surface in $\mathbb{R}^3.$ The tangent lines to the restrictions, 
$$
f_1(h) = g\left((x,y)+h(1,0)\right)\quad\text{and}\quad f_2(h) = g\left((x,y)+h(0,1)\right)
$$
 are a pair of lines in $\mathbb{R}^3$ that intersect at the point $(x,y,g(x,y))$ and their slopes are the partial derivatives if $g$. Just like two points in a plane determine a line, two intersecting lines in space determine a plane. The plane determined by these two tangent lines is called the \textit{tangent plane} to the graph of $g$ at the point $p$. The geometric content of the dot product formula for directional derivatives is that the tangent line to the restriction,
 $$
 f_3(h)=g(p+hu)
 $$
 lies in the tangent plane to $g$ at $p$ and the slope of the tangent line is just $u\cdot\nabla g(p).$
 
 \mpar
Not surprisingly, the idea of directional derivatives, gradients and the dot product formula extends from $2D$ to $3D$ and all higher dimensions, as long as you make suitable adjustments. For example, in $3D$ there are three partial derivatives corresponding to the $x,y$ and $z$ directions so there are three partial derivatives, hence the gradient $\nabla g$ is a three dimensional vector. If you calculate the directional derivative,
$$
\nabla_{(a,b,c)}\,g(x,y,z) = \lim_{h\to 0}\frac{g\left((x,y,z)+h(a,b,c)\right) - g(x,y,z)}{h}
$$
you find that,
$$
\nabla_{(a,b,c)}\,g(x,y,z) = (a,b,c)\cdot\nabla g(x,y,z).
$$
Now, the graph of $g$ is a smooth three dimensional manifold in $\mathbb{R}^4$ and it has a three dimensional tangent space at $(x,y,z, g(x,y,z))$ spanned by the tangent lines to the restriction of $g$ to the lines $(x+h,y,z), (x, y+h, z)$ and $(x,y, z+h).$ The slopes of these lines are the components of the gradient and the dot product formula says that the tangent line to the restriction of $g$ to the line $(x,y,z) + h(a,b,c)$ lies in the tangent space and has slope $(a,b,c)\cdot\nabla g(x,y,z).$

\mpar
This all gets harder and harder to visualize as the dimension grows but that doesn't really matter. The dot product formula,
$$
\nabla_u\,g(p)=u\cdot\nabla g(p)
$$
still holds whether $p$ and $u$ have $4$ dimensions, $100$ dimensions, or $2,940,000$ dimensions. That's what really matters.

\mpar
\subsection{Gradient Descent}
Let's go back to $2D$ and think about the directional derivative of a function $g$ at a point $p$ as a function of the direction $u.$ Since any vector $u$ can be written as $||u||\times u/||u||$ let's restrict our attention to vectors of unit length. If $||u||=1$ then it's endpoint lies on the unit circle and the numerical value of $\nabla_u g(p)$ changes as the endpoint spins around the circle. 

\mpar
In which unit direction is the directional derivative as large as it can possibly be? In other words, for which unit vector $u$ is $\nabla_u g(p)$ as large as possible? That would be the direction in which $g$ is changing the fastest in the neighborhood of the point $p$ which we call the \textit{direction of steepest ascent.} So, the question becomes how do we find it? 

\mpar
This turns out to be easier than you'd think because of our formula for the angle between two vectors. Recall that if $u$ and $v$ are any two vectors in the plane then,
$$
u\cdot v = ||u||\,||v|| \cos\theta,
$$
where $\theta$ is the angle between them. Letting $v=\nabla g(p)$ we have,
\begin{align*}
\nabla_u\,g(p) & = u\cdot\nabla g(p)\\
& = ||u||\times ||\nabla g(p)||\cos\theta\\
& =||\nabla g(p)||\cos\theta,
\end{align*}
because $u$ has unit length. The caption of Figure 20 explains why $\cos\theta$ is as large as it can be when $\theta=0$ and also that $\cos(0)=1.$ This just means that the direction of steepest ascent is,
$$
u=\frac{\nabla g(p)}{|| \nabla g(p)||}
$$
and for that choice of $u,$
$$
\nabla_u\,g(p) = ||\nabla g(p)||.
$$

It's helpful to say in words what we just figured out in symbols: \textit{the gradient is the direction of steepest ascent}, or, more explicitly,

\mpar
\textit{the direction in which $g$ is changing the fastest in the vicinity of the point $p$ is $u=\nabla g(p)/||\nabla g(p)||$ and the rate at which $g$ is changing in the direction $u$ is $||\nabla g(p)||.$}
\mpar

\begin{figure} 
\begin{tikzpicture}
\draw (0,0) circle (2);
\draw[thick, ->] (0,-3) -- (0,3);
\draw[thick, ->] (-3,0) -- (3,0);

\draw[very thick, -Latex] (0,0) -- (1.5, 1.35);
\draw[very thick] (1.5,0) -- (1.5, 1.35);
\node at (0.7, 0.25) {$\theta$};

\draw[fill=black] (1.5, 0) circle (2pt);
\draw[fill=black] (0, 1.35) circle (2pt);

\node at (0.5, 1.35) {$\sin\theta$};
\node at (1.5, -0.3) {$\cos\theta$};
\node at (3, 1.35 ) {$u =(\cos\theta, \sin\theta)$};

\end{tikzpicture}
\caption{This is a repetition and simplification of Figure 12. The endpoint of a vector $u$ of length 1 lies on the unit circle. It's coordinates are $(\cos\theta, \sin\theta)$ where $\theta$ is the angle between $u$ and $(1,0).$ This is because the definition of cosine is \textit{length of adjacent side over length of hypotenuse} and the definition of sine is \textit{length of opposite side over length of hypotenuse.} As $\theta$ increases from $0$ to $2\pi$ the endpoint of $u$ spins counterclockwise around the circle. As $\theta$ increases from $0$ to $\pi,\,\cos\theta$ decreases from $+1$ to $-1$ and as $\theta$ increases from $\pi$ to $2\pi\,\cos\theta$ increases from $-1$ to $+1.$ Evidently, the maximum value of $\cos\theta$ is $+1$ and this happens when $\theta=0.$}
\end{figure}
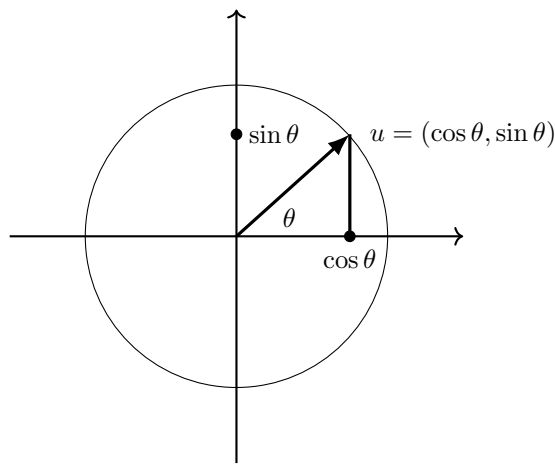

\mpar
To really grasp the meaning of this fact, let's do a little thought experiment. Let's think of the graph of $g$ as a landscape with hills and valleys scattered about and let's imagine we want to climb to the top of the nearest hill starting from the point $(p_0, g(p_0)).$ Which way should you go? 

\mpar
Well, one obvious strategy is to move in the direction of steepest ascent, namely $\nabla g(p)/||\nabla g(p)||,$ walk in that direction some distance then stop and reassess. Your new location is another point of the graph, say $(p_1, g(p_1)).$ The landscape may be different at your new location so the strategy would be to pick the direction of steepest ascent relative to your new location, namely $\nabla g(p_1)/||\nabla g(p_1)||,$ walk in that direction for a bit until you arrive at another new location, say $(p_2, g(p_2)),$ and reassess. You should do this for as many steps as you need to get to the top of the hill.  How do you know you're at the top of the hill? That's easy. At the top of a hill the gradient vanishes so the practical thing to do would be to stop when the gradient is pretty small.

\mpar
The image this generates in your mind is a zig-zag path of multiple steps that starts at $(p_0, g(p_0))$ and stops when the length of the gradient falls below some prescribed level. A conceptually equivalent procedure is the go in the direction proportional to the gradient, say  $\eta\nabla g(p_n),$ for a period of time rather than the direction of the gradient  $\nabla g(p)/||\nabla g(p)||$ for a certain distance. Moreover, he historical convention is to roll down the hill to find a valley, rather than climb the hill to find the top, so the end result is a set of zig-zig paths rolling down hill with velocity $-\eta\nabla g(p_n)$ until such time as the length of the gradient falls below some prescribed level. This method of finding the minimum of a function is called the \textit{method of steepest descent.}

\mpar
Steepest descent is an example of a class of optimization methods collectively called \textit{the greedy algorithm.} It's an algorithm because it presupposes you can't look up and just ``see'' the right step to take. It's called greedy because at every decision point you take the step that is best given the local information available to you - you work locally because you can't ``see'' the globally best step to take. 

\mpar
Although we discussed the method of steepest descent in $2D$ you can see that the method can be used in any dimension. It only depends of the mathematical fact that the negative of the gradient $-\nabla g(p)$ points in the direction in which $g$ is decreasing most rapidly. 

\mpar
There are many elaborate versions of steepest descent and using it in practice it's as much of an art than a science. For example, what value of $\eta$ should you choose? And what's a good prescribed level at which to stop? More importantly, if the graph of $g$ has many dimensions it becomes a heavy computational task to calculate the gradient if it has tens of thousand or even millions of coordinate directions. This may put restrictions on the total number of steps you can feasibly take on a given computational budget. Above all, if $g$ has a complicated graph with many hills and valleys and if your algorithm stops at the bottom of some valley, how do you know there's not a deeper valley somewhere else? There is a large amount of accumulated practical experience and some theoretical understanding of these questions but for the most part there are no truly definitive answers. It's an area of intense current interest in biology, chemistry, computer science, engineering, mathematics, and physics.

\mpar
Despite these questions, the method of steepest descent has an amazingly good track record of finding useful solutions to a large number of practical problems at a reasonable cost. A variant called \textit{stochastic gradient descent} is the method of choice in most neural networks which we discuss in the next section.

\mpar
\section{Probability}
To understand how neural networks are initialized and how they learn from data, we need to know something about the Gaussian, or normal, probability distribution. It's the familiar ``bell-shaped curve'' that describes many situations where a characteristic of a population, like a person's height or the noon-day temperature in Saskatoon in July, fluctuate around a typical value. 

\mpar
We also need to understand the Law of Large Numbers, sometimes called the law of averages. It quantifies the familiar fact that averages of the outcomes of many repetitions of the same  experiment tend to settle down to a common value. A simple way understand why this law holds is to calculate the mean and variance of the number of heads in many repeated coin flips. With a little extra work we can see that fluctuations of the average number of heads about the common value is described by the bell-shaped cure. This is known as the Central Limit Theorem.

\mpar
\subsection{Bernoulli Trials}
Because probability theory is less familiar to most people than calculus, we're going to include a few more calculations in this section than in the others. You're welcome to take these calculations on faith, but trying to work through them and seeing for yourself why they're correct will provide a deeper level of understanding that you may find satisfying and rewarding. 

\mpar
Imagine flipping a potentially biased coin. It can come up ``heads'' with probability $p$ or ``tails'' with probability $1-p.$ If the coin is fair then $p=\frac{1}{2}.$ This is the simplest probability distribution and it is the basic building block of many, many others. It's called the Bernoulli distribution, in honor of the Swiss mathematician Jacob Bernoulli (1655-1705).

\mpar
We say a \textit{random variable} $X$ has Bernoulli distribution provided $P(X=1)=p$ and $P(X=0)=1-p.$ This is convenient shorthand for the statement \textit{the probability a flipped coin comes up heads is $p$ and the probability it comes up tails is $1-p.$} The \textit{mean}, or average, or expected value of $X$ is,
$$
E(X)= 1\negthickspace \times \negthickspace p\,\, +\,\, 0\negthickspace\times\negthickspace (1-p) = p.
$$
It's worthwhile taking a minute to think about what this formula is saying. There are two possible outcomes in a coin flip experiment: heads or $1,$ which happens with probability $p,$ and tails or $0,$ which happens with probability $1-p.$ So, the formula just says the mean equals,
$$
\textit{possibility $\times$ probability} 
 $$
 summed over all the possibilities. This is the general rule for computing averages and it's a good image to keep it in mind.
 
 \mpar
The \textit{variance} of $X$ measures the extent to which $X$ deviates from its mean. The idea is that $X$ can be larger or smaller than $E(X),$ which means the difference $X-E(X)$ can be positive or negative. If we square the difference we get a new, non-negative random variable,
 $$
 Y=|X-E(X)|^2
$$
and its average is called the variance of $X,$ denoted $Var(X).$ Although it may seem complicated at first, it's actually easy to calculate the variance using the general rule above. There are two possibilities for $Y\negthickspace : $ either $Y=(1-p)^2,$ which happens if $X=1,$ or $Y=p^2,$ which happens if $X=0.$ Using a little bit of algebra, we find,
\begin{align*}
Var(X) & = E(Y)\\
& =\textit{sum of possibilities $\times$ probabilities}\\
& = (1-p)^2\negthickspace\times \negthinspace P(X=1) \,\, + \,\,p^2\negthickspace\times\negthickspace P(X=0)\\
& = (1-p)^2\negthickspace\times \negthinspace p \,\, + \,\,p^2\negthickspace\times\negthickspace (1-p)\\
& = (1-2p+p^2)p+ p^2(1-p)\\
& = (p-2p^2+p^3)+(p^2-p^3)\\
& = p-p^2\\
& =p(1-p).
\end{align*}
For any random variable, not just a Bernoulli random variable, it's handy to write $\mu = E(X)$ for the mean and $\sigma^2 = Var (X)$ for the variance.  It's natural to compare mean and variance because helps us understand if $X$ is tightly or loosely arranged about its mean. To make a proper comparison, though, we have make sure it's apples-to-apples. What does this mean?

\mpar
Suppose $X$ has units of measurement like length or mass. Then $\mu$ has units of length or mass but $\sigma^2$ has units of area or mass squared. This suggests looking at $\sigma,$ the square root of the variance, which is called the \textit{standard deviation} of $X.$ It has the same units as the mean and therefore the ratio $\mu /\sigma$ is a dimensionless number. 

\mpar
It's called the \textit{signal-to-noise ratio} because if $\mu /\sigma$ is large then $\sigma$ is small compared to $\mu$ and so $X$ must be tightly arranged about its mean. On the other hand, if $\mu /\sigma$ small then $\sigma$ is large compared to $\mu$ and therefore $X$ is loosely arranged about its mean. The upshot is this: a large signal-to-noise ratio means its easy to see the signal $X$ is sending you - $X$ is basically just $\mu$  plus a small error - while if the signal-to-noise ratio small then $X$ is swamped by the noise and the signal is hard to discern.

\mpar
Now imagine flipping a coin many times and keeping track of the number of heads. Suppose we flip it $N$ times, where $N$ could be a thousand, one hundred thousand, even a million or more, and let $X_1, X_2,\cdots, X_N$ denote the individual outcomes. This set of random variables is called \textit{Bernoulli trials} and $X_j$ is referred to as the outcome of the $j$-trial. Observe that the sum,
$$
S_N=X_1+X_2+\cdots + X_N
$$
is the total number of heads in $N$ trials. Why? Because $X_j=1$ if and only if the $j$-th flip came up heads and therefore $S_N$ is the exact number of times a head occurred in $N$ trials. Observe that $S_N$ is random and could be any number between 0 and $N.$

\mpar
We are assuming that the $X_j$'s have the same Bernoulli distribution - they're \textit{identically distributed} - so they have the same mean and variance. We're also assuming the flips are \textit{statistically independent} of one another, meaning the outcome of one flip has no effect on the outcome of any other flip. For example, it's obvious the  second flip can't affect the first flip because that flip has already happened. So we're basically just saying the outcome of a given flip does not affect the outcome of any subsequent flip. 

\mpar
To close this subsection let's look at the case $p=1/2$ where we're flipping a fair coin $n$ times. Assuming the flips are independent there are two equally likely outcomes in each of $n$ flips so there are $2^n$ equally likely outcomes altogether, where the probability, naturally, is $2^{-n}.$ It turns out to be quite useful to generalize this setup to equally likely outcomes for any finite set of possibilities.

\mpar
For example, a standard die has six possible outcomes: $1,2,3,4,5, 6.$ Assuming the die is fair each outcome is equally likely so they each have probability $1/6.$ If you roll it twice, each outcome such as  $(3,5), (1,2)$ or even $(5,3)$ is equally likely. Since there are $6\times 6=36$ possible outcomes they must all have probability $1/36.$ Now think of non-standard dice like sometimes occur in board games, like Scattergories which has a 20-sided die. Assuming the die is fair, each side is equally likely so each side each has probability $1/20.$ You get the picture.

\mpar
We say a \textit{uniform random variable} is an unbiased choice of one outcome among many. If there are $N$ possible outcomes then each one has probability $1/N.$ When we discuss mini-batch gradient descent we'll need to choose a permutation uniformly at random. If there are $n$ letters then there are $n!=n\times(n-1)\times(n-2)\times\,\dots\, 3\times 2\times 1$ permutations so each one has probability $1/n!.$

\mpar\mpar 
\subsection{Law of Large Numbers}
Let's calculate the mean and variance of $S_N$ under the assumption that $X_j,\, 1\leq j\leq N,$ are independent and identically distributed. It's easy to see the mean is,
\begin{align*}
E(S_N) & = E(X_1)+E(X_2)=\cdots +E(X_N)\\
& = p+p+\cdots +p\\
&=Np.
\end{align*}
The variance calculation is a bit more complicated so we break it down into two pieces. The easier piece is reducing the calculation of $Var(X)$ to the calculation of $E(S_N^2).$ The harder piece is calculating $E(S_N^2)$ which we do by multiplying out $(X_1+X_2+\cdots +X_N)^2$ using the binomial theorem then evaluating terms resulting terms $E(X_i X_j)$ for all possible choices of $i$ and $j$. 

\mpar
Let's begin with the first piece,
\begin{align*}
Var(S_N) & =E|S_N - E(S_N)|^2\\
& = E|S_N-Np|^2\\
& = E(S_N^2 -2Np\negthinspace\times\negthinspace S_N +(Np)^2)\\
& = E(S_N^2) - 2NpE(S_N) + N^2p^2\\
& =E(S_N^2) -2Np\negthinspace\times\negthinspace Np +N^2p^2\\
& = E(S_N^2)-N^2p^2.
\end{align*}
Before moving on to the second piece, let's remember the binomial theorem. It's a formula for the square of a sum, for example,
\begin{align*}
(a+b)^2 & = a^2+2ab + b^2\\
(a+b+c)^2 & = a^2+b^2 + c^2 +2ab + 2ac + 2bc
\end{align*}
and what it says in words is, 
$$
\textit{The square of a sum equals the sum of the squares plus the sum of the cross terms.}
$$
We're going to use the binomial theorem in the following form to calculate $ES_n^2,$
$$
(\,\sum_{i=1}^n x_i\,)^2 =\sum_{i=1}^nx_i^2 + \sum_{i\neq j}^nx_ix_j.
$$
The second term on the right contains cross terms, such as  $x_2 x_5$ and $x_5 x_2$ for example, so all cross terms appear exactly twice, as they should. We have,
\begin{align*}
E(S_n)^2 & = E (\,\sum_{i=1}^n X_i\,)^2\\
& = E(\,\sum_{i=1}^n X_i^2\,) + E(\,\sum_{i\neq j}^n X_i X_j\,)\\
& = \sum_{i=1}^nE(\, X_i^2\,) + \sum_{i\neq j}^n E(\,X_i X_j\,)\\
& = nE(\, X_1^2\,) + n(n-1) E(\, X_1 X_2\,).
\end{align*}

Why is the last step true? Well, the random variables are identically distributed so $E(\, X_i^2\,) = E(\, X_1^2\,)$ and $E(\,X_i X_j\,) = E(\,X_1 X_2\,).$ There are $n$ terms in the first sum and there are $n^2-n=n(n-1)$ terms in the second sum, because there are $n^2$ cross terms terms altogether and $n$ of them have the same index.

\mpar
Now, notice that $X_1^2 =1$ with probability $p$ and $X_1^2=0$ with probability $1-p$  so $E (\,X_1^ 2\,) =p.$ By the same logic, $X_1X_2=1$ if and only if both $X_1=1$ and $X_2=1$ and this happens with probability $p^2,$ since $X_1$ and $X_2$ are independent random variables. It follows that $X_1X_2=0$ with probability $1-p^2$ hence $E (\, X_1 X_2\,) = p^2.$ The end result is therefore, $E(\,S_N^2\,) = Np -N(N-1)p^2.$

\mpar
Putting this altogether we can calculate the variance of $S_N$ as
\begin{align*}
E(\ S_N - E(\,S_N\,)\,)^2 & = E(\,S_N^2\,) - N^2p^2\\
& = Np + N(N-1)p^2 -N^2p^2\\
& = N(p-p^2)\\
& = Np(1-p).
\end{align*}

If you recall that the variance of $X_1$ is $p(1-p)$ we've just shown that $Var(S_N) = N Var (X_1).$ This is true for independent Bernoulli random variables but it turns out to be true for any sequence of independent random variables or, in words,
$$
\textit{The variance of a sum of independent random variables is the sum of the variances.}
$$
The random variable $S_N$ the number of heads that occur in $N$ flips of a biased coin which could be any integer between $0$ and $N$. The fraction of heads in $N$ flips is $S_N/N$ and we just figured out its average, namely,
$$
E(\,\frac{S_N}{N}\,) = \frac{E(\,S_N\,)}{N} = \frac{Np}{N} =p = E(\,X_1\,).
$$
This formula is worth thinking about. It says,
$$
\textit{The average of the fraction of heads is the success probability.}
$$
What about a straight up comparison between $S_N/N$ and $p$, meaning how close is $S_N/N$ to $p,$ actually? One way to go is to look at the difference,
$$
\Delta_N=\frac{S_N}{N}-p,
$$
or better yet, the square of the difference,
$$
\Delta_N^2=|\frac{S_N}{N}-p|^2.
$$

Now $p$ is a fixed number but the fraction of heads is random so to put things on an equal footing, let's calculate $E(\,\Delta_N\,)$ and  $E(\,\Delta_N\,)^2.$ We've already seen that $E(\,S_N/N\,) = p,$ hence $E\,\Delta_N=0.$ So let's calculate the $E\Delta_N^2$ which is equal to the variance of $S_N/N.$ We have,

\begin{align*}
E\,\Delta_N^2 & = Var(\,\frac{S_N}{N}\,)\\
 & = E(\,|\frac{S_N}{N}-p|^2)\\ 
& = E(\,\frac{1}{N^2}|S_N-Np|^2\,)\\
& = \frac{1}{N^2}E(\,|S_N-E(\,S_N\,)|^2)\\
& =\frac{1}{N^2}\times Np(1-p)\\
& = \frac{p(1-p)}{N}.
\end{align*}

Stepping back from the calculations what we've found are formulas for mean and the variance of the fraction of heads in $N$ flips. Specifically,
$$
\mu_N= E(\,\frac{S_N}{N}\,) = p,
$$
and
$$
\sigma_N^2 = Var(\,\frac{S_N}{N}\,) = \frac{p(1-p)}{N}.
$$
It follows that the signal-to-noise ratio is,
$$
\frac{\mu_N}{\sigma_N} =\frac{p\,\sqrt{N}}{\sqrt{p(1-p)}},
$$
which gets larger and larger as the number of flips increases. Recall that a really large signal-to-noise ratio means that,
$$
\frac{S_N}{N} = p + \epsilon_N
$$ 
where $\epsilon _N$ is a negligible bit of mean-zero noise that is insignificant in size compared to $p.$ 

\mpar
Just how negligible and insignificant is $\epsilon_N$? With a bit of extra work involving something called \textit{Chebyshev's inequality} you can show that for any tolerance level $\lambda>0$ you choose, no matter how small, the probability that $|\epsilon_N|$ is greater than $\lambda$ is less than or equal to $p(1-p)/\lambda^2N$ which in symbols reads as,
$$
P(\,|\epsilon_N|>\lambda\,)\leq \frac{p(1-p)}{\lambda^2 N}.
$$
(The inequality is named in honor of the Russian mathematician Pafnuty Lvovich Chebyshev (1821-1894)). Since $p$ and $\lambda$ are fixed constants, but the number of flips $N$  is growing without bound, Chebyshev's inequality implies,
$$
\lim_{N\to\infty}P(\,|\frac{S_N}{N}-p|>\lambda\,)=0.
$$

This is the \textit{Law of Large Numbers} for Bernoulli random variables. It says in a very precise way that the long run fraction of heads is as close as you please to the success probability $p$ provided you flip enough coins.

\mpar\mpar
\subsection{Binomial Distribution}  We've already observed that the random variable $S_N$ has values that range from $0$ to $N.$ It's easy to calculate the probability that $S_N$ takes on one of the extreme values because $S_N=0$ if and only tails come up on every try while $S_N=N$ if an only if heads come up on every trial which implies,
$$
P(\,S_N=0\,) = (1-p)^N\quad\text{and}\quad P(\,S_N=N\,)= p^N.
$$
What about $P(\,S_N=1\,),$ the probability that there is exactly one head in $N$ flips? Or $P(\,S_N=2\,), \,P(\,S_N=3\,),$ or in general $P(\,S_N=k\,),$ for any $0\leq k\leq N$? How cab we calculate them? 

\mpar
Working out the $k=1$ case shows how to do it. If there was exactly one head in $N$ flips then there must have been $N-1$ tails and the lonely head must have occurred at some point: either the first flip, or the second, etc., or perhaps the $N^{th}$ and last flip. So, there are $N$ ways the head could have occurred and the for any given way, the probability of one heand and $N-1$ tails is $p(1-p)^{N-1},$ because the coin flips are statistically independent. Therefore,
$$
P(\, S_N=1\,) = Np(1-p)^{N-1}.
$$

So far so good. If there were two heads in $N$ flips they had to have occurred at two distinct times and any pair of times is one of the ways it could have happened. So, we just need to count the number of such pairs. Evidently, there are $N$ possible times for the first flip and $N-1$ times for the second flip so it looks like there are $N(N-1)$ ways to get two heads in $N$ flips. 

\mpar
But if you think about it for a second, you realize that can't be quite right. If $3$ and $7$ happened to be the two times you chose, that pair itself can happen in two ways: first 3 then 7 or first 7 then 3. So $N(N-1)$ counts every pair of times twice, meaning half that amount, $n(N-1)/2,$ is the correct count. By statistical independence, the probability of getting 2 heads and $N-2$ tails is $p^2(1-p)^{N-2}$ and therefore,
$$
P(\,S_N=2\,) =\frac{N(N-1)}{2}p^2(1-p)^{N-2}.
$$
So, the pattern is clear: calculate how many ways $k$ heads can occur in $N$ flips and multiply that number by $p^k(1-p)^{N-k}.$ You can pick $k$ times sequentially in 
$$
N(N-1)(N-2)\cdots (N-(k-1)) = N(N-1)\cdots (N-k+1)
$$
ways (you might have to think about that for a moment), but that overcounts things because we don't care about the order in which the times are chosen. There are,
$$
k!=k\times (k-1)\times (k-2)\cdots 3\times 2\times 1
$$
ways in which you can sequentially choose a given set of $k$ times so the correct number of ways must be $N(N-1)\cdots (N-k+1)/k!$

\mpar
There's a neat trick for cleaning up this expression. If $n$ is any positive number then $n!$ (pronounced ``$n$ factorial'') is the product of all numbers less than or equal to $n,$ 
$$
n!=n\times (n-1)\times (n-2)\cdots 3\times 2\times 1.
$$
With a little effort you can see that,
$$
N(N-1)\cdots (N-k+1) =\frac{N!}{(N-k)!}
$$
and therefore the number of ways you can choose $k$ times from a set of $N$ times is

$$
\binom{N}{k} = \frac{N!}{k!(N-k)!}.
$$
Putting it altogether we have,
$$
P(\, S_N=k\,) = \binom{N}{k}p^k(1-p)^{N-k}.
$$
This is called the \textit{binomial distribution}. It is the probability distribution of the number of successes in $N$ Bernoulli trials. It gets it's name from the binomial theorem which states,
$$
(a+b)^N = \sum_{k=0}^N\binom{N}{k}a^kb^{N-k}.
$$
The letters $a$ and $b$ are variables and it's a clever choice to pick $a=p$ and $b=(1-p)$ because then $a+b=1.$ In this case the binomial theorem says,
\begin{align*}
1 & =\sum_{k=0}^N\binom{N}{k}p^k(1-p)^{N-k}\\
& = \sum_{k=0}^N P(\, S_N=k\,),
\end{align*}
so all the probabilities add up to 1, just as they should.

\mpar
The binomial theorem is a powerful formula; here are just two impressive applications. We know that $E(\,S_N\,) =Np$ and $Var(\, S_N\,) = E(\,(S_N-Np))^2 =Np(1-p),$ right? If you recall that averages are computed by summing
$$
\textit{possibilities $\times$ probabilities}
$$
over all the possibilities we have,
$$
\sum_{k=0}^Nk\binom{N}{k}p^k(1-p)^{N-k} = Np,
$$
and
$$
\sum_{k=0}^N(k-Np)^2\binom{N}{k}p^k(1-p)^{N-k} = Np(1-p).
$$
Aren't those formulas amazing? Who knew such complicated sums had such simple answers.

\mpar\mpar
\subsection{Central Limit Theorem}
The Law of Large Numbers says that $S_N/N$ is very close to $p$ when $N$ is large; in fact it amounts to the statement that for any $\epsilon >0,$ no matter how small,
$$
\lim_{N\to\infty} P(\, 0\leq |\frac{S_N}{N} -p| \leq \epsilon\,) = 1.
$$

It turns out to be a good idea to ask how fast $S_N/N$ is approaching $p.$ It may seem counterintuitive at first but one way to find out is to multiply $|S_N/N-p|$ by a power of $N,$ say $\sqrt{N},$ square it and take the average. 

\mpar
You see, it's like a contest of strength: $|S_N/N-p|$ is shrinking to 0 but $\sqrt{N}$ is increasing to infinity so what's happening to their product? If the product still goes to zero then $|S_N/N-p|$ wins the contest and if the product blows up to infinity then $\sqrt{N}$ wins the contest. 

\mpar
But if the product does neither - not growing or shrinking - then it's a tie, and the two factors have the same strength. To find out which is which let's do the calculation,
\begin{align*}
E\left(\,\sqrt{N}\times |\frac{S_N}{N}-p|\,\right)^2 & = (\sqrt{N})^2 \times E(\,|\frac{S_N}{N}-p|^2\,)\\
& =N\times \frac{p(1-p)}{N}\\
& = p(1-p).
\end{align*}

So, it's a tie which means,
$$
\sqrt{N} (\,\frac{S_N}{N}-p\,) = \frac{S_N-Np}{\sqrt{N}}
$$
is a random variable with mean zero and variance $p(1-p)$ or, equivalently,
$$
Z_N = \frac{S_N-Np}{\sqrt{N p((1-p)}}
$$
is  random variable such that $E Z_N=0$ and $E Z_N^2 =1.$ Notice that we can rewrite this as,
$$
Z_N=\frac{S_N - E(\, S_N\,)}{\sqrt{Var (\, S_N\,)}}.
$$

\mpar
In words what we've discovered is: \textit{if you normalize $S_N$ by subtracting the mean and dividing by the square root of the variance, then the normalized random variables $Z_N$ have mean zero and variance 1. }

This is a very delicate balancing act. It brings the long run fraction of heads into sharp focus: $S_N$ is centered at its \textit{mean}, $Np,$ and distributed around the mean on a length scale determined by it's \textit{standard deviation}, $\sqrt{Np(1-p)}.$ \mpar
The \textit{Central Limit Theorem}, first proved by the English mathematician Abraham de Moivre (1667 - 1754), makes this observation precise. The theorem says that for any numbers $a <b,$
\begin{align*}
\lim_{N\to\infty} P(\,a&\leq \frac{S_N-Np}{\sqrt{N p(1-p)}} \leq b\,)\\
& = \lim_{N\to\infty} P(\,Np +a\sqrt{Np(1-p)}\leq S_N \leq Np +b\sqrt{Np(1-p)}\,)\\
& =\tfrac{1}{\sqrt{2\pi}}\int_a^b e^{-\tfrac{x^2}{2}}dx.
\end{align*}
In his proof, published in 1738, de Moivre demonstrated an early version of \textit{Stirling's formula} (James Stirling, Scottish mathematician (1692 - 1770)),
$$
\lim_{N\to\infty} \frac{N!}{N^{N -\tfrac{1}{2}}e^{-N}} =\sqrt{2\pi}
$$ and used it to approximate binomial probabilities as follows,
$$
\binom{N}{k}p^k(1-p)^{N-k} \,\cong\, \frac{1}{\sqrt{2\pi p(1-p)}}\times e^{-\tfrac{(k-Np)^2}{2Np(1-p)}}.
$$
With this approximation in hand, de Moivre made the substitution,
$$
x^2=\frac{(k-Np)^2}{Np(1-p)}
$$
and showed that an appropriate sum of binomial probabilities on the left hand side was approximately equal to an integral of the function $e^{-x^2/2}$ on the right hand side, and that the approximation became \textit{exact} in the limit $N\to\infty$. It was a mathematical \textit{tour de force} executed at a time when calculus was still a young subject.

\mpar\mpar
\subsection{Normal Distribution}
We say a random variable $Z$ has a \textit{standard normal distribution} provided for every $a < b,$
$$
P( a\leq Z\leq b\,) = \tfrac{1}{\sqrt{2\pi}}\int_a^b e^{-\tfrac{x^2}{2}}dx.
$$

This probability distribution has a long history and its roots are in empirical observation, for example of crop yields, or the heights or weights of human populations, or measurement errors in cartography, or observations of heads in coin flipping, for that matter. Such measurements tend to cluster around an average value and drop of symmetrically on either side following a characteristic ``bell-shaped curve''.  Data plots from distinct investigations across an enormous range of disciplines tended to have the same shape to the point that in 1894 the English statistician Karl Pearson (1857 - 1936) wrote of an ``error curve'' that he called the ``normal curve''. The name stuck.

\mpar
Here are five very important observations about the standard normal distribution and related distributions.  

\mpar
(1) First of all,
$$
\int_{-\infty}^\infty e^{-\tfrac{x^2}{2}}dx = \sqrt{2\pi}
$$ 
although it takes a bit of effort to prove. The reason we mention it is to verify that,
$$
\tfrac{1}{\sqrt{2\pi}}\int_{-\infty}^\infty e^{-\tfrac{x^2}{2}}dx =1,
$$
meaning that the probabilities really do add up to 1.  

\mpar
(2) The mean value of $Z$ is zero and its variance is one - that's why it's called the \textit{standard} normal distribution. In fact, 
$$
E(\,Z\,) = \tfrac{1}{\sqrt{2\pi}}\int_{-\infty}^\infty x e^{-\tfrac{x^2}{2}}dx =0,
$$
and,
$$
Var(\,Z\,) = \tfrac{1}{\sqrt{2\pi}}\int_{-\infty}^\infty x^2 e^{-\tfrac{x^2}{2}}dx =1.
$$
Once again, these formulas take some effort to prove (actually the first is easy if you know the trick) but they are very reminiscent of the formulas we already proved for the binomial distribution.

\mpar
(3) If $Z$ is standard normal then its possible values range across the whole real line. But possible values near 0 are much more likely than others even a modest distance away. For example, if  $x=\pm\sqrt{2},\, \pm\sqrt{4},\, \pm\sqrt{6},$ etc., then $e^{-x^2/2} = e^{-1},\, e^{-2},\,e^{-3},$ etc. This sequence of numbers is decreasing exponentially fast. So the function $e^{-x^2/2}$ decreases rapidly in $x$ hence the integral,
$$
P(\,Z\geq a\,) =\int_a^\infty e^{-\tfrac{x^2}{2}} dx 
$$
decreases rapidly in $a.$ Therefore, while $Z$ can in principle take on any real value, that value is very unlikely to be large.

\mpar
(4) A simple transformation of $Z$ changes the mean and variance in predictable ways. If $W=\mu +\sigma Z$ then it's easy to see that $E(\,W\,) = \mu$ and $Var(\,Z\,)=\sigma^2.$ Using this transformation, we can work out the probability distribution of $W$ because,
\begin{align*}
P(\,a\leq W\leq b\,) & = P(\,a\leq \mu +\sigma Z\leq b\,)\\
&  = P(\,\tfrac{a-\mu}{\sigma}\leq Z \leq \tfrac{b-\mu}{\sigma}\,)\\
& = \tfrac{1}{\sqrt{2\pi}}\int_{\tfrac{a-\mu}{\sigma}}^{\tfrac{b-\mu}{\sigma}} e^{-\tfrac{x^2}{2}} dx
\end{align*}
Using some standard calculus ideas to transform this integral we find,
$$
P(\,a\leq W\leq b\,)  = \tfrac{1}{\sqrt{2\pi\sigma}}\int_a^be^{-\tfrac{(x-\mu)^2}{2\sigma^2}} dx.
$$
The function,
$$
n_{\mu,\sigma}(x) = \tfrac{1}{\sqrt{2\pi\sigma}}\,e^{-\tfrac{(x-\mu)^2}{2\sigma^2}}
$$
is called the \textit{normal probability density} with mean $\mu$ and variance $\sigma$. The result of the calculation above is that,
$$
n_{\mu, \sigma}(x)=\tfrac{1}{\sqrt{\sigma}}\,n_{0,1}((x-\mu)/\sigma).
$$

\mpar
(5) What does statistical independence mean for two normal random variables where, say, $X$ is normal with mean $\mu$ and standard deviation $\sigma$ and $Y$ is normal with mean $\nu$ and standard deviation $\tau$? Well, it means for any numbers $a<b$ and $c<d,$

\begin{align*}
P(\, a\leq X\leq b\,\,\, \text{and} \,\,\,c\leq Y\leq d\,) & = P(\, a\leq X\leq b\,)\times P(\, c\leq Y\leq d\,)\\
& = \int_a^b\int_c^d n_{\mu,\sigma}(x)n_{\nu,\tau}(y) dxdy.
\end{align*}

In other words, $X$ and $Y$ are independent provided the probability that $X$ is in some interval and \textit{at the same time} $Y$ is in some other interval is just the product of the respective probabilities.

\mpar
Once you understand the meaning of independence for two random variables, it's easy to see what it means for three, or three thousand, or three million random variables to be independent. It means the probability that the inequalities $a_k\leq X_k\leq b_k$ for $1\leq k\leq N$  hold at the same time equals the product of the probabilities they hold individually. For normal random variable you can write it out in a formula like we did in the $N=2$ case above, but it's lengthy and just says in symbols what has already been said in words. Independence just means non-interference: knowing something about one random variable tells you nothing at all about another random variable.

\mpar
We're going to use normal random variable to initialize the weights and biases of our neural network prior to training it on data. It's important to be as non-judgmental as possible at the initial step because we want the network to learn from the data not from our choice of initial conditions. In particular, we don't want the initial weights and biases to know about each other or be influenced by one another so it makes sense that they are statistically independent, just as our coin flips were.

\mpar\mpar
\subsection{Stochastic Gradient Descent}
We've learned that gradient descent is an effective way, in principle, to find the minimum of a function of many variables. But a weakness of the method appears when there are a lot of variables and therefore a lot of partial derivatives to compute. Imagine it costs something to compute a partial derivative - either a cost in time or a cost in money - and you're on a computational budget. Of what practical use is the method of gradient descent if it costs you more resources to implement than you have available?

\mpar
This focuses attention on less expensive ways to approximate the direction of steepest descent in the hope that it still provides a good enough approximation of the minimum of the function. This is where probability enters the picture and it enters in a somewhat surprising way. This is where the adjective \textit{stochastic} comes from - it's a fancy synonym for random.

\mpar
In the case of training a neural network, the function we're trying to minimize - the total mean-square difference between the training data set and the computed answer - has a special form,
$$
\Delta=\Delta(W,B) =\frac{1}{N}\sum_{j=1}^N || \hat{y}(x_j,W,B)-y(x_j)||^2
$$ 
that we exploit in a random way. First of all, it is a function of the weights and biases of the network so these are the coordinates of the function to be minimized. Each data point in the training set adds a term 
$$
\delta^2(x_j) = || \hat{y}(x_j)-y(x_j)||^2
$$
to the total error and each term has it's own minimum. Since the total mean squared error is a sum of functions each with its own minimum, finding the values of the weights and biases that minimizes the total amounts to finding the right compromise among the terms, because you can't minimize them all simultaneously.

\mpar
The idea of \textit{stochastic gradient descent} is to randomly choose one data point in the training set, calculate the gradient of its error function, move the weights and biases in that direction at a certain speed (called the \textit{learning rate}), then stop. Next, randomly and independently pick another data point in the training set, calculate the gradient of the corresponding error function, move the weights and biases in that direction at a certain speed, then stop. Repeat as needed. (To make this work you need a sequence of independent uniform random variables on $\{1, 2,\,\dots,\, N\}).$ The advantage of this method is that it's computationally inexpensive but it's disadvantage is that progress towards the minimum is much slower that true gradient descent.

\mpar
A variant of stochastic gradient descent, called \textit{mini-batch descent}, is to pick a batch of data points in the training set, rather than just one, calculate the gradients of the corresponding terms, add them together, move the weights and biases in that direction, etc. To make this work you need a way of choosing a subset of size $R$ from $\{1, 2,\,\dots,\, N\}.$ (There are number of ways to do this but practical mini-batch relies on a sequence of independent uniform random variables on the set of permutations of $\{1, 2,\,\dots,\, N\}).$ If the batch size is not too big, the computational cost doesn't increase too much, hence the name mini-batch. What you gain for the additional computational cost is that the mini-batch gradient is closer to the true gradient than it is with stochastic gradient descent so the likelihood of convergence to the minimum is improved.

\mpar
In the Goldilocks analogy, gradient descent is like Papa Bear's bed - too hard; and stochastic gradient descent is like Mama Bear's bed - too soft; but mini-batch descent is like Baby Bear's bed - just right!

\end{document}